\documentclass{article}
\usepackage[utf8]{inputenc}
\usepackage[a4paper, total={6in, 8in}]{geometry}
\pdfoutput=1

\usepackage[pagebackref]{hyperref}
\usepackage{amsmath}
\usepackage{amssymb}
\usepackage{arydshln}
\usepackage{mathtools}
\usepackage{enumerate,enumitem}
\usepackage{amsfonts,mathrsfs}
\usepackage{subcaption}
\DeclareMathAlphabet{\mathcal}{OMS}{zplm}{m}{n}
\usepackage{amsthm,thmtools}
\makeatletter
\def\th@plain{%
  \thm@notefont{}
  \itshape 
}
\def\th@definition{%
  \thm@notefont{}
  \normalfont 
}
\makeatother
\usepackage{color,soul}
\usepackage{algorithm, algpseudocode, algorithmicx}
\usepackage{comment}

\newcommand*{\email}[1]{\href{mailto:#1}{\nolinkurl{#1}} } 

\usepackage{tikz}
\usetikzlibrary {positioning}
\usepackage{tikz-cd}
\usetikzlibrary{matrix, calc, arrows}

\newif\iftruncated
\newif\ifrelease
\newif\ifusedvrep

\truncatedtrue
\releasetrue
\usedvreptrue

\newcommand{\norel}[1]{\ifrelease \else #1 \fi}

\newcommand{\notrunc}[1]{\iftruncated \else #1 \fi}
\newcommand{\usedv}[1]{\ifusedvrep #1 \fi}
\newcommand{\nousedv}[1]{\ifusedvrep \else #1 \fi}

\newcommand{\ethan}[1]{\norel{{\color{blue} \textbf{Ethan:} #1}}}
\newcommand{\nithin}[1]{\norel{{\color{red} \textbf{Nithin:} #1}}}
\newcommand{\shivcheck}{\norel{{\color{magenta} \checkmark}}}
\newcommand{\ethancheck}{\norel{{\color{blue} \checkmark}}}
\newcommand{\ec}{\ethancheck}
\newcommand{\nithincheck}{\norel{{\color{red} \checkmark}}}

\newcommand{\real}{\mathbb{R}}
\DeclareMathOperator{\rank}{rank}
\newcommand{\twobytwo}[4]{\begin{bmatrix} #1 & #2 \\ #3 & #4 \end{bmatrix}}
\newcommand{\twobyone}[2]{\begin{bmatrix} #1 \\ #2 \end{bmatrix}}
\newcommand{\onebytwo}[2]{\begin{bmatrix} #1 & #2 \end{bmatrix}}
\renewcommand{\bar}{\overline}
\DeclareMathOperator{\diag}{diag}

\newcommand{\Ps}[1]{P_{#1}}
\newcommand{\Ptf}[2]{\Ps{#1,#2}}
\newcommand{\Qs}[1]{Q^T_{#1}}
\newcommand{\Qft}[2]{\Qs{#1,#2}}

\newcommand{\comp}[1]{\bar{#1}}
\newcommand{\intersect}[2]{#1#2}

\newcommand{\graph}{\mathbb{G}}
\newcommand{\subgraph}{\mathbb{H}}
\newcommand{\complementgraph}{\overline{\subgraph}}
\newcommand{\nodeset}{\mathbb{V}}
\newcommand{\edgeset}{\mathbb{E}}
\DeclareMathOperator{\nnz}{nnz}

\newcommand{\diagmat}[1]{\mathbf{#1}}

\newcommand{\dirpath}{\mathbb{P}}

\definecolor{myblue}{RGB}{220,220,255}
\definecolor{myred}{RGB}{255,220,220}

\definecolor{ggreen}{RGB}{59,179,0}
\definecolor{hred}{RGB}{255,25,25}
\definecolor{xblue}{RGB}{0,85,255}

\newcommand{\rg}{r^{\rm g}}
\newcommand{\rh}{r^{\rm h}}
\newcommand{\rghat}{\hat{r}^{\rm g}}
\newcommand{\rhhat}{\hat{r}^{\rm h}}
\DeclareMathOperator{\CB}{\mathcal{CB}}

\newcommand{\ropt}{r_{\rm opt}}

\newcommand{\PAbar}{\Ps{\comp{A}}}
\newcommand{\PBbar}{\Ps{\comp{B}}}
\newcommand{\PAB}{\Ps{\intersect{A}{B}}}
\newcommand{\QBbar}{\Qs{\comp{B}}}
\newcommand{\QCbar}{\Qs{\comp{C}}}
\newcommand{\QBC}{\Qs{\intersect{B}{C}}}

\newcommand{\QABB}{\Qft{\intersect{A}{B}}{B}}
\newcommand{\QABA}{\Qft{\intersect{A}{B}}{A}}
\newcommand{\QBbarB}{\Qft{\comp{B}}{B}}
\newcommand{\QAbarA}{\Qft{\comp{A}}{A}}
\newcommand{\PBBC}{\Ptf{B}{\intersect{B}{C}}}
\newcommand{\PCBC}{\Ptf{C}{\intersect{B}{C}}}
\newcommand{\PCCbar}{\Ptf{C}{\comp{C}}}
\newcommand{\PBBbar}{\Ptf{B}{\comp{B}}}

\newcommand{\RBbarAB}{R_{\comp{B},\intersect{A}{B}}}
\newcommand{\RBbarBbar}{R_{\comp{B},\comp{B}}}
\newcommand{\RBCAB}{R_{\intersect{B}{C},\intersect{A}{B}}}
\newcommand{\RBCBbar}{R_{\intersect{B}{C},\comp{B}}}

\newcommand{\FAbar}{F_{\comp{A}}}
\newcommand{\FCbar}{F^T_{\comp{C}}}

\newcommand{\BG}{B\!\! G}
\newcommand{\CH}{C\!\! H}

\newcommand{\PAbarbullet}{\Ps{\comp{A}\bullet}}
\newcommand{\PAG}{\Ps{\intersect{A}{\BG}}}
\newcommand{\PGbar}{\Ps{\comp{\BG}}}
\newcommand{\QGbar}{\Qs{\comp{\BG}}}
\newcommand{\QGH}{\Qs{\intersect{\BG}{\CH}}}
\newcommand{\QHbar}{\Qs{\comp{\CH}}}

\newcommand{\QAGG}{\Qft{\intersect{A}{\BG}}{\BG}}
\newcommand{\QAGA}{\Qft{\intersect{A}{\BG}}{A}}
\newcommand{\QGbarG}{\Qft{\comp{\BG}}{\BG}}
\newcommand{\QAbarbulletA}{\Qft{\comp{A}\bullet}{A}}
\newcommand{\PGGH}{\Ptf{\BG}{\intersect{\BG}{\CH}}}
\newcommand{\PHGH}{\Ptf{\CH}{\intersect{\BG}{\CH}}}
\newcommand{\PHHbar}{\Ptf{\CH}{\comp{\CH}}}
\newcommand{\PGGbar}{\Ptf{\BG}{\comp{\BG}}}

\newcommand{\QBCone}{\Qs{\intersect{B}{C}1}}
\newcommand{\QBCtwo}{\Qs{\intersect{B}{C}2}}
\newcommand{\QGHone}{\Qs{\intersect{G}{H}1}}
\newcommand{\QGHtwo}{\Qs{\intersect{G}{H}2}}
\newcommand{\QGbarone}{\Qs{\comp{G}1}}
\newcommand{\QGbartwo}{\Qs{\comp{G}2}}
\newcommand{\QGbarthree}{\Qs{\comp{G}3}}
\newcommand{\QHbarone}{\Qs{\comp{H}1}}
\newcommand{\QHbartwo}{\Qs{\comp{H}2}}
\newcommand{\QHbarthree}{\Qs{\comp{H}3}}

\newcommand{\PBBCone}{\Ptf{B}{\intersect{B}{C}1}}
\newcommand{\PBBCtwo}{\Ptf{B}{\intersect{B}{C}2}}
\newcommand{\PCBCone}{\Ptf{C}{\intersect{B}{C}1}}
\newcommand{\PCBCtwo}{\Ptf{C}{\intersect{B}{C}2}}
\newcommand{\PGGbartwo}{\Ptf{G}{\comp{G}2}}
\newcommand{\PGGbarthree}{\Ptf{G}{\comp{G}3}}
\newcommand{\PGGHtwo}{\Ptf{G}{\intersect{G}{H}2}}
\newcommand{\PHHbarthree}{\Ptf{H}{\comp{H}3}}
\newcommand{\PHGHtwo}{\Ptf{H}{\intersect{G}{H}2}}

\newcommand{\PGcapA}{\Ps{G\cap A}}
\newcommand{\PGcapAbar}{\Ps{\bar{G\cap A}}}
\newcommand{\QGcapAA}{\Qft{G\cap A}{A}}
\newcommand{\QGcapAG}{\Qft{G\cap A}{G}}
\newcommand{\QGcapAbarG}{\Qft{\comp{G\cap A}}{G}}

\newcommand{\QABbulletA}{\Qft{\intersect{A}{B}\bullet}{A}}

\newcommand{\RGbarAG}{R_{\comp{\BG},\intersect{A}{\BG}}}
\newcommand{\RGbarGbar}{R_{\comp{\BG},\comp{\BG}}}
\newcommand{\RGHAG}{R_{\intersect{\BG}{\CH},\intersect{A}{\BG}}}
\newcommand{\RGHGbar}{R_{\intersect{\BG}{\CH},\comp{\BG}}}

\newcommand{\FAbarbullet}{F_{\comp{A}\bullet}}
\newcommand{\FHbar}{F^T_{\comp{\CH}}}

\newcommand{\EA}{E\!\!A}
\newcommand{\FB}{F\!\!B}

\newcommand{\starorval}[1]{\ast}
\newcommand{\proofmarker}[1]{\vspace{3mm}\noindent\textit{#1}}

\newcommand{\PEbar}{\Ps{\comp{\EA}}}
\newcommand{\PEF}{\Ps{\intersect{\EA}{\FB}}}
\newcommand{\PFbar}{\Ps{\comp{\FB}}}
\newcommand{\QFbar}{\Qs{\comp{\FB}}}
\newcommand{\QFC}{\Qs{\intersect{\FB}{C}}}
\newcommand{\QCbarbullet}{\Qs{\comp{C}\bullet}}

\newcommand{\QEbarE}{\Qft{\comp{\EA}}{\EA}}
\newcommand{\QEFE}{\Qft{\intersect{\EA}{\FB}}{\EA}}
\newcommand{\QEFF}{\Qft{\intersect{\EA}{\FB}}{\FB}}
\newcommand{\QFbarF}{\Qft{\comp{\FB}}{\FB}}
\newcommand{\PCFC}{\Ptf{C}{\intersect{\FB}{C}}}
\newcommand{\PCCbarbullet}{\Ptf{C}{\comp{C}\bullet}}

\newcommand{\RFbarEF}{R_{\comp{\FB},\intersect{\EA}{\FB}}}
\newcommand{\RFCEF}{R_{\intersect{\FB}{C},\intersect{\EA}{\FB}}}
\newcommand{\RFbarFbar}{R_{\comp{\FB},\comp{\FB}}}
\newcommand{\RFCFbar}{R_{\intersect{\FB}{C},\comp{\FB}}}

\newcommand{\FCbarbullet}{F^T_{\comp{C}\bullet}}
\newcommand{\FEbar}{F_{\comp{\EA}}}

\DeclareMathOperator*{\argmax}{argmax}

\newcommand{\QABAnum}[1]{\Qft{\intersect{A}{B}}{A#1}}
\newcommand{\QEFEnum}[1]{\Qft{\intersect{\EA}{\FB}}{\EA #1}}
\newcommand{\QstarA}{\Qft{\ast}{A}}
\newcommand{\QABBnum}[1]{\Qft{\intersect{A}{B}}{B#1}}
\newcommand{\QEFFnum}[1]{\Qft{\intersect{\EA}{\FB}}{\FB #1}}
\newcommand{\QstarB}{\Qft{\ast}{B}}

\newcommand{\QFCnum}[1]{\Qs{\intersect{\FB}{C}#1}}
\newcommand{\QFbarnum}[1]{\Qs{\comp{\FB}#1}}

\newcommand{\PFFbar}{\Ptf{\FB}{\comp{\FB}}}
\newcommand{\PFFC}{\Ptf{\FB}{\intersect{\FB}{C}}}

\newcommand{\PCCbarnum}[1]{\Ptf{C}{\comp{C}#1}}
\newcommand{\FCbarnum}[1]{F^T_{\comp{C}#1}}

\newcommand{\QGHsym}[1]{\Qs{\intersect{\BG}{\CH}(#1)}}
\newcommand{\QGbarsym}[1]{\Qs{\comp{\BG}(#1)}}
\newcommand{\QHbarsym}[1]{\Qs{\comp{\CH}(#1)}}

\newcommand{\QBCnum}[1]{\Qs{\intersect{B}{C}#1}}
\newcommand{\QGHnum}[1]{\Qs{\intersect{G}{H}#1}}

\newcommand{\QGbarnum}[1]{\Qs{\comp{G}#1}}
\newcommand{\QHbarnum}[1]{\Qs{\comp{H}#1}}

\newcommand{\FHbarnum}[1]{F^T_{\comp{\CH}#1}}

\usepackage{tikz}

\DeclareMathOperator{\SSS}{SSS}
\DeclareMathOperator{\IDV}{DV}
\DeclareMathOperator{\CSS}{CSS}

\declaretheorem[name=Theorem,numberwithin=section]{theorem}
\declaretheorem[name=Proposition,numberlike=theorem]{proposition}
\declaretheorem[name=Conjecture]{conjecture}
\declaretheorem[name=Lemma,numberlike=theorem]{lemma}
\declaretheorem[name=Corollary,numberlike=theorem]{corollary}

\theoremstyle{remark}
\declaretheorem[name=Remark,qed={\lower-0.3ex\hbox{$\diamond$}},numberwithin=section]{remark}
\declaretheorem[name=Example,qed={\lower-0.3ex\hbox{$\triangleleft$}},numberwithin=section]{example}

\theoremstyle{definition}
\newtheorem{definition}{Definition}
\newtheorem{problem}{Problem}
\newtheorem{alg}{Algorithm}

\newcommand*{\MyPath}{.}

\numberwithin{equation}{section}
\numberwithin{problem}{section}

\title{Graph-Induced Rank Structures and their Representations}

\author{S. Chandrasekaran\thanks{Department of Electrical and Computer Engineering, University of California Santa Barbara, Santa Barbara, CA 93106-5070, USA. (\email{shiv@ucsb.edu})} \and E. N. Epperly \thanks{Departments of Mathematics and Computer Science, University of California Santa Barbara, Santa Barbara, CA 93106-5070, USA. (\email{epperly@ucsb.edu})} \and N. Govindarajan \thanks{Department of Electrical Engineering (ESAT), KU Leuven, Kasteelpark Arenberg 10, B-3001 Leuven, Belgium. (\email{Nithin.Govindarajan@kuleuven.be})}}
\date{\today}

\usepackage{graphicx}

\begin{document}

\maketitle

\begin{abstract}
A new framework is proposed to study rank-structured matrices arising from discretizations of 2D and 3D elliptic operators. In particular, we introduce the notion of a graph-induced rank structure (GIRS) which describes the fine low-rank structures which appear in sparse matrices and their inverses in relation to their adjacency graph $\graph$. We show that the GIRS property is invariant under inversion, and hence any effective representation of the inverse of a GIRS matrix $A$ would be an effective way of solving $Ax = b$. We then propose an extension of sequentially semiseparable (SSS) representations to $\graph$-semiseparable ($\graph$-SS) representations defined on arbitrary graph structures, which possess a linear-time multiplication algorithm. We show the construction of these representations to be highly nontrivial by determining the minimal $\graph$-SS representation for the cycle graph $\graph$. To obtain a minimal representation, we solve an exotic variant of a low-rank completion problem.

\end{abstract}

\tableofcontents

\section{Introduction}\label{sec:introduction}

The solution of linear systems of equations $Ax = b$ is ubiquitous in applications. Particular interest has been paid to matrices $A$ arising from the discretization of partial differential equations (PDEs), especially of the elliptic type. For a general nonsingular $N\times N$ matrix, the system $Ax = b$ can be solved in $\mathcal{O}(N^3)$ operations by using Gaussian elimination. This can be further improved to $\mathcal{O}(N^\omega) $operations using sophisticated fast matrix-matrix multiplication algorithms; $\omega \ge 2.372\ldots$ is the current fastest known algorithm \cite{AW21}. \ethancheck \nithincheck \shivcheck

For matrices possessing certain structural characteristics, significantly faster algorithms can be devised. For a sparse matrix $A$ with adjacency graph $\graph$ possessing many fewer than $N^2$ nonzero elements, significant performance improvements can be gained by reordering the matrix to reduce fill-in by Gaussian elimination. Finding the optimal ordering for a general graph $\graph$ is an NP-hard combinatorial problem \cite{articleFillNPcomplete}. For chordal graphs $\graph$ such as the line graph, Gaussian elimination can be done with no fill-in, resulting in complexity $\mathcal{O}(N)$. For the 2D mesh graph, the nested dissection ordering \cite{george1973nested} results in complexity $\mathcal{O}(N^{3/2})$. A combinatorial argument \cite{hoffman1973complexity} shows that any elimination ordering of the 2D mesh graph results in $\Omega(N^{3/2})$ complexity, showing the nested dissection ordering is asymptotically optimal. The nested dissection ordering produces $\mathcal{O}(N^2)$ complexity in three dimensions, and this elimintion ordering is also asymptotically optimal. \nithincheck

Due to the superlinear complexity of sparse Gaussian elimination on sparse matrices arising from discretization of 2D and 3D elliptic PDEs, there has been considerable focus on preconditioned iterative methods for these problems \cite{saad2003iterative}. For problems in which it is possible to construct a preconditioner $P$ for which the condition number $\kappa(P^{-1}A)$ remains uniformly bounded in $N$, the linear system $Ax = b$ can be solved to a given fixed accuracy by performing $\mathcal{O}(1)$ iterations, each of which requires multiplying by $P$ and $A$. Provided multiplying by $P$ and $A$ can be done in $\mathcal{O}(N)$ operations, as can be done by multigrid preconditioners for certain matrices $A$ such as the 2D discrete Poisson problem, this results in an overall linear time complexity. However, the performance of these methods is very dependent on the spectral properties of the matrix $A$. For particular classes of symmetric positive definite $A$, many of these methods work quite well, but the design of preconditioners for $Ax=b$ is significantly more challenging for a general $A$, particularly if $A$ is nonsymmetric
or indefinite.

An additional line of inquiry was started by noting that many matrices $A$ arising in applications have the property that certain off-diagonal blocks possess low (numerical) rank. The class of such matrices includes sparse matrices arising from the finite element discretization of PDE's as well as many dense matrices, such as those obtained from integral equations, inverses of sparse matrices, and (transformations of) classical structured matrices. There have been a proliferation of different and closely related rank structures and representations exploiting those lo- rank off-diagonal blocks to develop fast algorithms. For example, FMM \cite{greengard1987fast}, SSS \cite{chandrasekaran2002fast,chandrasekaran2005some}, HSS \cite{chandrasekaran2006fast,xia2010fast}, $\mathcal{H}$- and $\mathcal{H}^2$-matrices \cite{borm2003introduction,hackbusch2002data,hackbusch2015hierarchical}, HODLR \cite{ambikasaran2013mathcal,aminfar2016fast}, among numerous others. A summary of the differences and relations between many of these structures is provided in the first three sections of \cite{ambikasaran2014inverse}. These rank-structured solvers usually proceed in two steps. First, a compressed representation of the matrix is constructed by means of computing low rank factorizations of the off-diagonal blocks. Next, from this representation a (compressed) factorization ($LU$, $QR$, $ULV$) of $A$ is computed, and the system $Ax = b$ is solved. In cases where effective (compressed) factorization are not known (e.g.,\ in FMM), a fast matrix-vector multiply using the compressed representation can always be employed to accelerate an iterative solver.

For HSS and SSS matrices, the total time complexity of solving $Ax = b$ once the representation has been computed is $\mathcal{O}(Nr^2)$, where $r$ is the maximum rank among some collection of off-diagonal blocks. For the 2D Poisson equation on the square discretized according to a 5-point finite difference scheme with the natural ordering, the rank $r$ is of order $r = \mathcal{O}(N^{1/2})$ and the total time complexity is thus $\mathcal{O}(N^2)$, worse than the time complexity of solving $A$ using sparse Gaussian elimination.
For HSS, the sparse Gaussian elimination complexity of $\mathcal{O}(N^{3/2})$ can be recovered using the nested dissection ordering, but the complexity remains superlinear. The problem is only more severe in 3D, where the complexity of solving the 3D Poisson problem using either HSS or sparse Gaussian elimination in the nested dissection ordering jumps to $\mathcal{O}(N^2)$. 

\paragraph{Rank-structured Solvers in 2D and 3D} The problem of developing efficient rank-structured solvers for 2D and 3D problems have thus been an active area of research. One natural idea to extend rank structures to higher dimensions is to use a multilevel compression scheme in which dense matrices in an SSS or HSS representation are themselves represented as SSS or HSS matrices \cite{xiamulti}. However, in order to compute factorizations for such multilevel SSS and HSS matrices, one must compute sums and products of SSS and HSS matrices in Schur complement, which may increase the size of the off-diagonal blocks by a constant factor. Thus, during the course of factorization in which many of these sums and products occur, the off-diagonal block ranks may dramatically increase, leading to an overall superlinear time complexity. Such degradation in performance may be mollified by recompressing the off-diagonal blocks provided they still retain low rank \cite{qiu2015efficient}, though there may be no guarantees that these ranks will indeed be low enough to provably improve the time complexity. \nithincheck \shivcheck

Another family of methods that has received considerable attention have been methods based on recursive skeletonization \cite{martinsson2005fast,ho2012fast}. A formulation based on generalized $LU$ factorizations is presented in \cite{ho2016hierarchicalDE,ho2016hierarchicalIE}. In this family of methods, interpolative decompositions---rank factorizations where a column basis is selected from among the columns of the matrix being compressed---are used to produce a representation of a discretized PDE or integral equation, often without ever storing the full matrix $A$. The hierarchical interpolative factorization (HIF) \cite{ho2016hierarchicalDE,ho2016hierarchicalIE} takes this idea to its natural extent, repeatedly using the recurive skeletonization idea to reduce the degrees of freedom from 3D or 2D all the way down to 0D. Similar to the multilevel SSS and HSS examples, the performance of these methods depends on the off-diagonal ranks not growing too high during the computation of Schur complements, which appears to be the case for many problems of interest. For problems in which the Schur complements behave nicely, the HIF method is shown to have quasilinear time complexity. \shivcheck

Yet another strategy for solving high-dimensional linear problems involves computing sparse approximate factorizations without explicitly using low-rank compressions of off-diagonal blocks \cite{SSO20,SKO21}. These algorithms provably have quasilinear time complexity for inversion of kernel matrices corresponding the Green's function of a fixed elliptic operator in a fixed domain of a fixed dimension. The analysis of these methods are very delicate and use the ellipticity of the underlying differential operator in an essential way, and generalizing the ideas behind these methods to different problem classes is a challenging open problem.

Despite good empirical performance and growing theoretical understanding of many methods in the literature, we still believe the question of what the ``natural'' way of characterizing and exploiting rank structures for 2D and 3D problems remains open in an important way. Specifically, we are interested in the following two broad questions:

\begin{enumerate}
    \item \textit{What is the ``right'' algebraic  structure that naturally characterizing the rank structure properties of operators in 2D and 3D?}\nithin{what do you mean with operators? You specifically have elliptic operators in mind correct?} Much of the existing work on rank-structured solvers in 2D and 3D has focused on discretizations of PDEs or integral equations. In these contexts, it is often not clear what the exact algebraic properties of the discretized operator are that are being used in the algorithm, and whether the properties being exploited are only algebraic in nature or also analytic. A great strength of existing rank-structured techniques in 1D like HSS, SSS, and $\mathcal{H}$-matrices is that they are fully algebraic. This allows these techniques to be useful in problems unrelated to PDEs such as solving Toeplitz systems of linear equations \cite{martinsson2005fastb,chandrasekaran2008superfast}. We believe that clearly establishing the exact algebraic rank structure being exploited by an algorithm will greatly help clarify their scope of applicability. To us, an ideal algorithm will be purely algebraic in that it will perform equally effectively for PDEs, integral equations, or for general Toeplitz-block-Toeplitz matrices in Fourier space (which possess a natural 2D generalization of the rank structure of Toeplitz matrices in Fourier space).
    
    \item \textit{Does there exist an algebraic representation of the inverse of a 2D or 3D operator which can provably be applied in linear time?} Following from the last point, we are unaware of an algorithm which provably computes the action of the inverse of a general 2D or 3D rank-structured matrix on a vector in linear time (or even quasilinear time for 3D). Analysis of existing algorithms relies on bounding the rank growth involved in Schur complement expressions computed during factorization. In general, such rank growth estimates could possibly rely on analytic information about the underlying PDE or integral equation being discretized, possibly limiting the application of such methods to certain classes of PDEs and integral equations. We believe that a robust and provably linear time algorithm for multiplying by the inverse of a general high-dimensional rank-structured matrix remains an important open problem.
\end{enumerate}


\paragraph{Main contributions} 
In this paper, we provide a new framework characterizing the rank-structure of 2D and 3D operators (providing a conjectural answer to Question~1) and introduce representations that could potentially satisfy Question~2. Specifically, we introduce the notion of graph-induced rank structure (GIRS) that, for instance, precisely capture the exact low-rank structures of a sparse matrix \textit{and} its inverse in terms of its adjacency graph $\graph$  (Section \ref{subsec:GIRS}).

It will follow that \textit{if} we can compute a representation of the inverse of a sparse matrix (or any other GIRS matrix or its inverse) such that (i) the size of the representation, $M$, is linear in the matrix size $N$ and (ii) matrix-vector multiplication can be done in linear time in $M$, then repeated matrix-vector multiplications can \emph{subsequently} be done in near linear time in $N$, yielding a fast solver for multiple right-hand side problems. So in this new way of thinking about the problem, the crucial bottleneck is revealed to be the rapid computation of a compact matrix representation that captures the GIRS property.

Addressing this bottleneck appears to be a difficult problem. Inspired by the connections between GIRS and Sequentially Semiseparable (SSS) representations, we propose \nousedv{a generalization}\usedv{two generalizations} of SSS representations to arbitrary graph topologies which we call\usedv{Dewilde-van der Veen (DV) representations (Section~\ref{subsec:DV}) and} $\graph$-semi-separable ($\graph$-SS) representations (Section~\ref{sec:edv}). We describe some useful properties that these representations satisfy. For example, we show that $\graph$-SS representations possesses a linear time multiplication algorithm, satisfying (ii). We also show that matrices admitting these representations possess the GIRS property. The ultimate goal is to show the converse: matrices which satisfy the GIRS property admit a compact representation.

To this end, we study the converse problem for a particular example, showing how to construct $\graph$-SS representations for the cycle graph (Section~\ref{sec:css}). It is shown that to obtain a minimal representation, one needs to solve a highly non-trivial rank completion problem for which we present an algorithm in Section \ref{sec:hankelminimization}. It is expected that for general graphs, more complicated variants of this problem needs to be solved to obtain a minimal representations. Finally, we demonstrate the performance of the CSS representation numerically (Section \ref{sec:numerical_results}) and round out the paper with some concluding remarks (Section \ref{sec:overview}).

It is our hope that the ideas presented in this paper will throw further light on the general problem of constructing fast direct solvers for rank-structured matrices in 2D and 3D.


\section{Graph-induced rank structured matrices and their representations}\label{sec:GIRS}

We are interested in solving $Ax = b$ where, in the normative example, $A$ is obtained by discretizing a PDE or integral equation on a 2D or 3D domain. Our central thesis is that the low-rank structure of such matrices can be most precisely captured by considering the mesh graph $\graph$ on which this equation was discretized.\ethancheck



\subsection{Graph-induced rank structured matrices} \label{subsec:GIRS}
 
Consider a graph $\mathbb{G}= (\mathbb{V}_{\mathbb{G}},\mathbb{E}_{\mathbb{G}})$ with vertex set $\mathbb{V}_\graph$ of cardinality $|\mathbb{V}_{\graph}|=n$ and edge set $\mathbb{E}_{\graph}$. Associate with each node $i\in \mathbb{V}$ the vectors $x_i,b_i\in \mathbb{R}^{N_{i}}$ and let $A\in \mathbb{R}^{N\times N}$ be a square matrix of size $N=N_{1}+N_2 + \ldots + N_n$ such that
\begin{equation}\label{eq:block_matrix_multiplication}
 b_{i} = \sum_{j=1}^{n} A_{ij} x_j, 
\end{equation}
where $A_{ij} \in \mathbb{R}^{N_i\times N_j}$. The pair $(A,\graph)$ is referred together as a \emph{graph-partitioned matrix}. \ethancheck 

The introduction of this mathematical object is motivated by the following observation. Suppose that  $\subgraph =  (\mathbb{V}_{\subgraph},\mathbb{E}_{\subgraph}) $ is a induced subgraph of $\graph$ (i.e.   $\mathbb{V}_{\subgraph} \subset \mathbb{V}_{\graph}$ and $\mathbb{E}_{\subgraph}$ contains all edges in $\mathbb{E}_{\graph}$ whose vertices belong to $\mathbb{V}_{\subgraph}$) and suppose that $\complementgraph$ is its induced complement---that is, the induced subgraph with the vertex set $\nodeset_{\complementgraph} = \nodeset_{\graph} \setminus \nodeset_{\subgraph}$. The graphs $\subgraph$ and $\complementgraph$ naturally define a partition of $A$. That is, one can permute the rows and columns of $A$ such that
\begin{equation}
    PAP^T = \begin{bmatrix}
    A_{\subgraph,\subgraph} &  A_{\subgraph,\complementgraph}  \\
    A_{\complementgraph,\subgraph}  & A_{\complementgraph,\complementgraph} 
    \end{bmatrix}.\ethancheck \label{eq:Apermuted}
\end{equation}
Here $A_{\subgraph_1, \subgraph_2}$ contains all blocks $A_{ij}$ with $i \in \nodeset_{\subgraph_1}$ and $j$ in $\nodeset_{\subgraph_2}$ for $\subgraph_1,\subgraph_2 \in \{ \subgraph, \complementgraph \} $. We call $A_{\complementgraph,\subgraph}$ and $A_{\subgraph,\complementgraph}$ the \textit{Hankel blocks} induced by $\subgraph$ and $\complementgraph$, respectively. This naturally extends the established meaning of the term ``Hankel block'' in the theory of SSS matrices. Note that, by a classical result (proven in, e.g., \cite[Lem.~4.2(b)]{hackbusch2004hierarchical}), the ranks of the Hankel blocks are exactly preserved under inversion.
\begin{proposition}[The Hankel block property] \label{prop:rankinvariance}
Let $B$ denote the inverse of \eqref{eq:Apermuted} $P A P^T$  and write
$$ B =   \begin{bmatrix}
    B_{\subgraph,\subgraph} &  B_{\subgraph,\complementgraph}  \\
    B_{\complementgraph,\subgraph}  & B_{\complementgraph,\complementgraph}. 
    \end{bmatrix}$$
Then
$$ \rank B_{\subgraph,\complementgraph} =  \rank A_{\subgraph,\complementgraph}\quad \mbox{and} \quad \rank B_{\complementgraph,\subgraph} = \rank A_{\complementgraph,\subgraph}. $$
\end{proposition}

Proposition~\ref{prop:rankinvariance} shows that the ranks of the Hankel blocks of $A$ are invariant under inversion. This clarifies that any low-rank Hankel blocks in a matrix shall also be present in its inverse.
This motivates the following definition.

\begin{definition}[Graph-Induced Rank Structure] \label{def:GIRS}
The pair $(A,\graph)$ is said to have a \textit{{graph-induced rank structure} (GIRS) property} if there exists a constant $c > 0$ such that
\begin{equation}\label{eq:girs_rank_bound}
    \rank(A_{\complementgraph,\subgraph}) \le c \rho(\subgraph)
\end{equation}
for all induced subgraphs $\subgraph$ of $\graph$, where $\rho(\subgraph)$ denotes the number of border edges connecting nodes in $\subgraph$ to $\complementgraph$.. \ethancheck
\end{definition}

\begin{remark}
Note that this also implies that  $\rank(A_{\subgraph,\complementgraph}) \le c \rho(\subgraph)$. \ethancheck
\end{remark}

\begin{figure}
    \centering
    \includegraphics[width=0.4\textwidth]{\MyPath/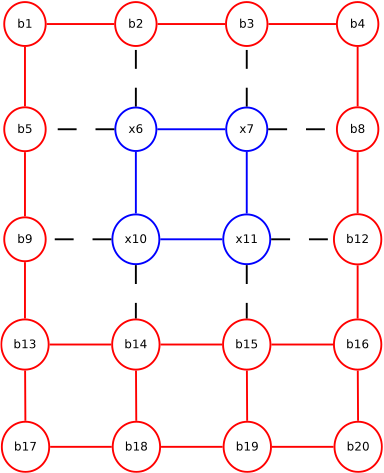}
    \caption{The $5\times 4$ 2D mesh graph $\graph$, the induced subgraph $\subgraph$ of $\{6,7,10,11\}$ (shown in blue), and the incluced complement $\bar{\subgraph}$ (shown in red). For this subgraph, there are eight border edges (shown dashed and in black), so $\rho(\subgraph) = 8$. If the matrix $A$ were to possess $\graph$ as its adjacency graph (such as the 2D discrete Poisson problem), then Proposition \ref{prop:sparse_girs} would show that $\rank A_{\bar{\subgraph},\subgraph} \le 8$. \ethancheck} 
    \label{fig:GIRSproperty}
\end{figure}

Figure~\ref{fig:GIRSproperty} illustrates the GIRS property for a 2D mesh graph. We shall refer to a pair $(A,\graph)$ as a GIRS-$c$ matrix, and shall omit the underlying graph $\graph$ when it is clear from context. Note that a matrix possessing graph-induced rank structure does not typically have low-rank off-diagonal blocks in the traditional sense. If $\graph$ is, for example, a 2D mesh graph and $\subgraph$ is generic rectangular subgraph, $\rho(\subgraph)$ will be $\rho(\subgraph) = \mathcal{O}(\sqrt{N})$ and it will not be algorithmically efficient to directly store a rank factorization of the Hankel block $A_{\complementgraph,\subgraph}$. What gives the GIRS property its strength is the requirement that the rank bound \eqref{eq:girs_rank_bound} hold for \textit{all} subgraphs $\subgraph$.

An important subclass of GIRS matrices are sparse matrices together with their adjacency graph.

\begin{proposition}\label{prop:sparse_girs}
Let $A$ be a sparse matrix with adjacency graph $\graph$. Then $(A,\graph)$ is a GIRS-1 matrix.
\end{proposition}

\begin{proof}
    For any subgraph $\subgraph$, 
    \begin{equation*}
        \rank (A_{\bar{\subgraph},\subgraph}) \le \nnz( A_{\bar{\subgraph},\subgraph}) = \rho(\subgraph),
    \end{equation*}
    where $\nnz$ denotes the number of nonzero entries in a matrix. The inequality follows from the rank that a matrix cannot have more linearly independent rows than nonzero entries and the equality is because every border edge in the adjacency graph is nonzero entry in $A_{\bar{\subgraph},\subgraph}$, by the definition of adjacency graph. \ethancheck 
\end{proof}

The following are natural and important extensions.

\begin{example}
Finite difference and finite element matrices satisfy the GIRS property with their adjacency graphs. \ethancheck 
\end{example}

\begin{example}
Banded matrices with bandwidth $k$ are GIRS-$k$ matrices. \ethancheck 
\end{example}

Proposition~\ref{prop:sparse_girs} shows that the GIRS property naturally expresses the low-rank structure possessed by sparse matrices. The upper bounds of the ranks of off-diagonal blocks in sparse matrices characterized by the GIRS property are precise for matrices arising from discretized PDEs for contiguous subgraphs $\subgraph$. \ethancheck 

\begin{example}[2D Poisson problem]
Consider the off-diagonal blocks of the 2D $\sqrt{N}\times \sqrt{N}$ Poisson problem in the natural order with node set $\nodeset_\graph = \{1,2,3,\ldots,N\}$ and subgraph $\subgraph_i$ with node set $\nodeset_{\subgraph_i} = \{1,2,\ldots,i\}$. Then
\begin{align*}
    \rank A_{\bar{\subgraph}_i,\subgraph_i} &= i < i+1 = \rho(\subgraph_i), & i < \sqrt{N}, \\
    \rank A_{\bar{\subgraph}_i,\subgraph_i} &= \sqrt{N} = \rho(\subgraph_i), & \sqrt{N} \le i \le N - \sqrt{N}, \\
    \rank A_{\bar{\subgraph}_i,\subgraph_i} &= N - i < N - i+1 = \rho(\subgraph_i), & N - \sqrt{N} < i.
\end{align*}
The bounds provided by the GIRS property for the off-diagonal Hankel blocks are very tight, off by at most one.\ethancheck
\end{example}

The fact that the GIRS property precisely captures the low-rank structure of the off-diagonal Hankel blocks of the 2D Poisson equation, the quintessential discretized PDE, demonstrates that the GIRS property may be a useful way of characterizing the rank structure of discretized PDEs. 

\begin{remark}
The GIRS property \textit{with $\graph$ taken to be the adjacency graph} is not excellent at expressing the low-rank structure of \textit{every} sparse matrix. However, a different graph other than the adjacency graph may capture this low-rank structure. This is demonstrated in Example~\ref{ex:arrowhead}. \ethancheck 
\end{remark}

\begin{example}[Arrowhead matrix]\label{ex:arrowhead}
Consider the arrowhead matrix $A$ where $A_{ij} \ne 0$ if, and only if, $i = j$, $i = 1$, or $j=1$. Then $A$'s adjacency graph is $\graph = (\nodeset_\graph,\edgeset_\graph)$ with $\nodeset_\graph = \{1,\ldots,N\}$ and $\edgeset_\graph = \{(1,j) : j \in \{2,\ldots,N\}\}$. However, for the contiguous subgraph $\subgraph_i$ with node set $\nodeset_{\subgraph_i} = \{1,2,\ldots,i\}$ for $1 < i < N$, we have $\rank A_{\bar{\subgraph_i},\subgraph_i} = 1 \ll N - i = \rho(\subgraph_i)$. Thus, the GIRS property with the graph $\graph$ is very poor at describing the ranks of the off-diagonal Hankel blocks for this matrix. However, $A$ is also GIRS-1 with the line graph $\mathbb{L}$ with $\nodeset_{\mathbb{L}} = \{1,\ldots,N\}$ and edge set $\edgeset_{\mathbb{L}} = \{\{ i,i+1\} : 1 \le i \le N-1 \}$, which does accurately predict the ranks of all the off-diagonal Hankel blocks.\ethancheck
\end{example}

The SSS and HSS rank structures are preserved under inversion, multiplication, and addition, and these properties are important in developing fast solvers for $Ax = b$. Thankfully, the GIRS property is preserved under these operations as well.

\begin{proposition}[Algebra of GIRS matrices]\label{prop:girs_algebra}
Let $(A,\graph)$ be a GIRS-$c$ matrix and $(B,\graph)$ be a GIRS-$d$. Then:
\begin{enumerate}[label={(\roman*)}]
    \item $(A^{-1},\graph)$ is a GIRS-$c$  matrix whenever $A$ is invertible,
    \item $(A+B,\graph)$ is a GIRS-($c+d$) matrix,
    \item $(AB,\graph)$ is a GIRS-($c+d$) matrix.
\end{enumerate}
\end{proposition}

\begin{proof}
    The conclusions of this theorem follow straightforwardly by multiplying the $2\times 2$ block partitions \eqref{eq:Apermuted} of $A$ and $B$ (similar to \cite[Lem.~4.2]{hackbusch2004hierarchical}). Statement (i) is an immediate consequence of Proposition~\ref {prop:rankinvariance}. \ethancheck Since $(A+B)_{\complementgraph,\subgraph} = A_{\complementgraph,\subgraph} + B_{\complementgraph,\subgraph}$,  statement (ii) follows from the fact that
    \begin{displaymath}
     \rank (A+B)_{\complementgraph,\subgraph}  \leq \rank A_{\complementgraph,\subgraph} + \rank B_{\complementgraph,\subgraph} \leq c \rho(\subgraph) + d \rho(\subgraph). \ethancheck 
    \end{displaymath}
    Finally for (iii), we note that
    \begin{displaymath}
    (AB)_{\complementgraph,\subgraph}  = A_{\complementgraph,\subgraph} B_{\subgraph, \subgraph} + A_{\complementgraph,\complementgraph} B_{\complementgraph,\subgraph}. \ethancheck
    \end{displaymath}
    Therefore,
    \begin{displaymath}
    \rank (AB)_{\complementgraph,\subgraph} \leq \rank A_{\complementgraph,\subgraph} B_{\subgraph, \subgraph} + \rank A_{\complementgraph,\complementgraph} B_{\complementgraph,\subgraph} \leq \rank A_{\complementgraph,\subgraph} + \rank B_{\complementgraph,\subgraph}, 
    \end{displaymath}
    from which the conclusion (iii) follows. \ethancheck 
    
\end{proof}

In particular, since the GIRS property is preserved under inversion, the GIRS property not only characterizes the low-rank structure of sparse matrices, but also their inverses. This closure of the rank structure property under inversion is usually necessary for a fast direct solver and was also exploit in the solvers HSS and SSS. \nithincheck

Our central question is whether there exists an efficient and effectively computable algebraic representation of GIRS matrices which can be leveraged to compute matrix-vector multiplications in linear time. If this were the case, then a representation of $A^{-1}$ could be computed and used as a linear-time direct solver for $Ax = b$. This would be particularly advantageous for repeated right-hand side problems, for even if it is expensive to compute such a representation for $A^{-1}$, the inverse operator could then be applied in linear time from that point forward. In the next section, we present a candidate representation which we conjecture may be able to answer this question in the affirmative.

\usedv{

\subsection{A motivating example: sequentially semi-separable matrices}\label{sec:sss}

In the case of GIRS matrix on the line graph, the formulation of an effective algebraic representation is well addressed by the theory of sequentially semi-separable (SSS) matrices \cite{chandrasekaran2002fast,chandrasekaran2005some}. A summary is provided in \cite[Sec.~3]{chandrasekaran2018fast}. Consider $(A,\graph)$ to be a GIRS matrix with $\mathbb{V}_{\graph} = \{1,2,\ldots,n \}$ and $\mathbb{E}_{\graph} = \{\{ i,i+1\} : 1 \le i \le n-1 \}$. An SSS representation for $A$ is given by a collection of matrices $(\{U_i \}^{n-1}_{i=1}$, $\{W_i  \}^{n-1}_{i=2}$, $\{V_i \}^{n}_{i=2}$, $\{ D_i \}^{n}_{i=1}$, $\{P_i\}^{n}_{i=2}$, $\{R_i  \}^{n-1}_{i=2}$, $\{Q_i \}^{n-1}_{i=1})$\ethancheck so that each block entry is expressed by
\begin{equation}
    \left[A\right]_{k\ell } = \begin{cases} 
    D_k & k=\ell  \\
    P_{k} R_{k-1} R_{k-2} \cdots R_{\ell+1} Q^T_{\ell}   & k> \ell \\
     U_{k} W_{k+1} W_{k+2} \cdots W_{\ell-1} V^T_{\ell}   & k< \ell 
    \end{cases} \ethancheck\label{eq:SSS}
\end{equation}
where $U_i\in \mathbb{R}^{N_i \times r^g_{i} }$, $W_i\in \mathbb{R}^{ r^g_{i-1} \times r^g_{i} }$, $V_i\in \mathbb{R}^{N_i \times r^g_{i-1} }$, $D_i\in \mathbb{R}^{N_i\times N_i}$, $P_i\in \mathbb{R}^{N_i \times r^{h}_{i-1} }$, $R_i\in \mathbb{R}^{r^h_{i}\times r^h_{i-1}}$,  and $Q_i\in \mathbb{R}^{N_i \times r^h_{i}}$.\ethancheck For example, in the case of $n=4$ the SSS representation reduces to
\begin{equation}
A = \begin{bmatrix}
{D_1  }              &          U_1V_2^T             &    U_1 W_2 V_3^T  &    U_1 W_2 W_3 V_4^T    \\
P_{2} Q_1^T               &  {D_2}          &  U_2V_3^T & U_2 W_3 V_4^T  \\
P_{3} R_{2} Q^T_1          &  P_{3} Q^T_2         &  {D_3} & U_3V_4^T      \\
P_{4} R_{3} R_{2} Q^T_1     &  P_{4} R_{3} Q^T_2  & P_{4} Q^T_3   & {D_4} 
\end{bmatrix}. 
\label{eq:SSSSn4}
\end{equation}
SSS matrices were first studied by Dewilde and van der Veen \cite{dewilde1998time} in the context of systems theory. In this framework, the entries of the SSS representation \eqref{eq:SSS} can be seen as the result of decomposing $A$ as the sum of a causal and anti-causal Linear Time Variant (LTV) system with input sequence $\{ x_i \}^{n=1}_{i=1}$ and outputs  $\{ b_i \}^{n=1}_{i=1}$. The anti-causal LTV system is described by the recursion 
\begin{equation} 
g_{k}  =  V^T_{k} x_{k} + W_{k}g_{k+1}, \quad k=n-1,n-2,\ldots,2 \ethancheck \ethancheck \label{eq:anticausal} 
\end{equation}
with terminal condition $g_n = V^T_n x_n$,\ethancheck whereas the causal LTV system is given by 
\begin{equation}
h_k = Q^T_k x_k + R_{k} h_{k-1}, \quad k=2,3,\ldots,n-1\ethancheck \label{eq:causal}  
\end{equation}
with initial condition $h_1 = Q^T_1 x_1$.\ethancheck The output equation reads
\begin{equation}
b_k  =   D_k x_k + U_{k} g_{k+1} + P_{k}h_{k-1}, \quad k=1,\ldots,n.  \label{eq:output}  
\end{equation}
\begin{remark}
Given that there is not really a natural ordering on arbitrary graphs, in later sections, where we generalize SSS to arbitrary graphs, the terms causal and anti-causal are replaced with ``upstream" and ``downstream", respectively.
\end{remark}
Direct execution of (\ref{eq:anticausal}-\ref{eq:output}) constitute the fast matrix-vector multiplication algorithm for a SSS representation. The computations involved in the multiplication algorithm can be depicted in a signal flow diagram as illustrated in figure~\ref{subfig:signalflowSSS}  for the case $n=5$. Due to its origins in systems theory, the collection of equations (\ref{eq:anticausal}-\ref{eq:output}) are referred to as the \textit{state-space equations} and the auxilliary variables $g_k$ and $h_k$ as the \textit{state space variables}.\ethancheck
\begin{figure}
    \centering
    \begin{subfigure}[b]{0.8\textwidth}
    \centering
    \includegraphics[width = 1\textwidth]{\MyPath/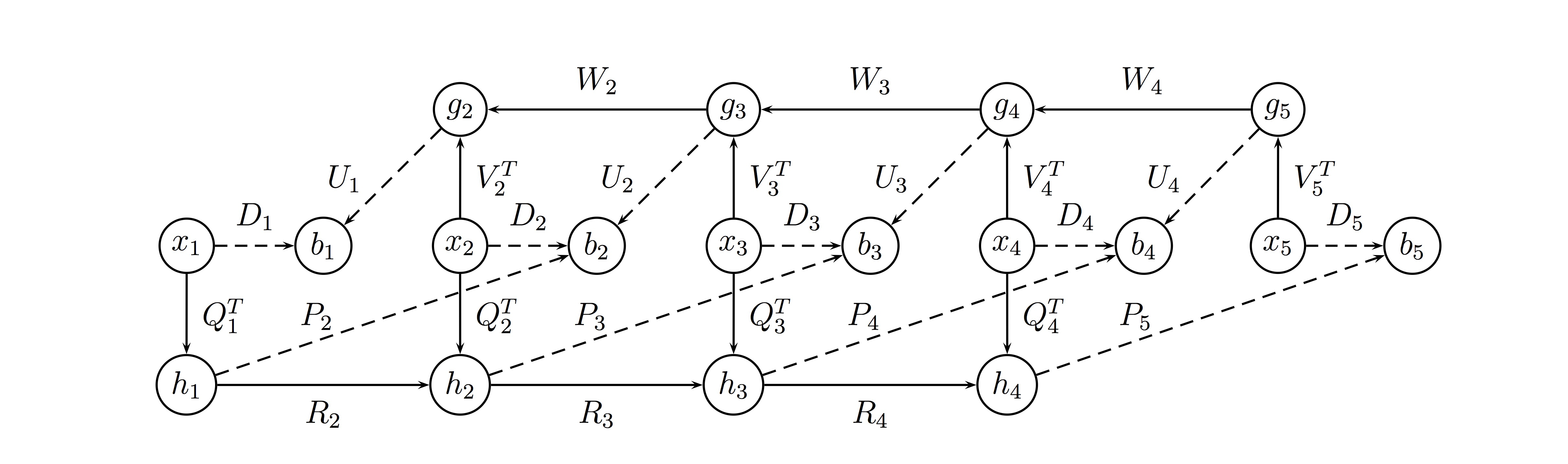}
    \caption{\nithincheck}
    \label{subfig:signalflowSSS}
    \end{subfigure}
    
    \begin{subfigure}[b]{0.8\textwidth}
        \begin{displaymath}
        \resizebox{0.99\textwidth}{!}{$
        \left[\begin{array}{ccccc|ccccc|ccccc}
        I & -W_1 & & &&&&&&& -V_1^T &&&& \\
        & I & -W_2 & &&&&&&&& -V_2^T &&& \\
        && I & -W_3 &&&&&&&&& -V_3^T &&\\
        &&& I & -W_4 &&&&&&&&& -V_4^T &\\
        &&&& I &&&&&&&&&& -V_5^T\\\hline 
        &&&&& I &&&&& -Q_1^T \\
        &&&&& -R_2 & I &&&&& - Q_2^T \\
        &&&&&& -R_3 & I &&&&& - Q_3^T \\
        &&&&&&& -R_4 & I &&&&& - Q_4^T \\
        &&&&&&&& -R_5 & I&&&&& - Q_5^T\\\hline
        & U_1 &&&&& &&&& D_1 \\
        && U_2 &&& P_2 &&&&&& D_2 \\
        &&& U_3 &&& P_3 &&&&&& D_3 \\
        &&&& U_4 &&& P_4 &&&&&& D_4 \\
        &&&&& &&& P_5 &&&&&& D_5
        \end{array} \right]
        $}
        \end{displaymath}
        \caption{} \label{subfig:SSS_matrix}
    \end{subfigure}
 
    \begin{subfigure}[b]{0.8\textwidth}
        \begin{displaymath}
        \resizebox{0.99\textwidth}{!}{$
        \left[\begin{array}{ccc;{2pt/2pt}ccc;{2pt/2pt}ccc;{2pt/2pt}ccc;{2pt/2pt}ccc}
        I && -V_1^T&-W_1&&&&&&&&&&& \\
        &I & -Q_1^T &&&&&&&&&&&&\\
        &&D_1&U_1 &&&&&&&&&&& \\ \hdashline[2pt/2pt]
         & && I & & -V_2^T&-W_2&&&&&&&&\\
        &-R_2&&& I & -Q_2^T &&&&&&&&& \\
         &P_2&&&& D_2 &U_2&&&&&&&&\\\hdashline[2pt/2pt]
        &&&  &&& I && -V_3^T & -W_3&&&&& \\
        &&&&-R_3&&& I & -Q_3^T &&&&&&  \\
        &&& &P_3&& && D_3 & U_3 &&&&&  \\ \hdashline[2pt/2pt]
        &&& &&& &&& I && -V_4^T & -W_4 \\
        &&& &&& & -R_4 && & I & -Q_4^T \\
        &&& &&& & P_4 && && D_4 & U_4 \\ \hdashline[2pt/2pt]
        &&& &&& &&& &&& I && -V_5^T \\
        &&& &&& &&& & -R_5 && & I & -Q_5^T \\
        &&& &&& &&& & P_5 && && D_5
        \end{array} \right]
        $}
        \end{displaymath}
        \caption{} \label{subfig:SSS_reordered}
    \end{subfigure}
    \caption{The signal flow graph of the SSS representation, shown in (\subref{subfig:signalflowSSS}), demonstrates how the state space equations (\ref{eq:anticausal}-\ref{eq:output}) can be used to compute the product $Ax$. An SSS matrix can be seen as the Schur complement of the lifted sparse matrix \eqref{eq:toplevelSSS}, shown in (\subref{subfig:DV_reordered}). The linear system $Ax = b$ can be solved by ordering \eqref{eq:toplevelSSS} to \eqref{eq:SSSreordered} by merging the nodes $x_i$, $g_i$, and $h_i$. The reordered system \eqref{eq:SSSreordered} is shown in (\subref{subfig:SSS_reordered}). Notice that the line graph structure of the multiplication signal flow diagram (\subref{subfig:signalflowSSS}) or equivalently the GIRS graph $\graph$ is reflected by the block tri-diagonal structure in \eqref{eq:SSSreordered}. }\label{fig:SSS}
\end{figure}

The equations (\ref{eq:anticausal}-\ref{eq:output}) may also be summarized by a single top level matrix notation. This is done by introducing the dummy variables  $g_n, h_1 \in\mathbb{R}^{0}$, $V_1 \in \mathbb{R}^{ N_1 \times 0}$, $W_1 \in \mathbb{R}^{0 \times r^g_1}$, $W_n\in \mathbb{R}^{r^g_{n-1} \times 0}$, $U_n \in \real^{N_n\times 0}$, $P_1\in \real^{N_1\times 0}$, $R_1 \in \mathbb{R}^{r^h_1\times 0} $, $R_n \in \real^{0\times r^h_{n-1}}$, $Q_n \in \real^{N_n\times 0}$\ethancheck and letting $\diagmat{U} = \diag \{ U_i \}^{n}_{i=1}$, $\diagmat{W} = \diag \{ W_i \}^{n}_{i=1}$, $\diagmat{D} = \diag \{ D_i \}^{n}_{i=1}$, etc.  we obtain
\begin{equation}
    \begin{bmatrix} I - \diagmat{W} Z^T &       0     & -\diagmat{V}^T \\
                 0           &   I-\diagmat{R}Z   & -\diagmat{Q}^T  \\
                   \diagmat{U}Z^T      &     \diagmat{P}Z     &  \diagmat{D}           \end{bmatrix} \begin{bmatrix} g \\ h \\ x  \end{bmatrix} = \begin{bmatrix} 0 \\ 0 \\ b \end{bmatrix},\ethancheck \label{eq:toplevelSSS}  
\end{equation}
where $Z$ denotes the downshift operator
\begin{equation*}
Z = \begin{bmatrix} 
0 & \\
I & 0\\
  & I & \ddots \\
  &   & \ddots & 0 \\
  &   &        & I & 0
\end{bmatrix}.\ethancheck
\end{equation*}
Performing Gaussian elimination on \eqref{eq:toplevelSSS}, one can express $A$ as the Schur complement of the matrix in \eqref{eq:toplevelSSS}, in effect deriving an alternative way of writing \eqref{eq:SSS}:
\begin{equation}
A = \diagmat{D} +  \diagmat{U} Z^T (I-\diagmat{W}Z^T)^{-1} \diagmat{V}^T +   \diagmat{P}Z (I- \diagmat{R}Z)^{-1} \diagmat{Q}^T,\ethancheck
\label{eq:diagrepSSS}
\end{equation}
which is referred to as the \emph{diagonal representation} of a SSS matrix. We may concisely denote the SSS representation by
\begin{equation}\label{eq:SSS_of}
A =  \SSS(\diagmat{U},\diagmat{W},\diagmat{V}, \diagmat{D}, \diagmat{P}, \diagmat{R}, \diagmat{Q}).\ethancheck
\end{equation}
The dimensions of the SSS representation are dictated by the numbers $r^g_{i}$ and $r^h_{i}$. Matrix vector multiplication $Ax$ can be evaluated in $\mathcal{O}(\max_{i}\{ r^g_{i}  , r^h_{i}   \}  N)$ operations.  When the value of $\max_{i}\{ r^g_{i}  , r^h_{i}   \}$ can be kept bounded for a family of matrices of increasing size (e.g. derived from a PDE discretization), multiplication can essentially be performed in linear time.\ethancheck Given a GIRS-$r$ matrix on the line graph, the following theorem shows that
\begin{displaymath}
\max_{i}\{ r^g_{i}  , r^h_{i}   \} \leq r.\ethancheck
\end{displaymath}
That is, GIRS matrices associated with the line graph are SSS matrices.

\begin{theorem} \label{thm:minimalSSS}
Let $(A,\graph)$ be a GIRS-$r$ matrix with $\mathbb{V}_{\graph} = \{1,2,\ldots,n \}$ and $\mathbb{E}_{\graph} = \{ (i,i+1) : 1\le i \le n-1\} $. \nithincheck and $\mathcal{G}_i = A_{\subgraph_i, \complementgraph_i}$ the Hankel blocks associated with the induced subgraphs $\mathbb{H}_i$ with vertex set $\mathbb{V}_{\subgraph_i}  = \{1,2,\ldots,i \}$.    Then, there exists a SSS representation 
\begin{equation*}
A = \SSS(\diagmat{U},\diagmat{W},\diagmat{V}, \diagmat{D}, \diagmat{P}, \diagmat{R}, \diagmat{Q})
\end{equation*}
whose dimensions satisfy
\begin{equation} {r^h_{i}}  = \rank \mathcal{H}_i \leq r , \quad  {r^g_{i}} = \rank \mathcal{G}_{i} \leq r, \quad i=1,2\ldots,n-1.  \label{eq:rankSSS}
\end{equation}
Furthermore, a representation of these dimensions is optimal in the sense that for any other representation
$
A = \SSS(\hat{\diagmat{U}},\hat{\diagmat{W}},\hat{\diagmat{V}}, \hat{\diagmat{D}}, \hat{\diagmat{P}}, \hat{\diagmat{R}}, \hat{\diagmat{Q}}),
$
the dimensions must satisfy the inequality
\begin{displaymath}
{r^h_{i}} \leq \hat{r}^h_i, \quad {r^g_{i}} \leq \hat{r}^g_i.
\end{displaymath}
\end{theorem}	

\begin{proof}
	The fact that $\rank \mathcal{H}_i \leq r$ and $\rank \mathcal{G}_{i} \leq r$ follows from the GIRS property and the fact that $\rho(\subgraph) = 1$.\ethancheck  To show the existence of an SSS representation such that ${r^h_{i}} = \rank \mathcal{H}_{i}$, first note that
	$$
	\mathcal{H}_i = \begin{bmatrix}
	P_{i+1} \\ 
	P_{i+2}R_{i+1} \\
	P_{i+3} R_{i+2} R_{i+1} \\
	\vdots \\
	P_{n} R_{n-1} R_{n-2} \cdots R_{i+1}
	\end{bmatrix}  \begin{bmatrix} R_{i} \cdots R_{3} R_{2} Q_{1}^T  & \cdots &  R_{i} R_{i-1}  Q_{i-2}^T & R_{i} Q_{i-1}^T & Q_i^T
	\end{bmatrix}.\ethancheck
	$$
	which shows $\rank \mathcal{H}_i \le r_i^h$ for any SSS representation of $A$.  To achieve the equality ${r^h_{i}} = \rank \mathcal{H}_{i}$, consider the following construction. Compute a low-rank factorization of $\mathcal{H}_i$ as
	$$ \mathcal{H}_i = X_i Y^T_i = 
	\begin{bmatrix} X_{i,1}  \\ X_{i,2} \end{bmatrix} 
	\begin{bmatrix} Y^T_{i,1}  & Y^T_{i,2} \end{bmatrix} $$
	where $X_{i,1}\in\mathbb{R}^{ N_{i+1} \times  \rank \mathcal{H}_i}$ and $Y^T_{i,2}\in\mathbb{R}^{   \rank \mathcal{H}_i\times N_i }$. Set $P_{i+1} = X_{i,1}  $ and $Q_{i} = Y_{i,2} $ for $i=1,2,\ldots, n-1$. Furthermore, denote $X^{\dagger}_{i+1} $  to be the pseudo-inverse of $ X_{i+1}$ and set $R_{i+1} = X^{\dagger}_{i+1} X_{i,2}$ for $i=1,2,\ldots, n-2$. \nithincheck This procedure produces an SSS representation of the desired size
	$$ {r^h_{i}}  = \rank \mathcal{H}_i, \quad  {r^g_{i}} =  \rank \mathcal{G}_i, \quad i=1,2\ldots,n-1. 
    $$
    The analogous results for the upper triangular Hankel blocks follow the exact procedure.
\end{proof}

The previous theorem shows that GIRS matrices are SSS. The following theorem shows the converse: SSS matrices of size at most $r$ coincide exactly with GIRS-$r$ matrices on the line graph. \ethancheck

\begin{theorem}\label{thm:SSStoGIRS}
Suppose that a graph-partitioned matrix $(A,\graph)$ with $\mathbb{V}_{\graph} = \{1,2,\ldots,n \}$ and $\mathbb{E}_{\graph} = \{ (i,i+1) : 1\le i \le n-1\}$ has  a SSS representation 
$
A = \SSS(\diagmat{U},\diagmat{W},\diagmat{V}, \diagmat{D}, \diagmat{P}, \diagmat{R}, \diagmat{Q})
$ 
whose dimensions are given by
${r^h_{i}}$ and  ${r^g_{i}}$ for $ i=1,2\ldots,n-1$. Then $(A,\graph)$   is a GIRS-$r$ matrix with
\begin{equation}
r = \max_{i=1,2,\ldots,n-1}\{ r^g_{i}  , r^h_{i}   \}. \label{eq:rdef}
\end{equation}
Furthermore, if the representation 
$
A = \SSS(\diagmat{U},\diagmat{W},\diagmat{V}, \diagmat{D}, \diagmat{P}, \diagmat{R}, \diagmat{Q})
$ is optimal, then the constant $r$, as defined above, presents the best possible GIRS constant for $(A,\graph)$.
\end{theorem}

\begin{proof}
To prove our claim, we must confirm that  $\rank(A_{\complementgraph,\subgraph}) \le r \rho(\subgraph)$ for each subgraph  $\subgraph \subset \graph$. This fact is easy to confirm for a subgraph consisting of a single connected component, which then later can be generalized to a general subgraph.\ec  

Indeed, if we  have $\subgraph = \{ s,\ldots, t \}$    with $1 \leq s \leq t \leq n$, we may break up the complement graph $\complementgraph$ further into two disjoint sub-graphs: the  "upstream" complement $u(\complementgraph, \subgraph)$ consisting of all nodes $k\in\complementgraph$ with $ k > t $, and the  "downstream" complement ${d(\complementgraph,\subgraph)}$   consisting of all nodes $k\in\complementgraph$ with $k < s$.
We have the inequality
\begin{equation*}
\rank(A_{\complementgraph,\subgraph})   =   \rank \begin{bmatrix}A_{u(\complementgraph, \subgraph), \subgraph } \\ A_{d(\complementgraph, \subgraph), \subgraph}  \end{bmatrix} \leq \rank(A_{u(\complementgraph, \subgraph),  \subgraph} ) +  \rank(A_{d(\complementgraph, \subgraph) , \subgraph}).\ec 
\end{equation*} 
If $u(\complementgraph, \subgraph)$ is not the empty graph, we may factorize $A_{u(\complementgraph, \subgraph),  \subgraph}$ as
\begin{displaymath}
A_{u(\complementgraph, \subgraph), \subgraph}  =  \begin{bmatrix}
	P_{t+1} \\ 
	P_{t+2}R_{t+1} \\
	P_{t+3} R_{t+2} R_{t+1} \\
	\vdots \\
	P_{n} R_{n-1} R_{n-2} \cdots R_{t+1}
	\end{bmatrix}  \begin{bmatrix} R_{t} \cdots R_{s+2} R_{s+1} Q_{s}^T  & \cdots &  R_{t} R_{t-1}  Q_{t-2}^T & R_{t} Q_{t-1}^T & Q_t^T
	\end{bmatrix},\ec 
\end{displaymath}
showing that $\rank(A_{u(\complementgraph, \subgraph), \subgraph}) \le  r^g_t$.   Similarly,  if $d(\complementgraph, \subgraph)$ is not the empty graph, we may factorize $\rank(A_{d(\complementgraph, \subgraph), \subgraph})$ as
\begin{displaymath}
A_{d(\complementgraph, \subgraph), \subgraph}	=  \begin{bmatrix}
	U_{1} W_{2} W_{3} \cdots W_{s-1}  \\ 
		\vdots \\
	U_{s-3} W_{s-2} W_{s-1}\\
	U_{s-2} W_{s-1}  \\
	U_{s-1} 
	\end{bmatrix}  \begin{bmatrix} V_s^T & W_{s} V_{s+1}^T    & W_{s} W_{s+1}  V_{s+2}^T  &  \cdots  & W_{s}  W_{s+1} \cdots W_{t-1} V_{t}^T    \end{bmatrix}, \ec 
\end{displaymath}
showing that  $\rank(A_{d(\complementgraph, \subgraph), \subgraph }) \le  r^h_s $.\ec Overall, we have
\begin{displaymath}
\rank(A_{\complementgraph,\subgraph})  \leq r^h_{t}  + r^g_{s}  \leq 2  \max \{ r^h_{t} , r^g_{s}   \}   \leq \rho(\subgraph) r.\ec 
\end{displaymath}
Note that in the special edge cases where either $u(\complementgraph, \subgraph)$ or $d(\complementgraph, \subgraph)$ are empty, the desired inequality still holds. \ec

For a general subgraph  $\subgraph \subset \graph$, we observe that $\subgraph$ can be broken into its connected components: 

%
\begin{displaymath}
\subgraph = \bigcup^{p}_{i=1} \subgraph_i,\ec
\end{displaymath}
with $\subgraph_i$ consisting of the vertex set $\mathbb{V}_{\subgraph_i} =  \{ s_i,\ldots, t_i \}$ and edge set $\mathbb{E}_{\subgraph_i} = \{ (s_i,s_i+1),\ldots ,(t_i-1,t_i)  \}$. By partitioning $A_{\complementgraph,\subgraph}$ in conjunction with the connected components, it easily follows that
\begin{displaymath}
\rank(A_{\complementgraph,\subgraph})  =  \rank \begin{bmatrix} A_{\complementgraph,\subgraph_1} \\
A_{\complementgraph,\subgraph_2} \\
\vdots \\
A_{\complementgraph,\subgraph_p}
\end{bmatrix}   \leq \sum_{i=1}^{p} \rank(A_{\complementgraph,\subgraph_i}).\ec 
\end{displaymath}
Let $u(\complementgraph,\subgraph_i)$ denote a subgraph of $\complementgraph$ consisting of the nodes $k\in \complementgraph$ with $k>t_i$. Similarly, let $d(\complementgraph,\subgraph_i)$ denote a subgraph of $\complementgraph$ consisting of the nodes $k\in \complementgraph$ with $k<s_i$. We recognize that
$$ u(\complementgraph,\subgraph_i) \subseteq u(\complementgraph_i,\subgraph_i) , \qquad  d(\complementgraph,\subgraph_i) \subseteq d(\complementgraph_i,\subgraph_i) .\ec $$
%
%
This allows us to derive the following set of inequalities:
\begin{eqnarray*}
\rank(A_{\complementgraph,\subgraph}) &  \leq & \sum_{i=1}^{p} \rank(A_{\complementgraph,\subgraph_i}) \\
& \leq & \sum_{i=1}^{p} \rank(A_{u(\complementgraph,\subgraph_i),\subgraph_i}) +  \rank(A_{d(\complementgraph,\subgraph_i),\subgraph_i})  \\
& \leq & \sum_{i=1}^{p} \rank(A_{u(\complementgraph_i,\subgraph_i),\subgraph_i}) +  \rank(A_{d(\complementgraph_i,\subgraph_i),\subgraph_i}).\ec
\end{eqnarray*}
Applying the result for a single connected component and observing that $ \rho(\subgraph ) = 2p $,\ec we have
\begin{displaymath}
\rank(A_{\complementgraph,\subgraph}) \leq \sum_{i=1}^{p} \rank(A_{u(\complementgraph_i,\subgraph_i),\subgraph_i}) +  \rank(A_{d(\complementgraph_i,\subgraph_i),\subgraph_i})  \leq 2p   \max_{i=1,\ldots,p} \{ r^h_{t} , r^g_{s}   \}   \leq    \rho(\subgraph) r\ec
\end{displaymath}
Note we did not handle edge cases where $d$ and $u$ are empty, but, similar to above, the argument still passes through. \ec

To show that the constant $r$ presents the best possible GIRS constant for $(A,\graph)$, we simply observe that if there would exist a better bound, then this would immediately contradict the definition of $r$ in  \eqref{eq:rdef}. \nithincheck
\end{proof}

Theorems~\ref{thm:minimalSSS} and \ref{thm:SSStoGIRS} show an equivalence between GIRS on the line and SSS representations. That is, if a matrix is GIRS, it has a SSS representation. Vice versa, if it has a compact SSS representation, then the  matrix is GIRS on the line graph. This in effect relates a constructive description of SSS matrices (the existence of a compact SSS representation \eqref{eq:SSS}) to a purely algebraic characterization of SSS matrices (the GIRS property with the line graph).  This is a pattern we hope to extend to general GIRS matrices.  \ec

SSS representations possess very nice algebraic properties. In particular, the inverse of an SSS matrix is an SSS matrix with the same Hankel block ranks.\ec 

\begin{proposition}[SSS algebra] \label{prop:SSSalgebra}
Let: 
\begin{displaymath}
A = \SSS(\diagmat{U}_A,\diagmat{W}_A,\diagmat{V}_A, \diagmat{D}_A, \diagmat{P}_A, \diagmat{R}_A, \diagmat{Q}_A), \quad \SSS(\diagmat{U}_B,\diagmat{W}_B,\diagmat{V}_B, \diagmat{D}_B, \diagmat{P}_B, \diagmat{R}_B, \diagmat{Q}_B)
\end{displaymath}
be  SSS representations for the matrices $A$ and $B$.  Then,

\begin{enumerate}[label={(\roman*)}]
    \item there exists a SSS representation for $A^{-1} = \SSS(\diagmat{U}_{A^{-1}},\diagmat{W}_{A^{-1}},\diagmat{V}_{A^{-1}}, \diagmat{D}_{A^{-1}}, \diagmat{P}_{A^{-1}}, \diagmat{R}_{A^{-1}}, \diagmat{Q}_{A^{-1}})$  with dimensions:
    \begin{displaymath}
    r^h_{A,i} =  r^h_{A^{-1},i} \quad\mbox{and}\quad r^g_{A,i} =  r^g_{A^{-1},i}
    \end{displaymath}
    for $i=1,2,\ldots,n-1$.
    \item there exists a SSS representation for $C=A+B= \SSS(\diagmat{U}_C,\diagmat{W}_C,\diagmat{V}_C, \diagmat{D}_C, \diagmat{P}_C, \diagmat{R}_C, \diagmat{Q}_C)$ with dimensions:
      \begin{displaymath}
    r^h_{C,i} \leq   r^h_{A,i} + r^h_{B,i} \quad\mbox{and}\quad  r^g_{C,i} \leq   r^g_{A,i} + r^g_{B,i}
    \end{displaymath}
    for $i=1,2,\ldots,n-1$.
    \item there exists a SSS representation  for $C= AB =  \SSS(\diagmat{U}_C,\diagmat{W}_C,\diagmat{V}_C, \diagmat{D}_C, \diagmat{P}_C, \diagmat{R}_C, \diagmat{Q}_C)$ with  dimensions:
      \begin{displaymath}
    r^h_{C,i} \leq   r^h_{A,i} + r^h_{B,i} \quad\mbox{and}\quad  r^g_{C,i} \leq   r^g_{A,i} + r^g_{B,i}
    \end{displaymath}
    for $i=1,2,\ldots,n-1$.
\end{enumerate}
\end{proposition}

\begin{proof}
    By Theorem~\ref{thm:minimalSSS}, it is sufficient to bound the Hankel block ranks associated with the subgraph $\subgraph_i = \{1,\ldots,i\}$. Statement (i) follows from the Hankel block property $\rank (A^{-1})_{\complementgraph,\subgraph} = \rank A_{\complementgraph,\subgraph}$ (see Proposition~\ref{prop:rankinvariance}). \ec Statement (ii) follows again from Theorem~\ref{thm:minimalSSS} and the fact that
	\begin{displaymath}  \rank (A+B)_{\complementgraph,\subgraph}  \leq \rank A_{\complementgraph,\subgraph} + \rank B_{\complementgraph,\subgraph}. \ec
	\end{displaymath}
    Finally, for statement (iii) we  note that
	\begin{displaymath}
	\rank (AB)_{\complementgraph,\subgraph} = \rank (A_{\complementgraph,\complementgraph}B_{\complementgraph,\subgraph} + A_{\complementgraph,\subgraph}B_{\subgraph,\subgraph}) \leq \rank B_{\complementgraph,\subgraph} + \rank A_{\complementgraph,\subgraph}. \ec
	\end{displaymath} 
\end{proof}

SSS representations provide a framework for fast inversion of GIRS matrices on the line graph. The fast inversion algorithm can be derived from  \eqref{eq:toplevelSSS} and the structure of the signal flow diagram in Figure~\ref{subfig:signalflowSSS}. By merging the nodes in signal flow diagram, we may introduce the vectors $\xi = (\xi_i)^{N}_{i=1}$ and $\beta = (\beta_i)^{N}_{i=1}$ where
\begin{displaymath}
\xi_i = \begin{bmatrix} g_i \\ h_i \\ x_i \end{bmatrix}, \quad \beta_i = \begin{bmatrix} 0\\0 \\ b_i \end{bmatrix}.\ec
\end{displaymath}
Through this re-ordering,  \eqref{eq:toplevelSSS} can be re-expressed as
\begin{equation}
    \left( \boldsymbol{\Sigma} + \boldsymbol{\Theta} Z + \boldsymbol{\Gamma} Z^T  \right) \xi = \beta, \ec  \label{eq:SSSreordered}
\end{equation}
where
\begin{displaymath}
 \boldsymbol{\Sigma} =  \diag \left\lbrace \begin{bmatrix} 
 I & 0 & -V^T_i \\
 0 & I & -Q^T_i \\
 0 & 0 & D_i \end{bmatrix}
 \right\rbrace^{n}_{i=1}, \: \boldsymbol{\Theta} = \diag \left\lbrace \begin{bmatrix} 
 0 & 0 & 0 \\
 0 & -R_i & 0 \\
 0 & P_i & 0 \end{bmatrix}
 \right\rbrace^{n}_{i=1}   , \: \boldsymbol{\Gamma} =   \diag \left\lbrace \begin{bmatrix} 
 -W_i & 0 & 0 \\
 0 & 0 & 0 \\
 U_i & 0 & 0 \end{bmatrix}
 \right\rbrace^{n}_{i=1}. 
\end{displaymath}
With this re-ordering, it becomes evident  that the (block) adjacency graph of the lifted sparse system, as described by \eqref{eq:SSSreordered}, is the line graph, the same graph $\graph$ of the GIRS matrix (see figures \ref{subfig:SSS_matrix}  and \ref{subfig:SSS_reordered}). Hence, the complexity of the solver will be equivalent to doing Gaussian elimination on the graph $\graph$, which is dictated by the fastest elimination order of the graph. Since the line graph, can be eliminated in linear time, with the SSS representation, a GIRS-$r$ matrix on the line graph can be solved in $ \mathcal{O}(r^2 N)$ complexity.  \ec

\subsection{Dewilde-van der Veen representations} \label{subsec:DV}

Since SSS representation completely characterizes GIRS matrices on the line graph (see Theorems \ref{thm:minimalSSS} and \ref{thm:SSStoGIRS}), it is natural to seek SSS-like representations for more general GIRS matrices. To this end, we introduce two candidate representations which we shall call Dewilde-van der Veen (DV) and $\graph$-semi-separable ($\graph$-SS) matrices (to be introduced later in Section \ref{sec:edv}). Both of these representations, if constructed, would give rise to fast linear solves in time complexity commensurate to doing sparse Gaussian elimination on the underlying graph $\graph$.  \ec 

The key idea behind SSS matrices is to introduce a flow on the nodes of the line graph using the state variables $h_i$ and $g_i$. In the case of SSS, this flow is decomposed into two explicit ones (upstream/casual and downstream/anti-causal). For a line graph, it is pretty straightforward  how these flows should be defined, but this is less so for arbitrary graphs. Nevertheless, we would like to extend SSS to arbitrary graphs, and it is possible to also consider implicit state-space equations which arise from a single implicit flow. For line graphs, on which SSS is defined, these implicit equations take the form
\begin{eqnarray*}
g_{k}  & = &  V^T_{k} x_{k} + W^-_{k} g_{k-1}  + W^+_{k} g_{k+1} \\ 
b_k  & = &   D_k x_k + U^{-}_{k} g_{k-1} + U^{+}_{k}g_{k+1} 
\end{eqnarray*}
for $k=2,\ldots,n$ with boundary conditions
\begin{displaymath}
\begin{array}{lcl}
g_{1}  & = &  V^T_{1} x_{1}  + W^+_{1} g_{2} \\ 
b_1  & = &   D_1 x_1  + U^{+}_{1} g_2
\end{array},\qquad \begin{array}{lcl}
g_{n}  & = &  V^T_{n} x_{n} + W^-_{n} g_{n-1}  \\ 
b_n  & = &   D_n x_n + U^{-}_{n} g_{n-1}  
\end{array}.
\end{displaymath}
Implicit representation such as these can be easily generalized to general graphs, which gives rise to the following definition. 

\begin{definition}[Dewilde-van der Veen Representation]
Let $(A,\graph)$ be a graph-partitioned matrix and let $\mathbb{N}_{\graph}(i)$ denote all nodes $j\in\nodeset_\graph$ adjacent to the node $i\in\nodeset_\graph$. A \textit{Dewilde-van der Veen (DV) Representation} for $A$ is a collection of matrices such that the \textit{implicit state-space equations} 
\begin{subequations}\label{eq:idv}
\begin{align}
 g_i = V^T_i x_i + \sum_{j\in \mathbb{N}_{\graph}(i)} W_{i,j} g_j \label{eq:idv_1}\ec \\
 b_i = D_i x_i + \sum_{j\in \mathbb{N}_{\graph}(i)} U_{i,j} g_j \label{eq:idv_2}\ec 
\end{align}
\end{subequations}
for $i\in \nodeset_{\graph}$ are uniquely solvable and consistent with \eqref{eq:block_matrix_multiplication}.\ec 
\end{definition} 

Similar to \eqref{eq:toplevelSSS} for SSS, equation \eqref{eq:idv} can be expressed more compactly in matrix notation. Denote  $\diagmat{V} := \diag \{V_i\}_{i\in\nodeset_\graph}$ and $\diagmat{D} := \diag  \{ D_i \}_{i\in \nodeset_\graph}$. Furthermore denote $\diagmat{W} := \diag \{W_{i,j}\}_{(i,j)\in\edgeset_\graph}$ and $\diagmat{U} := \diag \{U_{i,j}\}_{(i,j)\in\edgeset_\graph}$. We may introduce a matrix-valued operator $Z_\graph[\cdot]$ such that the $(i,j)$-th block entry $Z_\graph[\diagmat{W}]$ is given by
\begin{displaymath}
\left[ Z_\graph[\diagmat{W}] \right]_{i,j} = \begin{cases}
W_{i,j} & \mbox{if }j\in \mathbb{N}_{\graph}(i)\\
0, & \mbox{if }j\notin \mathbb{N}_{\graph}(i). \ec\end{cases} 
\end{displaymath}
With the help of $Z_\graph[\cdot]$, \eqref{eq:idv_1} can compactly be expressed as $\diagmat{D}^T x + Z_\graph[\diagmat{U}] g = b$, and likewise, $g = \diagmat{V}^Tx + Z_\graph[\diagmat{W}]g$ is a compact expression for \eqref{eq:idv_2}. Overall,  \eqref{eq:idv} is placed on the same footing as \eqref{eq:toplevelSSS} does for SSS:
\begin{equation}\label{eq:idv_matrix_form}
    \twobytwo{I - Z_\graph[\diagmat{W}]}{\diagmat{V}^T}{Z_\graph[\diagmat{U}]}{\diagmat{D}}\twobyone{g}{x} = \twobyone{0}{b}.\ec 
\end{equation}
\begin{figure}
    \centering
\begin{subfigure}[b]{0.3\textwidth}    
\begin{center}
\begin{tikzpicture}
    \node[circle,draw=black, fill=white, inner sep=0pt,minimum size=20pt] (a) at (0,0)  {$1$};
    \node[circle,draw=black, fill=white, inner sep=0pt,minimum size=20pt] (b) at (1,0)  {$2$};
    \node[circle,draw=black, fill=white, inner sep=0pt,minimum size=20pt] (c) at (2,0)  {$3$};
    \node[circle,draw=black, fill=white, inner sep=0pt,minimum size=20pt] (d) at (0,1)  {$4$};
    \node[circle,draw=black, fill=white, inner sep=0pt,minimum size=20pt] (e) at (1,1)  {$5$};
    \node[circle,draw=black, fill=white, inner sep=0pt,minimum size=20pt] (f) at (2,1)  {$6$};
    \node[circle,draw=black, fill=white, inner sep=0pt,minimum size=20pt] (g) at (0,2) {$7$};
    \node[circle,draw=black, fill=white, inner sep=0pt,minimum size=20pt] (h) at (1,2) {$8$};
    \node[circle,draw=black, fill=white, inner sep=0pt,minimum size=20pt] (i) at (2,2)  {$9$};
    \draw (a) edge[-] (b);
    \draw (b) edge[-] (c);
    \draw (a) edge[-] (d);
    \draw (b) edge[-] (e);
    \draw (c) edge[-] (f);
    \draw (d) edge[-] (e);
    \draw (e) edge[-] (f);
    \draw (d) edge[-] (g);
    \draw (e) edge[-] (h);
    \draw (f) edge[-] (i);
    \draw (g) edge[-] (h);
    \draw (h) edge[-] (i);
\end{tikzpicture}
\end{center}
\caption{\ec} \label{subfig:3x3}
\end{subfigure}

    \begin{subfigure}[b]{0.9\textwidth}
        \begin{displaymath}
        \resizebox{0.99\textwidth}{!}{$\left[\begin{array}{ccc;{2pt/2pt}ccc;{2pt/2pt}ccc|ccc;{2pt/2pt}ccc;{2pt/2pt}ccc}
        I & -W_{1,2} &  & -W_{1,4} &  &  &  &  &  &              V^T_{1} &  &  &  &  &  &  &  &  \\ 
        -W_{2,1} & I & -W_{2,3} &  & -W_{1,5} &  &  &  &  &         & V^T_{2} &  &  &  &  &  &  &  \\
         & -W_{3,2} & I &  &  & -W_{3,6} &  &  &  &               &  & V^T_{3} &  &  &  &  &  &  \\  \hdashline[2pt/2pt]
        -W_{4,1} &  &  & I & -W_{4,5} &  & -W_{4,7} &  &  &         &  &  & V^T_{4} &  &  &  &  &  \\
         & -W_{5,2} &  & -W_{5,4} & I & -W_{5,6} &  & -W_{5,8} &  &   &  &  &  & V^T_{5} &  &  &  &  \\
         &  & -W_{6,3} &  & -W_{6,5} & I &  &  & -W_{6,9}  &         &  &  &  &  & V^T_{6} &  &  &  \\  \hdashline[2pt/2pt]
         &  &  & -W_{7,4} &  &  & I & -W_{7,8} &   &              &  &  &  &  &  & V^T_{7} &  &  \\
         &  &  &  & -W_{8,5} &  & -W_{8,7} & I & -W_{8,9} &         &  &  &  &  &  &  & V^T_{8} &  \\
         &  &  &  &  & -W_{9,6} &  & -W_{9,8}  & I  &             &  &  &  &  &  &  &   & V^T_{9} \\ \hline
         & U_{1,2} &  & U_{1,4} &  &  &  &  &    &            D_{1} &  &  &  &  &  &  &  &  \\
        U_{2,1} &  & U_{2,3} &  & U_{1,5} &  &  &  &  &         & D_{2} &  &  &  &  &  &  &  \\
         & U_{3,2} &  &  &  & U_{3,6} &  &  &   &              &  & D_{3} &  &  &  &  &  &  \\ \hdashline[2pt/2pt]
        U_{4,1} &  &  &  & U_{4,5} &  & U_{4,7} &  &  &         &  &  & D_{4} &  &  &  &  &  \\
         & U_{5,2} &  & U_{5,4} &  & U_{5,6} &  & U_{5,8} &  &   &  &  &  & D_{5} &  &  &  &  \\
         &  & U_{6,3} &  & U_{6,5} &  &  &  & U_{6,9} &         &  &  &  &  & D_{6} &  &  &  \\ \hdashline[2pt/2pt]
         &  &  & U_{7,4} &  &  &  & U_{7,8} &  &               &  &  &  &  &  & D_{7} &  &  \\
         &  &  &  & U_{8,5} &  & U_{8,7} &  & U_{8,9} &         &  &  &  &  &  &  & D_{8} &  \\
         &  &  &  &  & U_{9,6} &  & U_{9,8}  &  &              &  &  &  &  &  &  &   & D_{9} 
        \end{array} \right]$}
        \end{displaymath}
        \caption{\ec} \label{subfig:DV_matrix}
    \end{subfigure}
 
    \begin{subfigure}[b]{0.9\textwidth}
        \begin{displaymath}
        \resizebox{0.99\textwidth}{!}{$\left[\begin{array}{cc;{2pt/2pt}cc;{2pt/2pt}cc|cc;{2pt/2pt}cc;{2pt/2pt}cc|cc;{2pt/2pt}cc;{2pt/2pt}cc}
        I        &   V^T_1 &  -W_{1,2} &   &   &  & -W_{1,4} &   &  &  &  &  &  &  &  &  &  &  \\
                &  D_1    &     U_{1,2} &   &    &  &  U_{1,4} &   &  &  &  &  &  &  &  &  &  &  \\ \hdashline[2pt/2pt]
        -W_{2,1} &        & I &  V^T_2 &  -W_{2,3} &    &  &  &  -W_{2,5} &   &  &  &  &  &  &  &  & \\
        U_{2,1}  &        &  &  D_2 &  U_{2,3} &   &    &  &  -U_{2,5} &   &  &  &  &  &  &  &  &  \\ \hdashline[2pt/2pt]
         &  & -W_{3,2} &   & I &  V^T_3  &   &  &   &   &  -W_{3,6} &    &  &   &   &    &   &   \\
         &  & U_{3,2} &   &  &  D_3 &   &  &   &   & U_{3,6} &   &  &   &   &    &   &   \\  \hline
        
         -W_{4,1} &   &   &    &  &   & I &  V^T_4  & -W_{4,5} &   &  &  & -W_{4,7} &    &   &    &   &   \\ 
        U_{4,1} &   &   &    &  &   &  &  D_4  & U_{4,5} &   &  &  & U_{4,7} &    &   &    &   &   \\  \hdashline[2pt/2pt]
         &  &   -W_{5,2} &     &  &   &  -W_{5,4} &   & I &  V^T_5 &  -W_{5,6} &   &   &   &   -W_{5,8} &   &  &  \\ 
         &  &    U_{5,2} &     &  &   & U_{5,4} &     &  &  D_5   &   U_{5,6} &   &   &   &  U_{5,8} &   &  &  \\  \hdashline[2pt/2pt]
         &  &   &    &  -W_{6,3} &   &  &  & -W_{6,5} &   & I &  V^T_6 &  &  &   &   &  -W_{6,9} &     \\ 
         &  &   &    &   U_{6,3} &   &  &  &  U_{6,5} &    &  &  D_6   &  &   &   &    &  U_{6,9} &    \\  \hline
         &  &   &    &  &   &  -W_{7,4} &   &  &  &    &   & I &  V^T_7 &  -W_{7,8} &   &   &      \\ 
         &  &   &    &  &   &  U_{7,4} &    &  &   &  &     &   &  D_7    &  U_{7,8} &   &   &    \\  \hdashline[2pt/2pt]
         &  &   &    &  &   &  &   &  -W_{8,5} &   &  &  &  -W_{8,7} &  & I &  V^T_8  &   -W_{8,9} &    \\ 
         &  &   &    &  &   &  &   &  U_{8,5}  &   &  &  &  U_{8,7} &   &  &  D_8  &  U_{8,9} &    \\  \hdashline[2pt/2pt]
         &  &   &    &  &   &  &   &   &  & -W_{9,6} &   &   &   &   -W_{9,8} &   & I &  V^T_9 \\ 
         &  &   &    &  &   &  &   &   &  & U_{9,6} &   &   &   &    U_{9,8} &   &  &  D_9
        \end{array} \right] $}
        \end{displaymath}
        \caption{\ec} \label{subfig:DV_reordered}
    \end{subfigure}
    
    \caption{The matrix form of the DV representation \eqref{eq:idv_matrix_form}, shown in (\subref{subfig:DV_matrix}), for the $3\times 3$ mesh graph (\subref{subfig:3x3}). Upon reordering \eqref{eq:idv_matrix_form} to \eqref{eq:solveDV}, the (block) adjacency graph of the reordered system \eqref{eq:solveDV} coincides with the graph $\graph$ (here, the $3\times 3$ mesh graph); see (\subref{subfig:DV_reordered}). This allows for a fast solver in the same time complexity as performing Gaussian elimination on $\graph$.
    }\label{fig:DV}
\end{figure}
The block sparsity pattern is shown in Figure \ref{fig:DV}. This gives rise to the Schur complement expression
\begin{equation}\label{eq:idv_schur_complement}
    A = \diagmat{D} + Z_\graph[\diagmat{U}](I-Z_\graph[\diagmat{W}])^{-1}\diagmat{V}^T,  
\end{equation}
which is the analogue of \eqref{eq:diagrepSSS} for DV representations. We may denote DV representations concisely by 
\begin{displaymath}
A = \IDV(\graph; \diagmat{V}, \diagmat{W},  \diagmat{D}, \diagmat{U}).
\end{displaymath}
The main differences are the replacement of the downshift operator $Z$ with a more general matrix valued operator $Z_\graph[\cdot]$, and the use of a single state variable $g_i$, as opposed to two: $g_i$ and $h_i$. Nevertheless, a SSS representation is also a DV representation in the broadest sense. Indeed, we may simply merge the two state variables into one, leading to the supposedly ``implicit" equations:
\begin{eqnarray*}
\begin{bmatrix}
g_{k}
\\
h_k
\end{bmatrix}  & = &  \begin{bmatrix} V^T_{k} \\  Q^T_k  
\end{bmatrix}   x_{k} + \begin{bmatrix}   W_{k} & 0 \\  0 & 0 \end{bmatrix} \begin{bmatrix}
g_{k-1}
\\
h_{k-1}
\end{bmatrix}  +  \begin{bmatrix} 0 & 0
 \\ 0 &  R_{k} \end{bmatrix} \begin{bmatrix}  g_{k+1} \\  h_{k+1} \end{bmatrix} \\ 
b_k  & = &   D_k x_k + \begin{bmatrix}  U_{k} & 0 \end{bmatrix} \begin{bmatrix}
g_{k-1}
\\
h_{k-1}
\end{bmatrix}  + \begin{bmatrix}   0 &  P_{k} \end{bmatrix} \begin{bmatrix}
g_{k+1}
\\
h_{k+1}
\end{bmatrix}
\end{eqnarray*}\ec 
for $k=2,\ldots,n-1$ with boundary conditions
\begin{displaymath}
\begin{array}{lcl}
\begin{bmatrix}
g_{1}
\\
h_{1}
\end{bmatrix} & = &  \begin{bmatrix} V^T_{1} \\ 0  \end{bmatrix}  x_{1}  +  \begin{bmatrix} 0 & 0 \\
  0 & 0
\end{bmatrix}  \begin{bmatrix}
g_{2}
\\
h_{2}
\end{bmatrix} \\ 
b_1  & = &   D_1 x_1  + \begin{bmatrix}  0 & P_1  \end{bmatrix}  \begin{bmatrix} g_{2} \\ h_2 \end{bmatrix} 
\end{array},\qquad \begin{array}{lcl}
\begin{bmatrix}
g_{n}
\\
h_{n}
\end{bmatrix}   & = &  \begin{bmatrix} 0   \\ Q^T_{n}  \end{bmatrix}  x_{n}  + \begin{bmatrix} 0 & 0 \\
  0 & 0
\end{bmatrix}  \begin{bmatrix}
g_{n-1}
\\
h_{n-1}
\end{bmatrix}   \\ 
b_n  & = &   D_n x_n + \begin{bmatrix}  U_{n} & 0 \end{bmatrix}  \begin{bmatrix} g_{n-1} \\ h_{n-1} \end{bmatrix}   
\end{array}. \ec 
\end{displaymath}
DV representations present a generalization of SSS representations for more general graph partitioned matrices. Given this fact, we are interested in addressing the following questions:

\begin{enumerate}
    \item When do graph partitioned matrices have (efficient) DV representation and how is this related to GIRS property?
    \item How can (efficient) DV representations be constructed?
    \item To what extent are the properties of SSS representations inherited by DV representations?
\end{enumerate}

For the first two questions, we only know incomplete answers at this stage. In general, it is not yet clear how DV representation can be found (let alone finding minimal or compact ones, see remark~\ref{rmk:optimality_DV}), except for some special cases. For example, we know that in the case of sparse matrices, the problem of construction is relatively straightforward as the following example shows.

\begin{remark}\label{rmk:optimality_DV}
In the case of SSS representations, we can produce a single representation for which the dimensions of each state space variables $g_i$ and $h_i$ are as small as they can be: that is, the SSS representation is \textit{uniformly minimal} in the sense that $r_i^g \le \min \hat{r}_i^g$ and $r_i^h \le \min \hat{r}_i^h$ for all $I \in \{1,\ldots,n-1\}$, where the minimum is taken over all SSS representations. The matrices involved in a DV representation for a graph partitioned matrix $(A,\graph)$ are of the dimensions $V_i\in\mathbb{R}^{N_i \times r_i}$, $W_{i,j} \in \mathbb{R}^{r_i \times r_j}$, $D_i \in \mathbb{R}^{N_i \times N_i}$, and $U_{i,j}\in\mathbb{R}^{ N_i \times r_j}$.   Hence, we say that a DV representation for a matrix $A$ is \textit{uniformly minimal} if $r_i \le \min \hat{r}_i$, where the minimum is taken over all DV representations of $A$. A priori, uniformly minimal DV representations may not exist in which case we will have to settle for either \textit{minimal} DV representations---which minimize the total size $\sum_{i\in\nodeset_\graph} r_i$ of the representation---or, even more loosely, merely compact DV representation---for which $\sum_{i\in\nodeset_\graph} r_i \ll N^2$.
\end{remark}

\begin{example}\label{prop:sparse_idv}
Let $A$ be a sparse matrix with adjacency graph $\graph$. Then $D_i := A_{ii}$, $U_{i,j} = A_{ij}$ for $i\ne j$ and $A_{ij} \ne 0$, $V_i = 1$, and $W_{i,j} = 0$ gives a DV representation of $A$.\ethancheck
\end{example}

Given that our motivation for considering GIRS matrices was that they precisely characterized the rank-structure properties of sparse matrices, the existence of compact DV representations for sparse matrices lends credence to the idea that DV representations may hold promise for representing general GIRS matrices. In Theorems \ref{thm:minimalSSS} and \ref{thm:SSStoGIRS} we showed how SSS matrices are closely related GIRS matrices on the line graph and vice versa. For general DV representation, we only have the following result.\ethancheck

\begin{proposition}\label{prop:idv_girs}
Suppose that a graph-partitioned matrix $(A,\graph)$ possesses a DV representation $A=\IDV(\graph; \diagmat{V}, \diagmat{W},  \diagmat{D}, \diagmat{U})$  with dimensions $r_i$ for $i \in \nodeset_{\graph}$. Then $(A,\graph)$ is GIRS-$2r$ where $r = \max_{i \in \nodeset_{\graph}}  r_{i}$.
\end{proposition}

\begin{proof}
Let $\subgraph\subset \graph$ be a sub-graph and observe that
\begin{eqnarray*}
\rank \left( Z_\graph[\diagmat{W}] \right)_{\complementgraph,\subgraph} & \leq &  \sum_{\stackrel{(i,j)\in \edgeset_{\graph}}{i \in \nodeset_{\complementgraph},j \in \nodeset_{\subgraph}}} \rank  W_{i,j}  \\ 
& \leq &  \rho(\subgraph) \max_{\stackrel{(i,j)\in \edgeset_{\graph}}{i \in \nodeset_{\complementgraph},j \in \nodeset_{\subgraph}}} 
\left(\rank  W_{i,j} \right) = \rho(\subgraph) r. \ethancheck 
\end{eqnarray*}
It hence follows that $I- Z_\graph[\diagmat{W}]$ is also GIRS-$r$. By Proposition~\ref{prop:girs_algebra}-(i), $(I- Z_\graph[\diagmat{W}])^{-1}$  is also GIRS-$r$. The same implies for $(I- Z_\graph[\diagmat{W}])^{-1}\diagmat{V}^T$.  With a similar argument as for $Z_\graph[\diagmat{W}]$, note that
\begin{displaymath}
\rank \left( Z_\graph[\diagmat{U}] \right)_{\complementgraph,\subgraph} \leq \rho(\subgraph) r, \ethancheck 
\end{displaymath}
which shows that $Z_\graph[\diagmat{U}]$ is also GIRS-$r$. By Proposition~\ref{prop:girs_algebra}-(ii), the product of two  GIRS-$r$ matrices is GIRS-$2r$, hence $Z_\graph[\diagmat{U}](I-Z_\graph[\diagmat{W}])^{-1}\diagmat{V}^T$
is GIRS-$2r$, \ethancheck and therefore
\begin{displaymath}
A = \diagmat{D}+ Z_\graph[\diagmat{U}](I-Z_\graph[\diagmat{W}])^{-1}\diagmat{V}^T
\end{displaymath}
is GIRS-$2r$.\ethancheck 
\end{proof}

Proposition~\ref{prop:idv_girs} shows that a compact DV representation with respect to some partitioning $(A,\graph)$ implies that $A$ is GIRS on the corresponding $\graph$. The converse of this statement is still an open question and will be discussed in Section~\ref{sec:overview}. In the case of SSS, we relied upon the construction algorithm to prove the converse statement in  Theorem~\ref{thm:minimalSSS}.  However, the construction for general DV representations with arbitrary graphs appear to be nontrivial. Constructing the optimal representation requires determining the minimal size $r_i$ of the state-space variables $g_i$ in addition to determination of the weights $V_i$, $W_{i,j}$, $U_{i,j}$  and $D_i$. \ethancheck

Despite of the difficulties in construction, the algebraic properties of DV representation very nicely generalize those of SSS. The following proposition is the analogue of Proposition~\ref{prop:SSSalgebra} for general DV representations.\ec

\begin{proposition}[DV algebra]\label{prop:idv_algebra}
Let
\begin{displaymath}
A =\IDV(\graph; \diagmat{V}_A, \diagmat{W}_A,  \diagmat{D}_A, \diagmat{U}_A),\quad B = \IDV(\graph; \diagmat{V}_B, \diagmat{W}_B,  \diagmat{D}_B, \diagmat{U}_B)  .
\end{displaymath}
Then the following hold:

\begin{enumerate}[label={(\roman*)}]
    \item There exists a DV representation for $A^{-1} = 
    \IDV(\graph; \diagmat{V}_{A^{-1}}, \diagmat{W}_{A^{-1}},  \diagmat{D}_{A^{-1}}, \diagmat{U}_{A^{-1}})$ with dimensions
    
    \begin{displaymath} 
    r_{A^{-1},i} =   r_{A,i  }
    \end{displaymath}
    
    for $i\in \nodeset_\graph$.
    
    \item There exists a DV representation for $C = A+B =\IDV(\graph; \diagmat{V}_C, \diagmat{W}_C,  \diagmat{D}_C, \diagmat{U}_C) $ with the dimensions
    
    \begin{displaymath} 
    r_{C,i} \leq   r_{A,i} + r_{B,i} 
    \end{displaymath}
    
    for $i\in \nodeset_\graph$.
    
    \item There exists a DV representation for $C = AB = \IDV(\graph; \diagmat{V}_C, \diagmat{W}_C,  \diagmat{D}_C, \diagmat{U}_C)$ with the dimensions
    \begin{displaymath} 
    r_{C,i} \leq   r_{A,i} + r_{B,i}
    \end{displaymath}
    
    for $i\in \nodeset_\graph$.
\end{enumerate}
\end{proposition}

\begin{proof}
    Statement (i) is validated by applying the Sherman-Morrison-Woodbury identity
    \begin{equation}\label{eq:SMW}
        (B+UCV)^{-1} = B^{-1} - B^{-1} U (C^{-1} + V B^{-1} U)^{-1} V B^{-1}\ethancheck 
    \end{equation}
    to the diagonal representation \eqref{eq:idv_schur_complement}. 
    This gives
    \begin{eqnarray*}
     A^{-1}  & = & \left( \diagmat{D} + Z_\graph[\diagmat{U}](I-Z_\graph[\diagmat{W}])^{-1}\diagmat{V}^T \right)^{-1}  \\
       & =  &  \diagmat{D}^{-1} +   \diagmat{D}^{-1} Z_\graph[\diagmat{U}] \left( I-Z_\graph[\diagmat{W}]  +  \diagmat{V}^T \diagmat{D}^{-1}  Z_\graph[\diagmat{U}] \right)^{-1}  \diagmat{V}^T \diagmat{D}^{-1},\ethancheck 
    \end{eqnarray*}
    which leads to
    \begin{equation} 
    A =  \diagmat{\hat{D}} +    Z_\graph[\diagmat{\hat{U}}] \left( I-Z_\graph[\diagmat{\hat{W}}] \right)^{-1}  \diagmat{\hat{V}}^T, \label{eq:inverseDV}\ethancheck 
    \end{equation}
    where
    \begin{gather*}
    \diagmat{\hat{D}} = \diag  \{ D^{-1}_i \}_{i\in \nodeset_\graph}, \quad  \diagmat{\hat{U}} = \diag \{ D^{-1}_i U_{i,j} \}_{(i,j)\in\edgeset_\graph},\\
    \diagmat{\hat{W}} = \diag \{ W_{i,j} + V^T_i D^{-1}_i U_{i,j} \}_{(i,j)\in\edgeset_\graph}, \quad \diagmat{\hat{V}} = \diag  \{ D^{-T}_i V_i \}_{i\in \nodeset_\graph}.\ethancheck 
    \end{gather*}
   For statement  (ii), we may simply set 
    \begin{gather*}
	  V_{C,i} =  \begin{bmatrix} 
		  V_{A,i} & V_{B,i} \end{bmatrix}, \quad W_{C,i,j} = \begin{bmatrix} W_{A,i} & \\
		  &    W_{B,i,j} \end{bmatrix}, \\
		  D_{C,i} = D_{A,i} + D_{B,i}, \quad  U_{C,i,j} = \begin{bmatrix} U_{A,i,j} & U_{B,i,j} \end{bmatrix}.\ethancheck 
	\end{gather*}
	to obtain and DV representation for $C=A+B$.
	Similarly, for $C=AB$, write $b = AB x = A z$. We have the set of equations
	\begin{align}
    g_i = V^T_{B,i} x_i + \sum_{j\in \mathbb{N}_{\graph}(i)} W_{B,i,j} g_j  \\
    z_i = D_{B,i} x_i + \sum_{j\in \mathbb{N}_{\graph}(i)} U_{B,i,j} g_j  \label{eq:thisone}
 \end{align}
	and 
	\begin{align}
    h_i = V^T_{A,i} z_i + \sum_{j\in \mathbb{N}_{\graph}(i)} W_{A,j} h_j \label{eq:andthisone} \\
    b_i = D_{A,i} z_i + \sum_{j\in \mathbb{N}_{\graph}(i)} U_{A,i,j} h_j \label{eq:donforgetthisone}. 
 \end{align}
By substitution, we can merge these two sets of equations into
\begin{align}
    \begin{bmatrix} g_i \\ h_i \end{bmatrix} = \begin{bmatrix} V^T_{B,i} \\   V^T_{A,i} D_{B,i}\end{bmatrix} x_i + \sum_{j\in \mathbb{N}_{\graph}(i)} 
    \begin{bmatrix} W_{B,i,j} & \\ V^T_{A,i} U_{B,i,j} &    W_{A,i,j} \end{bmatrix}   \begin{bmatrix} g_j \\ h_j \end{bmatrix}  \\
    z_i = D_{A,i}D_{B,i} x_i + \sum_{j\in \mathbb{N}_{\graph}(i)} \begin{bmatrix} D_{A,i} U_{B,i,j} & U_{A,i,j} \end{bmatrix}  \begin{bmatrix} g_j \\ h_j \end{bmatrix}   
 \end{align}	
which shows that
	\begin{gather*}
	  V_{C,i} =  \begin{bmatrix} 
		  V_{B,i} & D^T_{B,i} V_{A,i} \end{bmatrix}, \quad W_{C,i,j} = \begin{bmatrix} W_{B,i,j} & \\
		  V^T_{A,i} U_{B,i,j} &    W_{A,i,j} \end{bmatrix},\\
		  D_{C,i} = D_{A,i}D_{B,i}, \quad  U_{C,i,j} = \begin{bmatrix} D_{A,i} U_{B,i,j} & U_{A,i,j} \end{bmatrix}.
	\end{gather*}
	is a DV representation for $C=AB$, hence proving (iii). \nithincheck
\end{proof}

Once a DV representation of $A$ has been computed, solving $Ax = b$ can be done in a similar way to SSS. Introducing the vectors $\xi = (\xi_i)^{N}_{i}$ and $\beta = (\beta_i)^{N}_{i}$ with
\begin{displaymath}
\xi_i = \begin{bmatrix}  g_i \\ x_i \end{bmatrix}, \quad \beta_i = \begin{bmatrix} 0 \\ b_i \end{bmatrix},
\end{displaymath}
we obtain an equation similar to \eqref{eq:SSSreordered}:
\begin{equation}
    \left( \boldsymbol{\Sigma} +  Z_\graph[ \boldsymbol{\Theta} ] \right) \xi = \beta,      \label{eq:solveDV}
\end{equation}
where
\begin{displaymath}
 \boldsymbol{\Sigma} =  \diag \left\lbrace \begin{bmatrix} 
 I &  V^T_i \\
 0 &  D_i \end{bmatrix}
 \right\rbrace^{n}_{i=1}, \quad \boldsymbol{\Theta} =  \diag \left\lbrace \begin{bmatrix} 
 -W_{i,j} &  0 \\
 U_{i,j} &  0 \end{bmatrix}
 \right\rbrace_{{(i,j)\in\edgeset_\graph}}.\ethancheck 
\end{displaymath} 
Interesting, with regard to the complexity required to evaluate $A^{-1}b$ or $Ax$, both are equally expensive operations for DV representations---both requiring time necessary to do sparse Gaussian elimination on $\graph$.

\begin{proposition} \label{prop:DVsolve}
Let $A = \IDV(\graph; \diagmat{V}, \diagmat{W},  \diagmat{D}, \diagmat{U}) $ be an DV representation for a graph-partitioned matrix $(A,\graph)$. The products $Ax$ and $A^{-1}b$ can be evaluated in the time complexity of sparse Gaussian elimination on $\graph$.\ethancheck 
\end{proposition}

\begin{proof}
    For  matrix-vector multiplication $Ax$, this result immediately  follows from \eqref{eq:idv_schur_complement}, which involves computing a product of the form $(Z_\graph[\diagmat{W}]-I)^{-1}\eta$ which can be done by performing sparse Gaussian elimination on $\graph$. For the multiplication by the inverse $A^{-1}b$, we perform sparse Gaussian elimination on \eqref{eq:solveDV}. Alternatively, one can also derive the result from \eqref{eq:inverseDV}.  \ethancheck 
\end{proof}

}

\subsection{\texorpdfstring{$\graph$-}{Graph }semi-separable representations}\label{sec:edv}

\nousedv{In the case of GIRS matrix on the line graph, the formulation of an effective algebraic representation is well-addressed by the theory of sequentially semi-separable (SSS) matrices \cite{chandrasekaran2002fast,chandrasekaran2005some}. A summary is provided in \cite[Sec.~3]{chandrasekaran2018fast}. In an expanded preprint of this article, we discuss the SSS representation in the context of the GIRS property, establishing a one-to-one correspondence between GIRS-$c$ matrices and SSS matrices with off-diagonal blocks of rank uniformly bounded by $c$ \cite[Sec.~2.2]{chandrasekaran2019graph}. We summarize the main results regarding SSS matrices within the GIRS language in the following theorem.

\begin{theorem}\label{thm:sss}
Let $(A,\mathbb{L})$ be a graph-partitioned matrix with respect to the $\mathbb{L}$ a line (path) graph on $n$ nodes. Then $A$ is described by an SSS representation consisting of matrices $U_i \in \real^{N_i \times \rg_i}$, $W_i \in \real^{\rg_{i-1}\times \rg_i}$, $V_i \in \real^{N_i\times \rg_{i-1}}$, $D_i \in \real^{N_i\times N_i}$, $P_i \in \real^{N_i \times \rh_{i-1}}$, $R_i\in \real^{\rh_i\times \rh_{i-1}}$, and $Q_i \in \real^{N_i\times \rh_i}$ for $1\le i\le n$ such that
\begin{equation}\label{eq:sss}
    A = \resizebox{0.82\textwidth}{!}{$\begin{bmatrix}
D_1 &   U_1 V_2^T             &    U_1 W_2 V_3^T  &    U_1 W_2 W_3 V_3^T & \cdots &  U_1 W_2 W_3\cdots W_{n-1}V_n^T \\
\tilde{P}_{2} Q_1^T   & D_2  & U_2 V_3^T & U_2 W_3 V_4^T  & \cdots &          U_2 W_3 \cdots W_{n-1} V_n^T   \\
P_{3} R_{2} Q_1^T          &  P_{3} Q^T_2 & D_3 & U_3 V^T_4 & \cdots &   U_3W_4 \cdots W_{n-1}V_n^T   \\
P_{4} R_{3} R_{2} Q_1^T     &  P_{4} R_{3} Q^T_2  & P_{4} Q^T_3  & D_4  & \cdots &  U_4 W_5 \cdots W_{n-1}V_n^T  \\
\vdots   & \vdots & \vdots & \vdots & \ddots & \vdots  \\
P_{n} R_{n-1} \cdots R_{2} Q_1^T     &  P_{n} R_{n-1} \cdots R_{3}  Q^T_2         &  P_{n} R_{n-1} \cdots R_{4}  Q^T_3 & P_{n} R_{n-1} \cdots R_{5}  Q^T_4 & \cdots & D_n
\end{bmatrix}$}.
\end{equation}
The state dimensions $\rg_i$ and $\rh_i$ are given by the ranks of the Hankel blocks
\begin{equation*}
    \rg_i = \rank A_{\subgraph_i,\bar{\subgraph}_i}, \quad \rh_i = \rank A_{\bar{\subgraph}_i,\subgraph_i} \quad \textnormal{for } 1\le i \le n
\end{equation*}
where $\subgraph_i$ is the induced subgraph with nodes $\nodeset_{\subgraph_i} = \{1,2,\ldots,i\}$. In particular, if $(A,\mathbb{L})$ is a GIRS-$r$ matrix then the state dimensions $\rg_i$ and $\rh_i$ are bounded by $r$. In such case, given an SSS representation of $A$, there exist algorithms to multiply by $A$ in time $\mathcal{O}(nr^2)$ and to solve $Ax = b$ in time $\mathcal{O}(nr^3)$.
\end{theorem}}

A key ingredient to the SSS representation is the ability to compute the action of the matrix by a two-sweep algorithm, with one front-to-back sweep computing the action of the upper triangular portion of the matrix and a back-to-front sweep computing the action of the lower triangular portion. \usedv{Matrix vector multiplications with the Dewilde-van der Veen representation require solving a linear system of equations and don't have this attractive property of possessing an efficient, explicit multiplication algorithm. To compensate for this,}\nousedv{Mimicking this,} we induce an ordering of the nodes of the graph $\graph$ artificially and seek a representation in which information ``flows'' only from front-to-back (i.e.\ upstream) or back-to-front (i.e.\ downstream). This is leads to the following \usedv{different }generalization of the SSS representation. \ethancheck \nithincheck

\newcommand{\upstream}{{\rm u}}
\newcommand{\downstream}{{\rm d}}


\begin{definition}[$\graph$-semi-separable representations] \label{def:GSS}
Let $(A,\graph)$ be a graph-partitioned matrix, where $\graph$ admits a Hamiltonian path $\dirpath$ inducing a total order (denoted by $\prec$) on the node set  $\nodeset_\graph$:
\begin{equation*}
    u \prec v \iff \mbox{there exists a directed path from $u$ to $v$ along $\dirpath$}.
\end{equation*}
Let
$
\dirpath_{\upstream}(i) := \{ j \in \nodeset_\graph: (i,j) \in \edgeset_{\graph}\mbox{ and }j \prec i   \}
$ and  $
\dirpath_{\downstream}(i) :=  \{ j \in \nodeset_\graph: (i,j) \in \edgeset_{\graph} \mbox{ and } j \succ  i   \}
$. 
A \textit{$\graph$-semi-separable representation} for $A$ is a collection of matrices such that the \textit{state-space equations} 
\begin{subequations} \label{eq:gss}
\begin{align}
	g_i &= V_i^Tx_i + \sum_{j\in \dirpath_{\downstream}(i)} W_{i,j} g_j \label{eq:gss_upstream} \\
	h_i &= Q_i^T x_i +\sum_{j \in  \dirpath_{\upstream}(i)} R_{i,j} h_j  \label{eq:gss_downstream}  \\
	b_i &= D_i x_i + \sum_{j\in \dirpath_{\downstream}(i)} U_{i,j} g_j + \sum_{j \in \dirpath_{\upstream}(i) } P_{i,j} h_j \label{eq:gss_out}
\end{align}
\end{subequations}
for $i \in \nodeset_{\graph}$ are consistent with \eqref{eq:block_matrix_multiplication}.
\end{definition}

\begin{figure}
    \centering
    \begin{subfigure}[b]{0.23\textwidth}
        \includegraphics[width =\textwidth]{\MyPath/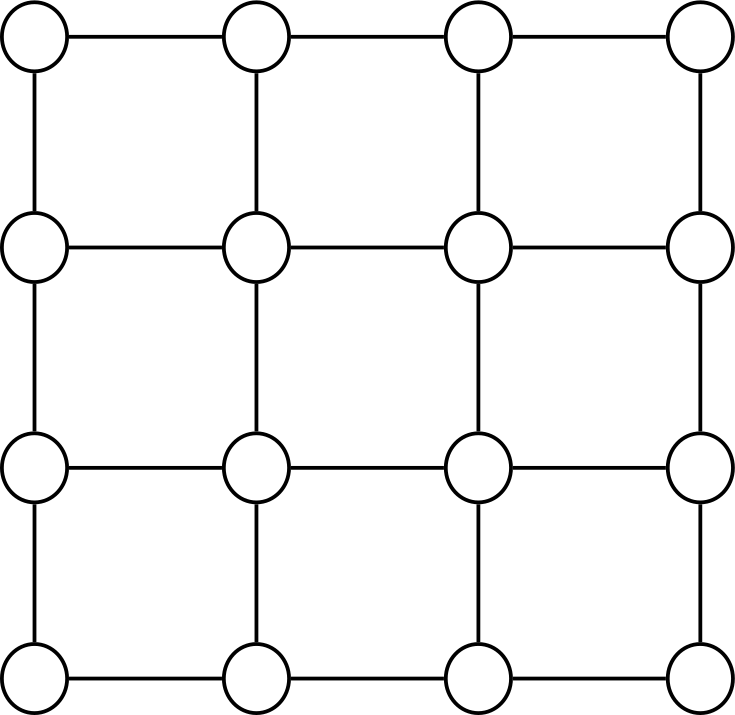} 
        \caption{} \label{subfig:2dmesh}
    \end{subfigure}
    ~ 
    \begin{subfigure}[b]{0.23\textwidth}
        \includegraphics[width = \textwidth]{\MyPath/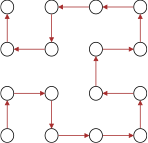}
        \caption{} \label{subfig:2dmesh_ham_path}
    \end{subfigure}
    ~
    \begin{subfigure}[b]{0.23\textwidth}
        \includegraphics[width = \textwidth]{\MyPath/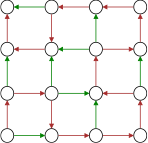}
        \caption{}
        \label{subfig:2dmesh_induced_order}
    \end{subfigure}
    ~ 
    \begin{subfigure}[b]{0.23\textwidth}
        \includegraphics[width = \textwidth]{\MyPath/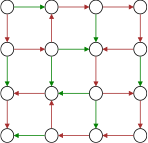}
        \caption{}
        \label{subfig:2dmesh_reverse_order}
    \end{subfigure}
    \caption{The $4\times 4$ 2D mesh graph $\graph$, shown in (\subref{subfig:2dmesh}), admits many different Hamiltonian paths. One is the Hilbert space-filling curve, shown in (\subref{subfig:2dmesh_ham_path}). This imbues every edge in graph $\graph$ with a natural orientation, where the orientation of the edge $(i,j)$ is chosen such that $i \rightsquigarrow j$, shown in (\subref{subfig:2dmesh_induced_order}). This represents the ``direction of information flow'' in computing the recurrence \eqref{eq:gss_downstream}. The corresponding flow for \eqref{eq:gss_upstream} is obtained by reversing the direction of all of the edges, shown in (\subref{subfig:2dmesh_reverse_order}). \ethancheck}\label{fig:hamiltonian_path_2dmesh}
\end{figure}

Definition~\ref{def:GSS} applies to a very general class of graph partitioned matrices. Any graph which contains a Hamiltonian path permits a $\graph$-semiseparable ($\graph$-SS, for short) representation. Figure~\ref{fig:hamiltonian_path_2dmesh} illustrates  how one can construct a causal and anti-causal flow on a two-dimensional mesh graph with the help of a Hamiltonian path. \ethancheck

The state-space equations \eqref{eq:gss} for $\graph$-SS representations are entirely explicit in that the $g$'s and $h$'s (and thus the product $Ax = b$) can be computed in sequence by summing matrix-vector products. If all the transition operators $V_i$, $W_{ij}$, $Q_i$, $R_{i,j}$, $D_i$, $U_{i,j}$, and $P_{i,j}$ have no more than $r$ columns, the recurrences \eqref{eq:gss} can be computed in $\mathcal{O}(r(m+n)(r+\max_i N_i))$ time, where $m$ is the number of edges in $\graph$. If the degree of $\mathbb{G}$ is bounded indepently of $n$ for a family of graphs and $n = \mathcal{O}(N/r)$ since $N_i = \mathcal{O}(r)$ for all $i \in \nodeset_\graph$, then $\graph$-SS representations can be multiplied in $\mathcal{O}(Nr)$ time.

Helpfully, the block entries of a $\graph$-SS matrix $A$ are given by an explicit expression. Define the ``state-transition" matrices by 
\begin{equation}   
\Phi_{\downstream}(k,\ell) = \begin{cases}  0  & k \succ  \ell \\
                   I  & k = \ell \\
                  \sum_{  s \in \dirpath_{\downstream}(k)  }  W_{k,s}  \Phi_{\downstream}(s,\ell)      &  k \prec \ell
                \end{cases}, \qquad  \Phi_{\upstream}(k,\ell) = \begin{cases}  \sum_{  s \in \dirpath_{\upstream}(k)   }  R_{k, s}  \Phi_\upstream (s,\ell)    & k \succ  \ell \\ 
                   I  & k = \ell \\
                    0  &  k \prec \ell
                \end{cases}   . \label{eq:transition}
\end{equation}
Then
\begin{equation}  
    [A]_{k\ell} = \begin{cases} \sum_{s \in \dirpath_{\upstream}(k)}  P_{k,s}  \Phi_{\upstream}(s,\ell) Q^T_\ell    & k \succ \ell, \\
    D_k & k = \ell,  \\
    \sum_{s \in \dirpath_{\downstream}(k)}  U_{k,s}  \Phi_{\downstream}(s,\ell) V^T_\ell     & k \prec \ell. 
    \end{cases}  \label{eq:GSSexplicit} \ethancheck
\end{equation}
A simple evaluation shows that when $\graph$ is the line graph, the above expression reduces to the SSS representation. For a general graph, the block entries of the matrix are sums with multiple terms, because there are generally multiple paths to reach one node from another. Expanding the transition matrices in terms of their basic components $W_{i,j}$ or $U_{i,j}$ can be a monumental task. \ethancheck 

The $\graph$-SS representations can be written compactly using matrix notation. Introduce matrix valued operators $ Z_\downstream [\cdot]$ and $Z_\upstream [\cdot]$ such that
\begin{displaymath}
\left[ Z_\downstream [\diagmat{W}] \right]_{i,j} =  \begin{cases}
W_{i,j} & \mbox{if }j \in  \dirpath_{\downstream}(i) \\
0 & \mbox{if }j\notin \dirpath_{\downstream}(i) \end{cases},\qquad \left[ Z_\upstream [\diagmat{R}] \right]_{i,j} =  \begin{cases}
R_{i,j} & \mbox{if }j \in  \dirpath_{\upstream}(i) \\
0 & \mbox{if }j\notin \dirpath_{\upstream}(i) \end{cases}. \ethancheck 
\end{displaymath}
One can then show that if one stacks $(g_i)_{i\in\nodeset_\graph}$ and $(h_i)_{i\in\nodeset_\graph}$ into long vectors $g$ and $h$, one has
\begin{equation}\label{eq:EDV_sparse}
    \begin{bmatrix}
    I - Z_\upstream [\diagmat{R}] & 0 & -\diagmat{Q}^T \\
    0 & I - Z_\downstream[\diagmat{W}] & -\diagmat{V}^T \\
    Z_\upstream[\diagmat{P}] & Z_\downstream[\diagmat{U}] & \diagmat{D}
    \end{bmatrix}\begin{bmatrix}
    h \\ g \\ x
    \end{bmatrix} = \begin{bmatrix}
    0 \\ 0 \\ b
    \end{bmatrix}, \ethancheck
\end{equation}
which gives to rise to the expression
\begin{equation} A = \diagmat{D} +  Z_{\downstream}[\diagmat{U}] (I- Z_{\downstream}[\diagmat{W}])^{-1} \diagmat{V}^T +  Z_{\upstream}[\diagmat{P}] (I-  Z_{\upstream}[\diagmat{R}])^{-1} \diagmat{Q}^T.
\end{equation}
This should be compared with the \nousedv{classical }diagonal form of an SSS representation \usedv{\eqref{eq:diagrepSSS}}\nousedv{\cite[Eq.~(3)]{chandrasekaran2018fast}}.
We are interested in the following general questions about these $\graph$-SS representations.

\begin{enumerate}[label=($\mathbb{G}$\arabic*)]
    \item When does a graph partitioned matrix have an efficient $\graph$-SS representation and how does this relate to GIRS? \label{item:GSS1}
    \item How can efficient $\graph$-SS representations be constructed? \label{item:GSS2}
    \item To what extent are the properties of SSS inherited by $\graph$-SS representations? \label{item:GSS3}
\end{enumerate}
%
%
In a partial answer to \ref{item:GSS1}, we show every graph-partitioned matrix possesses \textit{a} $\graph$-SS representation, but not necessarily one which is efficient. 

\begin{proposition}\label{prop:edv_universal}
Every graph partitioned matrix $(A,\graph)$ with $\graph$ admitting a Hamiltonian path $\dirpath$ can be described by a $\graph$-SS representation.\ethancheck 
\end{proposition}

\begin{proof}
    Let $P$ be the permutation which puts $PAP^T$ consistent with the total order induced by $\dirpath$. We may then construct a SSS representation for $PAP^T$, which subsequently can be convert in $\graph$-SS representation by setting matrix entries associated with the induced edges equal to zero.  Thus, in this way, existence of the SSS representation guarantees existence of the $\graph$-SS representation.\ethancheck 
\end{proof}

The result above shows that $\graph$-SS representations are universal. However, while Proposition \ref{prop:edv_universal} guarantees the existence of \textit{a} $\graph$-SS representation, it is highly unlikely that this representation furnished by the construction in theorem will be optimal or even efficient in any sense. After all, the matrix entries associated with the induced edges are not utilized at all.

One reason we are optimistic about the potential of $\graph$-SS in representing GIRS matrices is that a converse result holds: every matrix possessing an efficient $\graph$-SS representation is GIRS with a small constant. Refer to the state dimensions of $g_i$ and $h_i$ as $r^{\rm g}_i$ and $r^{\rm h}_i$ for $i \in \nodeset_\graph$. \usedv{Now notice that every $\graph$-SS representation can be converted into an equivalent DV representation. This is an immediate consequence of merging the two states:
 \begin{subequations}
 \begin{align}
	\begin{bmatrix} g_i \\  h_i \end{bmatrix} & = \begin{bmatrix} V_i^T \\ Q^T_i \end{bmatrix} x_i + \sum_{j\in \dirpath_{\downstream}(i)}  \begin{bmatrix} W_{i,j}  & 0 \\ 0 & 0
	\end{bmatrix} \begin{bmatrix} g_j \\  h_j \end{bmatrix}  + \sum_{j\in \dirpath_{\upstream}(i)}  \begin{bmatrix} 0  & 0 \\ 0 & R_{i,j}
	\end{bmatrix} \begin{bmatrix} g_j \\  h_j \end{bmatrix}   \\
	b_i &= D_i x_i + \sum_{j\in \dirpath_{\downstream}(i)}  \begin{bmatrix} U_{i,j}  & 0 \
	\end{bmatrix} \begin{bmatrix} g_j \\  h_j \end{bmatrix}  + \sum_{j\in \dirpath_{\upstream}(i)}  \begin{bmatrix}  0 & P_{i,j}
	\end{bmatrix} \begin{bmatrix} g_j \\  h_j \end{bmatrix}.  \ethancheck .
\end{align} \label{eq:GSSISDV}
\end{subequations}
Thus, from Proposition~\ref{prop:idv_girs}, we can immediate conclude the following.}\nousedv{Then we have the following.}

\begin{proposition}\label{prop:GSS_is_GIRS}
Let $(A,\graph)$ be a graph-partitioned matrix possessing a $\graph$-SS representation and define the maximal state dimension $r = \max_{i \in \nodeset(\graph)} \{r_i^g,r_i^h\}$. Then $A$ is GIRS-$4r$.
\end{proposition}

\nousedv{In an extended preprint version of this article, we also investigate a different class of representations for GIRS matrices which we call Dewilde-van der Veen (DV) Representations \cite[Sec.~2.3]{chandrasekaran2019graph}. These representations have two interesting properties: they can be very easily constructed for sparse matrices and from a DV representation of $A$ one can efficiently compute a DV representation for $A^{-1}$ of the same size. The downside of these representations is they are implicit, requiring the solution of a sparse linear system of equations to compute matrix-vector products.

The $\graph$-SS representations introduced in this section are a subclass of these more general DV representations, but we are unaware of a general construction method for them (even for sparse matrices). Moreover, while one can easily construct a \textit{DV representation} of the inverse $A^{-1}$ of a matrix $A$ represented as a $\graph$-SS representation, we are unaware of a general algorithm for constructing a \textit{$\graph$-SS representation} for the inverse $A^{-1}$. In fact, we construct in the next section an example where the minimal $\graph$-SS representation of $A^{-1}$ is a constant factor larger size than the representation for $A$.}

\usedv{How do $\graph$-SS representations compare to DV representations? Under the present state of our understanding, they suffer from a couple of important deficits. Firstly, we are unaware of any general algorithm for constructing practically useful $\graph$-SS representations beyond the SSS case, even for sparse matrices. Further, given a $\graph$-SS representation of $A$, we know no algorithm for computing a $\graph$-SS algorithm for $A^{-1}$. In fact, we construct in the next section an example where the minimal $\graph$-SS representation of $A^{-1}$ is a constant factor larger size than the representation for $A$. All of these deficits are not true for DV representations. However, $\graph$-SS representations possess a highly useful property that DV representations do not: a linear-time multiplication algorithm. This makes searching for $\graph$-SS representations an appealing prospect, as they could prove quite useful if constructable.}

\nousedv{Still, we}\usedv{We} take the connection between GIRS matrices and $\graph$-SS representations in Proposition~\ref{prop:GSS_is_GIRS}, the existence of implicit $\graph$-SS-like representations for sparse matrices and their inverses, and the forthcoming results on $\graph$-SS representations for the special case when $\graph$ a cycle graph to be promising signs that $\graph$-SS may be effective algebraic representations for general GIRS matrices.

\section{A case study: the cycle semi-seperable representation} \label{sec:css}

As a beachhead to tackling more complicated graphs, we consider the problem of constructing a $\graph$-SS representation for a graph-partitioned matrix $(A,\mathbb{G})$ with $\mathbb{G}$ being the cycle graph consisting of $n$ nodes $\nodeset_\graph = \{1,2,3,\ldots,n\}$. Our aim in doing so is not to propose these cycle semi-separable (CSS) representations as an improvement over SSS for practical applications, but rather to investigate the questions \ref{item:GSS1}, \ref{item:GSS2}, and \ref{item:GSS3} in the simplest example (other than the line graph). 

\begin{figure}
    \centering
    \begin{subfigure}[b]{0.3\textwidth}
        \includegraphics[width =\textwidth]{\MyPath/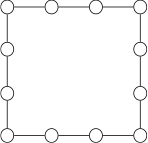} 
        \caption{} \label{subfig:circular}
    \end{subfigure}
    ~ 
    \begin{subfigure}[b]{0.3\textwidth}
        \includegraphics[width = \textwidth]{\MyPath/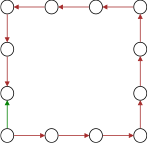}
        \caption{} \label{subfig:circular_upstream}
    \end{subfigure}
    ~
    \begin{subfigure}[b]{0.3\textwidth}
        \includegraphics[width = \textwidth]{\MyPath/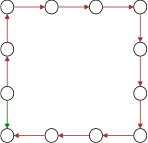}
        \caption{}
        \label{subfig:circular_downstream}
    \end{subfigure}
    \caption{The $10$ node cycle graph $\graph$, shown in (\subref{subfig:circular}). Choosing the Hamiltonian path $1\to2\to\cdots \to9\to 10$ leads to the induced causal/upstream (\subref{subfig:circular_upstream}) and anticausal/downstream (\subref{subfig:circular_downstream}) digraphs. \ethancheck}\label{fig:circularSSS}
\end{figure}

Taking our Hamiltonian path to be $1 \to 2 \to \cdots \to n$, the explicit flow is illustrated in Figure~\ref{fig:circularSSS}.  According to \eqref{eq:gss}, the state-space equations reduce to
\begin{equation*} 
g_{k}  =  V^T_{k} x_{k} + W_{k} g_{k+1} \ethancheck
\end{equation*}
for $k=2,3,\ldots,n-1$  with $g_n = V^T_n x_n$ in the case of the downstream flow \eqref{eq:gss_upstream}, and
\begin{equation*}
h_k = Q^T_k x_k + R_{k} h_{k-1} \ethancheck
\end{equation*}
for $k=n-1,n-2,\ldots,2$ with $h_1 = Q^T_1 x_1$ in the case of the upstream flow \eqref{eq:gss_downstream}. The output equation \eqref{eq:gss_out} is then
\begin{equation*}
b_k  = \begin{cases}  D_1 x_1 + U_1 g_2 + U_0 g_n &  k= 1 \\
                     D_n x_n + P_{n} h_{n-1} + P_{0} h_1 & k= n     \\
                     D_k x_k + U_{k} g_{k+1} + P_{k}h_{k-1} & \mbox{otherwise}.\\
\end{cases}\ethancheck
\end{equation*}
The resulting matrix $A\in \mathbb{R}^{N\times N}$, which is of dimension $N=N_1+\ldots + N_n$, has a relatively simple structure where the $(k,\ell)$-th block entry of $A$ has the expression
\begin{equation}
    [A]_{k\ell} = \begin{cases} 
    D_k & k=\ell  \\
    U_{k} W_{k-1} W_{k-2} \cdots W_{\ell+1} V^T_{\ell}   & k< \ell,\quad  \ell \ne n  \\
    U_{n} W_{n-1} W_{n-2} \cdots W_{2} V^T_{1} + U_0 V_n^T   & k=1,\quad \ell=n  \\
     P_{k} R_{k+1} R_{k+2} \cdots R_{\ell-1} Q^T_{\ell}   & k> \ell,\quad  k \ne n  \\
     P_{1} R_{2} R_{3} \cdots R_{n-1} Q^T_{n} + P_0 Q_1^T   & k=n,\quad  \ell = 1 
    \end{cases}. \ethancheck \label{eq:CSS} 
\end{equation}
The sizes of the generators of the representation are given by $U_i\in \mathbb{R}^{N_i \times \rg_{i} }$, $W_i\in \mathbb{R}^{ \rg_{i-1} \times \rg_{i} }$, $V_i\in \mathbb{R}^{N_i \times \rg_{i-1} }$, $D_i\in \mathbb{R}^{N_i\times N_i}$, $P_i\in \mathbb{R}^{N_i \times r^{h}_{i-1} }$, $R_i\in \mathbb{R}^{\rh_{i}\times \rh_{i-1}}$,  and $Q_i\in \mathbb{R}^{N_i \times \rh_{i}}$ where $\rg_i$ and $\rh_i$ are the sizes of $g_i$ and $h_i$ respectively. As boundary cases, we have $U_0 \in\mathbb{R}^{N_1\times \rg_{n-1}}$ and $P_0 \in\mathbb{R}^{N_n   \times \rh_{1}}$. In the case of $n=4$, \eqref{eq:CSS} reduces to
\begin{equation}
A = \begin{bmatrix}
{D_1  }              &          U_1V_2^T             &    U_1W_2V_3^T  &    U_1W_2W_3V_4^T + U_0 V_4^T     \\
P_{2} Q_1^T               &  {D_2}          &  U_2 V_3^T & U_2W_3V_4^T  \\
P_{3} R_{2} Q^T_1          &  P_{3} Q^T_2         &  {D_3} & U_3 V_4^T      \\
P_{4} R_{3} R_{2} Q^T_1  + P_0 Q_1^T    &  P_{4} R_{3} Q^T_2  & P_{4} Q^T_3   & {D_4} 
\end{bmatrix}. \ethancheck 
\label{eq:CSSSn4}
\end{equation}
The form of $A$ is identical to an SSS representation, as shown in \nousedv{\eqref{eq:sss}}\usedv{\eqref{eq:SSSSn4}}, except for the additional terms $P_0Q_1^T$ and $U_0V^T_n$ in the bottom left and top right corners of the matrix. 
\subsection{The CSS construction algorithm}\label{sec:circular_sss_construction}
CSS matrices are universal as any matrix $A$ can be put into this representation. We may simply construct a SSS representation of $A$ and then set $U_0$ and $P_0$ to zero matrices of appropriate sizes. More generally, we may construct a CSS representation as follows. Define:
\begin{equation}
    A(X,Y) = \begin{bmatrix}
    A_{11} & A_{12} & A_{13} & \cdots & Y  \\
    A_{21} & A_{22} & A_{23} & \cdots & A_{2n}  \\
    A_{31} & A_{32} & A_{33} & \cdots & A_{3n}   \\
    \vdots & \vdots & \vdots & \ddots & \vdots  \\
    X & A_{n2} & A_{n3} & \cdots & A_{nn}
\end{bmatrix},\qquad     \CB(E_1,E_2) = \begin{bmatrix}
    0 & 0 & \cdots & 0 & E_1 \\
    0 & 0 & \cdots & 0 & 0 \\
    \vdots & \vdots & \ddots & \vdots & \vdots \\
    0 & 0 & \cdots & 0 & 0 \\
    E_2 & 0 & \cdots & 0 & 0
    \end{bmatrix} \label{eq:corner_block}
\end{equation}
and express
\begin{equation} A  =    A(X,Y)  + \CB(A_{1n} - Y,    A_{n1} - X ). \ethancheck \label{eq:splitting}
\end{equation}
The block entries $A_{n1}\in\mathbb{R}^{N_n \times N_1}$ and $A_{1n}\in\mathbb{R}^{N_1 \times N_n} $ are replaced by arbitrary placeholders $X\in\mathbb{R}^{N_n \times N_1} $ and $Y\in\mathbb{R}^{N_1 \times N_n}$, respectively. Observe that the off-diagonal (Hankel) blocks are now dependent on the choice of $X$ and $Y$ \nithincheck Specifically, the $k$th upper and lower Hankel blocks are given by
\begin{equation}
    \mathcal{H}_{k}(X) := \begin{bmatrix}
    A_{(k+1)1} & \cdots & A_{(k+1)k} \\
    \vdots &     & \vdots \\
    X & \cdots & A_{nk} 
    \end{bmatrix},  \quad \mathcal{G}_{k}(Y) := \begin{bmatrix}
    A_{1(k+1)} & \cdots & Y \\
    \vdots &     & \vdots \\
    A_{k(k+1)} & \cdots & A_{n(k+1)} 
    \end{bmatrix}  \ethancheck    \label{eq:hankelblock}
\end{equation}
for $k = 1,2,\ldots,n-1$. To compute the CSS representation, we proceed by first constructing a SSS representation for $A(X,Y)$. Recognizing that we would like to set $P_0 Q^T_1 = A_{n1}- X$ and $U_0 V^T_n = A_{n1}- Y$,  the SSS representation can be morphed into a CSS representation to include the "perturbations" described by the second term in \eqref{eq:splitting}. We have the following algorithm.\ethancheck 
\begin{alg}[CSS construction]\label{alg:CSSconstruct} 
Let $(A,\mathbb{G})$ be a graph-partitioned matrix with $\mathbb{G}$ as the cycle graph consisting of $n$ nodes.\ethancheck 
\begin{enumerate}
\item  Select a $X\in\mathbb{R}^{N_n \times N_1} $ and $Y\in\mathbb{R}^{N_1 \times N_n}$ and express $A \in\mathbb{R}^{N \times N}$ as per \eqref{eq:splitting}.\ethancheck 
\item  Construct a SSS representation for $A(X,Y)$, i.e.
$$
A(X,Y) = \resizebox{0.82\textwidth}{!}{$\begin{bmatrix}
D_1 &   U_1 V_2^T             &    U_1 W_2 V_3^T  &    U_1 W_2 W_3 V_3^T & \cdots &  U_1 W_2 W_3\cdots \tilde{W}_{n-1}\tilde{V}_n^T \\
\tilde{P}_{2} \tilde{Q}_1^T   & D_2  & U_2 V_3^T & U_2 W_3 V_4^T  & \cdots &          U_2 W_3 \cdots \tilde{W}_{n-1} \tilde{V}_n^T   \\
P_{3} \tilde{R}_{2} \tilde{Q}_1^T          &  P_{3} Q^T_2 & D_3 & U_3 V^T_4 & \cdots &   U_3W_4 \cdots \tilde{W}_{n-1}\tilde{V}_n^T   \\
P_{4} R_{3} \tilde{R}_{2} \tilde{Q}_1^T     &  P_{4} R_{3} Q^T_2  & P_{4} Q^T_3  & D_4  & \cdots &  U_4 W_5 \cdots \tilde{W}_{n-1}\tilde{V}_n^T  \\
\vdots   & \vdots & \vdots & \vdots & \ddots & \vdots  \\
P_{n} R_{n-1} \cdots \tilde{R}_{2} \tilde{Q}_1^T     &  P_{n} R_{n-1} \cdots R_{3}  Q^T_2         &  P_{n} R_{n-1} \cdots R_{4}  Q^T_3 & P_{n} R_{n-1} \cdots R_{5}  Q^T_4 & \cdots & D_n
\end{bmatrix}$.}\ethancheck 
$$
\item
Replace the terms denoted with a tilde, i.e. $\tilde{Q}_1^T$, $\tilde{P}_2$, $\tilde{R}_2$, $\tilde{U}_{n-1}$, $\tilde{V}^T_n$ and $\tilde{W}_{n-1}$, to reflect perturbations caused by $A_{n1}-X$ and $A_{1n}-Y$. Since the rows of matrix $\tilde{Q}_1^T$ are only a row basis for the first Hankel $\mathcal{H}_1(X)$ but do not necessarily span the rows of $A_{n1}-X$,  compute a low rank factorization 
$$  \begin{bmatrix} \tilde{Q}_1^T \\ A_{n1}-X \end{bmatrix}  = Z Q_1^T$$

to obtain $Q_1$ whose columns are a row basis for the space spanned by the rows of $\mathcal{H}_1(X)$ and $A_{n1}-X$. Since $Q_1^T$ has full row rank, its pseudo-inverse is a right inverse $(Q_1^T)^{\dagger}$. After replacing with $\tilde{Q}_1$ with $Q_1$, we replace $\tilde{P}_2 $  with $P_2 = \tilde{P}_2 \tilde{Q}^T_1 (Q_1^T)^{\dagger} $  and $\hat{R}_2$ with ${R}_2 = \hat{R}_2 \tilde{Q}^T_1 (Q_1^T)^{\dagger}$. Finally, set
$$ P_0 =  (A_{n1}-X)(Q_1^T)^{\dagger} $$
Similar formulas apply also for $U_{n-1}$, ${V}^T_n$, ${W}_{n-1}$ and $U_0$. \ethancheck \nithincheck
\end{enumerate}
\end{alg}

\subsection{Finding the minimal CSS representation}

Similar to SSS, for a given matrix $A$, there will be many CSS representations. Algorithm~\ref{alg:CSSconstruct} provides a means for computing many CSS representations with different choices of $X$ and $Y$. It is of interest to construct the representation of smallest dimensions, i.e.\ requiring the minimum amount numbers to be stored. \ethancheck 

The CSS representation involves only two terms more than the SSS representation, namely $P_0 \in\mathbb{R}^{N_n\times \rh_1}$ and $U_0 \in\mathbb{R}^{N_1\times \rg_{n-1}}$. One can convince oneself that regardless of the choice of $X$ and $Y$ in Algorithm~\ref{alg:CSSconstruct}, the state dimensions $\rh_1$ and $\rg_{n-1}$ can be made no less than the ranks of their corresponding (unpeturbed) Hankel blocks $\rh_1 \ge \rank \mathcal{H}_1(0)$ and $\rg_{n-1} \ge \rank \mathcal{G}_{n-1}(0)$, with equality if, and only if, $X$ and $Y$ lie in the row spaces of their corresponding Hankel blocks $\mathcal{H}_1(0)$ and $\mathcal{G}_{n-1}(0)$. Since the state dimensions of these terms cannot be reduced, the total size of the CSS representation will be minimized if $X$ and $Y$ are chosen so that the other state dimensions $\rh_i$ and $\rg_i$ are made as small as possible. \textit{A priori}, it may be possible that there exists no CSS representation which is \textit{uniformly minimal} in the sense that for any CSS representation of $A$ with ranks $\rghat_i, \rhhat_i$,
\begin{equation*}
    \rg_i \le \rghat_i,\quad \rh_i \le \rhhat_i, \quad \textnormal{for }i = 1,\ldots,n-1.\ethancheck
\end{equation*}
If no such uniformly minimal representation existed we would need to settle for minimizing some measure of the total size of the representation, such as the sum of state dimensions $\sum^{n-1}_{i=1} \left({\rg_{i}} + {\rh_{i}}\right)  $. Fortunately, we shall show that, in fact, it is possible to produce a uniformly minimal CSS representation using the construction Algorithm~\ref{alg:CSSconstruct}. \ethancheck
\begin{theorem}\label{thm:CSS_construction}
Let matrices $\hat{X}\in \mathbb{R}^{N_n \times N_1}$ and  $\hat{Y}\in \mathbb{R}^{N_1 \times N_n}$ satisfy
\begin{equation}\label{eq:Xhat_Yhat}
    \rank  \mathcal{H}_i(\hat{X}) = \min_{ X \in \mathbb{R}^{N_n \times N_1} } \rank \mathcal{H}_i(X) \quad \mbox{and} \quad  \rank  \mathcal{G}_i(\hat{Y}) = \min_{ Y \in \mathbb{R}^{N_1 \times N_n} } \rank \mathcal{G}_i(Y)
\end{equation}
for  $i=1,2,\ldots, n-1$. For these matrices, the CSS construction process, as described in Algorithm~\ref{alg:CSSconstruct}, will generate a CSS representation which is uniformly minimal: given any other representation with state dimensions $\rghat_i$ and $\rhhat_i$ for $1\le i \le n-1$, we must have
$$ \rh_i \leq \hat{r}^h_i \quad \mbox{and} \quad  \rg_i \leq \hat{r}^g_i.   $$
for $i=1,2,\ldots, n-1$. \ethancheck
\end{theorem}
\begin{proof}
We shall only prove $\rh_i \leq \hat{r}^h_i$ for $i=1,2,\ldots, n-1$ since the proof for $\rg_i \leq \hat{r}^g_i$ follows an analogous path.  The values for $\rh_i$ are equal to the number of rows in $Q_i^T$. For $i=2,\ldots,n-1$, the minimal number of columns in $Q_i$ is given by
$$  \rh_i= \rank \mathcal{H}_i(X).\ethancheck $$
This is a consequence of the properties of the SSS representation (Theorem~\ref{thm:sss}). Henceforth, the statement $\rh_i \leq \rhhat_i$ holds true by definition of $\hat{X}$ for $i=2,\ldots,n-1$. On the other hand, we know that for $i=1$, we have
\begin{eqnarray*}
\rh_1   & = &  \rank Q_1 \\
        & = & \rank \begin{bmatrix} \tilde{Q}_1 & (A_{n1}-X)^T \end{bmatrix} \\
    & = & \rank \begin{bmatrix} \mathcal{H}^T_1(X) & (A_{n1}-X)^T \end{bmatrix} \\
    & = &  \rank \begin{bmatrix} A^T_{11} & A^T_{21} & \cdots & A^T_{(n-1)1} & X^T &  A^T_{n1}-X^T  \end{bmatrix} . \ethancheck
\end{eqnarray*}
Since $X$ minimizes the rank of the first Hankel block, we have
$$ R(X^T) \subseteq R\left(\begin{bmatrix} A^T_{11} & A^T_{21} & \cdots &  A^T_{(n-1)1} \end{bmatrix}\right). \ethancheck  $$
Thus
\begin{equation*}
    R\left(\begin{bmatrix} A^T_{11} & A^T_{21} & \cdots & A^T_{(n-1)1} & X^T &  A^T_{n1}-X^T ) \end{bmatrix} \right) = R(\mathcal{H}^T_1(A_{n1})).
\end{equation*}
Note that for any CSS representation of $A$, $Q_1^T$ has to be a row basis for $\mathcal{H}_1(A_{n1})$, so $\rhhat_i \ge \rank \mathcal{H}_1(A_{n1}) = \rh_i$. This completes the proof of uniform minimality.\ethancheck \nithincheck
%
\end{proof}
\begin{remark}
The existence of $\hat{X}$ and $\hat{Y}$ satisfying \eqref{eq:Xhat_Yhat} for  $i=1,2,\ldots, n-1$ is a fairly delicate question and is studied in Section~\ref{sec:hankelminimization}. We show in Theorem~\ref{thm:rankcompletion} that the existence of these matrices are always guaranteed. 
\end{remark}
\begin{corollary}
For a uniformly minimal CSS representation, the state dimensions are given by
\begin{align*}
    \rh_1 &=  \rank \mathcal{H}_1(A_{n1}),&\quad \rh_i &= \min_{ X \in \mathbb{R}^{N_n \times N_1} } \rank \mathcal{H}_i(X) \mbox{ for  $i=2,3,\ldots, n-1$},\\
    \rg_{n-1} &=  \rank \mathcal{G}_{n-1}(A_{1n}),&\quad \rg_i &= \min_{ Y \in \mathbb{R}^{N_1 \times N_n} } \rank \mathcal{G}_i(Y) \mbox{ for  $i=1,2,\ldots, n-2$}.\ethancheck
\end{align*}
\end{corollary}
Thus, the essential difficulty of the CSS construction boils down to a specific matrix completion problem: given a block triangular array, how may the bottom left corner block be chosen so that the ranks of all rectangular subblocks containing the corner block are minimized simultaneously? As the following example shows, for some basic matrices, this problem can be solved from a quick direct inspection.

\begin{example} \label{ex:CSSconstruction}
Consider the matrix
\begin{equation}\label{eq:circular_tridiagonal}
A = \left[
\begin{array}{c|c|c }
\begin{array}{c c}  b & a \\ a & b \end{array} & \begin{array}{c c} 0 & 0 \\ a & 0  \end{array} & \begin{array}{c c} 0 & a \\ 0 & 0 \end{array}  \\
\hline
\begin{array}{c c} 0 & a \\ 0 & 0 \end{array} & \begin{array}{c c}  b & a \\ a & b \end{array} & \begin{array}{c c} 0 & 0 \\ a & 0 \end{array}  \\
\hline
\begin{array}{c c} 0 & 0 \\ a & 0 \end{array} & \begin{array}{c c} 0 & a \\ 0 & 0 \end{array}  & \begin{array}{c c}  b & a \\ a & b \end{array} \end{array} \right],
\end{equation}
where $a$ and $b$ are scalars. We must solve the overlapping Hankel block minimization problem for
\begin{equation*}
A(X,Y) = \left[
\begin{array}{c|c|c }
\begin{array}{c c}  b & a \\ a & b \end{array} & \begin{array}{c c} 0 & 0 \\ a & 0  \end{array} & Y  \\
\hline
\begin{array}{c c} 0 & a \\ 0 & 0 \end{array} & \begin{array}{c c}  b & a \\ a & b \end{array} & \begin{array}{c c} 0 & 0 \\ a & 0 \end{array}  \\
\hline
X & \begin{array}{c c} 0 & a \\ 0 & 0 \end{array}  & \begin{array}{c c}  b & a \\ a & b \end{array} \end{array} \right].
\end{equation*}
The solution, for this particular example is simple. We may set
\begin{displaymath}
\hat{X} =  \begin{bmatrix} 0 & \alpha_1 \\ 0 & 0 \end{bmatrix}, \quad \hat{Y} = \begin{bmatrix} 0 & 0 \\ \alpha_2 & 0 \end{bmatrix}.
\end{displaymath}
where $\alpha_1,\alpha_2\in\mathbb{R}$ can be chosen freely.
To keep things simple, we can pick $\alpha_1 = \alpha_2 = 0$ for step 1 of Algorithm~\ref{alg:CSSconstruct}. Proceeding with step 2, the SSS representation for $A(\hat{X},\hat{Y})$ can be written as
\begin{equation*}
A(\hat{X}, \hat{Y}) = \left[
\begin{array}{ccc }
D_1 & {U}_1 {V}^T_2  & {U}_1 \tilde{W}_2 \tilde{V}^T_3   \\
\tilde{P}_2 \tilde{Q}^T_1 & D_2 &  \tilde{U}_2 \tilde{V}^T_3    \\
P_3 \tilde{R}_2 \tilde{Q}^T_1 & P_3 Q^T_2 & D_3 \end{array} \right]
\end{equation*}
with
\begin{displaymath}
D_1 =D_2 = D_3 = \begin{bmatrix}  b & a \\ a & b \end{bmatrix},\quad \tilde{Q}_1 ={Q}_2 = \begin{bmatrix} 
0 \\ 1
\end{bmatrix}, \quad \tilde{P}_2 = {P}_3 = \begin{bmatrix} 
a \\ 0
\end{bmatrix}, \quad \tilde{R}_2 = \begin{bmatrix} 0 & 0 \end{bmatrix},
\end{displaymath} 
\begin{displaymath}
\tilde{V}_2 = {V}_3 = \begin{bmatrix} 
1 \\ 0
\end{bmatrix}, \quad \tilde{U}_2 = {U}_3 = \begin{bmatrix} 
0 \\ a
\end{bmatrix}, \quad \tilde{W}_2 = \begin{bmatrix} 0 & 0 \end{bmatrix}.
\end{displaymath} 
Finally step 3 may lead to
\begin{equation*}
A = \left[
\begin{array}{ccc }
D_1 & {U}_1 {V}^T_2  & {U}_1 {W}_2 {V}^T_3 + U_0 V^T_3   \\
{P}_2 {Q}^T_1 & D_2 &  {U}_2 {V}^T_3    \\
P_3 {R}_2 {Q}^T_1 + P_0 {Q}^T_1 & P_3 Q^T_2 & D_3 \end{array} \right]
\end{equation*}
where
\begin{displaymath}
{Q}_1 = \begin{bmatrix} 
0  & 1 \\ 1 & 0
\end{bmatrix}, \quad {P}_0 =  \begin{bmatrix} 
0 & 0 \\ 0 & a
\end{bmatrix}, \quad P_2 =  \begin{bmatrix} 
0 & a \\ 0 & 0
\end{bmatrix}, \quad  R_2 = \begin{bmatrix} 
0 & 0
\end{bmatrix},
\end{displaymath} 
\begin{displaymath}
{V}_3 = \begin{bmatrix} 
1  & 0 \\ 0 & 1
\end{bmatrix}, \quad {U}_0 =  \begin{bmatrix} 
0 & a \\ 0 & 0
\end{bmatrix}, \quad U_2 =  \begin{bmatrix} 
0 & 0 \\ a & 0
\end{bmatrix}, \quad  W_2 = \begin{bmatrix} 
0 & 0
\end{bmatrix}.
\end{displaymath} 
As can be seen from the steps of algorithm~\ref{alg:CSSconstruct} and the rank completion problem, the solution for uniformly minimal CSS representation can be highly non-unique. \nithincheck
\end{example}
The general problem is highly nontrivial. In Section~\ref{sec:hankelminimization}, we address this problem fully by proving existence of such a matrix through the formulation of a construction algorithm for finding it. In Appendix~\ref{app:css_properties}, we discuss some properties possessed by CSS representations, in particular their stability properties under inversion (see Theorem~\ref{thm:CSSInversion}).

\begin{remark}\label{rmk:CSS_vs_SSS}
We caution the reader that CSS is not to be taken as a better representation than SSS in practice. As shown in Appendix~\ref{app:css_properties}, CSS can be significantly better than SSS for examples of the form $A + \CB(E_1,E_2)$ where $A$ is SSS-$r$ and $\rank E_1 = \rank E_2 = R \gg r$. In this case, the SSS representation has size proportional to $nR$ whereas the CSS representation has size proportional to $nr$. However, if we perform a permutation and repartition, we can write
\begin{equation*}
    P(A+\CB(E_1,E_2))P^T = \begin{bmatrix}
    \twobytwo{A_{11}}{A_{1n}+E_1}{A_{n1}+E_2}{A_{nn}} & \twobyone{A_{12}}{A_{n2}} & \twobyone{A_{13}}{A_{n3}} & \cdots & \twobyone{A_{1(n-1)}}{A_{n(n-1)}} \\
    \onebytwo{A_{21}}{A_{2n}} & A_{22} & A_{23} & \cdots & A_{2(n-1)} \\
    \onebytwo{A_{31}}{A_{3n}} & A_{32} & A_{33} & \cdots & A_{3(n-1)} \\
    \vdots & \vdots & \vdots & \ddots & \vdots \\
    \onebytwo{A_{(n-1)1}}{A_{(n-1)n}} & A_{(n-1)2} & A_{(n-1)3} & \cdots & A_{(n-1)(n-1)}
    \end{bmatrix},
\end{equation*}
which is an SSS-$r$ matrix and can thus can be stored with size proportional $Nr$. This, in asymptotic terms, CSS and permuted SSS have the same storage complexity for representing such matrices and without the need for the solution of the low rank completion problem.

We reiterate that our goal is this paper is not to propose CSS as a practical alternative to SSS for most problems. However, it is our conjecture that the techniques used to construct the CSS representation will generalize to allow us to construct $\graph$-SS representations for more complicated graphs $\graph$, which will ultimately be more efficient than SSS.
\end{remark}

\section{The overlapping Hankel low-rank completion problem} \label{sec:hankelminimization}

This section addresses the low rank completion problem of finding $\hat{X}$ which simultaneously minimizes the ranks of all Hankel blocks \eqref{eq:hankelblock}. Since the lower and upper triangular parts of the matrix are equivalent, we shall focus on minimizing all of the lower Hankel blocks.
%
%
We will work towards the following general result:

\begin{theorem} \label{thm:rankcompletion}
As defined by \eqref{eq:hankelblock}, let $\mathcal{H}_k(X)$ denote the Hankel blocks corresponding to the lower triangular part of $A$. There exists a $\hat{X}\in \mathbb{R}^{N_n \times N_1}$ solving the overlapping Hankel block low-rank completion problem
\begin{equation}\label{eq:ohblrcp}
\rank  \mathcal{H}_k(\hat{X}) = \min_{ X \in \mathbb{R}^{N_n \times N_1} } \rank \mathcal{H}_k(X),\qquad k=1,2,\ldots, n-1. \tag{$\mathcal{H}$-LRCP}
\end{equation}
That is, the ranks of all Hankel blocks are simultaneously minimized. \ethancheck
\end{theorem}
In \cite{EGC21}, we present a proof of this result with a different construction than the one we originally discovered. We shall leave the present section in its current form to document our original approach to this problem.

Our proof to this theorem will be constructive and will generate a particular solution $\hat{X}$. 
The broad outline of our construction strategy is as follows. First, note that any two consecutive Hankel blocks $\mathcal{H}_k(X)$ and $\mathcal{H}_{k+1}(X)$ can be obtained from one another by first removing rows off of the top of $\mathcal{H}_k(X)$ and then adding columns to the right. Thus, if we construct the complete set of all $X$ minimizing the rank of the first Hankel block and are able to update the set when rows are removed or columns are added to the block, then we can iteratively sweep through the Hankel blocks in sequence, constructing a set of common solutions to the first $k$ Hankel blocks:
\begin{equation}
    \mathcal{S}_k := \left\{ \hat{X}\in\mathbb{R}^{m_2\times m_1}:\quad  \rank  \mathcal{H}_k(\hat{X}) = \min_{ X \in \mathbb{R}^{N_n \times N_1} } \rank \mathcal{H}_\ell(X),\qquad \ell=1,2\ldots,k  \right\} \label{eq:solutionsets}
\end{equation}
It will be shown that $\mathcal{S}_k$ remains nonempty after all $n-1$ Hankel blocks have been considered, and hence, we can simply select an arbitrary element of $\mathcal{S}_{n-1}$ as our candidate solution $\hat{X}$. As Example~\ref{ex:CSSconstruction} already highlighted, this set can contain more than element, i.e.\ the solution to the overlapping Hankel low-rank completion problem is non-unique.  \ethancheck \nithincheck

\subsection{Supporting lemmas} \label{sec:supportinglemmas}

Before we discuss the details of the main proof for Theorem~\ref{thm:rankcompletion} in Section~\ref{sec:proof_of_main_result}, we shall first derive several supporting lemmas. The proof of Theorem~\ref{thm:rankcompletion} is constructed based on the following intermediate results:
\begin{enumerate}[label={(\roman*)}]
    \item \emph{Lemma~\ref{thm:2x2_low_rank_completion}}. We provide an exact characterization for the full solution set of the block two-by-two low-rank completion problem
    \begin{equation}
    \min_{X\in \mathbb{R}^{m_2\times m_1}} \rank \begin{bmatrix}
        A & B \\
        X & C
        \end{bmatrix}, \tag{LRCP} \label{eq:lrcp}
    \end{equation}
     where $A \in \real^{m_1\times n_1}$, $B \in \real^{m_1\times n_2}$, $X \in \real^{m_2\times n_1}$, $C \in \real^{m_2\times n_2}$.
    This result is stated and derived in Section~\ref{sec:2x2_low_rank_completion}.
    
      \item \emph{Lemma~\ref{lem:addition_of_cols} and Lemma~\ref{lem:removal_of_rows}}. We develop a method for constructing a restricted \textit{nonempty} common solution set for \eqref{eq:lrcp} and the low-rank completion problem with additional columns
    \begin{equation}
    \min_{X\in \mathbb{R}^{m_2\times m_1}} \rank \begin{bmatrix}
        A & \BG \\
        X & \CH
        \end{bmatrix},\tag{LRCP$+$Cols} \label{eq:lrcpcols}
    \end{equation}
    where
    \begin{equation}\label{eq:BGCH}
        \BG =\begin{bmatrix} B & G \end{bmatrix},\qquad \CH =\begin{bmatrix} C & H \end{bmatrix}.
    \end{equation}
    A similar technique works to find common solutions for a low-rank completion problem 
    \begin{gather}
        \min_{X\in \mathbb{R}^{m_2\times m_1}} \rank \begin{bmatrix}
        \EA & \FB \\
        X & C
        \end{bmatrix}, \tag{LRCP$+$Rows} \label{eq:lrcprows} \\
        \EA = \begin{bmatrix} E\\ A \end{bmatrix},\qquad   \FB = \begin{bmatrix}  F \\ B \end{bmatrix} \label{eq:EAFB}.
    \end{gather}
    and the original low-rank completion problem \eqref{eq:lrcp}, which amounts to a removal of rows from \eqref{eq:lrcprows}.
    This result is stated and derived in Section~\ref{sec:addcolremrows}.
    
    \item \textit{Lemma \ref{lem:equivhankelproblem}}. Our solution to the problems \eqref{eq:lrcprows} and \eqref{eq:lrcpcols} shall use the additional assumption that $R(B^T) \cap R(F^T) = \{0\}$ and $R(B) \cap R(G) = \{0\}$, respectively. In Lemma \ref{lem:equivhankelproblem} of Section~\ref{sec:hankel_block_modification}, we show that we can modify the Hankel blocks of the matrix $A$ such that every adjacent pair of Hankel blocks $\mathcal{H}_k(X)$ and $\mathcal{H}_{k+1}(X)$ involves a removal of rows followed by an addition of columns satisfying the hypotheses of Lemmas~\ref{lem:addition_of_cols} and \ref{lem:removal_of_rows}.


\end{enumerate}

\subsubsection{Solution set of the two-by-two low-rank completion problem} \label{sec:2x2_low_rank_completion}

Consider the two-by-two low-rank completion problem \eqref{eq:lrcp} whose solution set is denoted by
%
%
%
    \begin{equation}
    \mathscr{S}_{A,B,C} := \left\{ \hat{X}\in\mathbb{R}^{m_2\times m_1}:\quad \rank \twobytwo{A}{B}{\hat{X}}{C} = \min_{X\in \mathbb{R}^{m_2\times m_1}} \rank \twobytwo{A}{B}{X}{C}  \right\}.  \label{eq:lrcp_solution_set}
\end{equation}
The complete solution of a generalization of this problem was derived in \cite{Woe87,KW88,Woe89}. An alternate construction yielding some solutions based on rank factorizations is given in \cite{EGH14}, but not all solutions are provided (nor is this claimed). Here, we provide characterization to the complete solution set in the same spirit as \cite{EGH14} by means of rank factorizations and intersection of subspaces. We emphasize that none of the results in Section~\ref{sec:2x2_low_rank_completion} are new. \ethancheck

Defining $r_{\rm opt} := \min_{X\in \mathbb{R}^{m_2\times m_1}} \rank M(X)$, we shall see that the solution to \eqref{eq:lrcp} can be characterized in terms of the column spaces of $A$ and $B$, and the row spaces of $B$ and $C$. For this reason, we define the spaces $V_{AB} := R(A) \cap R(B)$ and $W_{BC} = R(B^T) \cap R(C^T)$. We then choose complementary subspaces $V_{\bar{A}}$, $V_{\bar{B}}$, $W_{\bar{B}}$, and $W_{\bar{C}}$ satisfying
\begin{equation}\label{eq:complementary_subspaces}
    \begin{split}
    &R(A) = V_{AB} \oplus V_{\bar{A}}, \: R(B) = V_{AB} \oplus V_{\bar{B}},\\
    &R(B^T) = W_{BC} \oplus W_{\bar{B}}, \: R(C^T) = W_{BC} \oplus W_{\bar{C}}.\ethancheck
    \end{split}
\end{equation}
It follows from these definitions that $R(B^T) + R(C^T) = W_{\bar{B}} \oplus W_{BC} \oplus W_{\bar{C}}$ and $R(A) + R(B) = V_{\bar{A}} \oplus V_{AB} \oplus V_{\bar{B}}$. As a result, we see
\begin{equation*}
    \rank \onebytwo{A}{B} = \dim V_{\bar{A}} + \dim V_{AB} + \dim V_{\bar{B}},\quad  \rank(B) = \dim V_{AB} + \dim V_{\bar{B}},\quad  \textrm{ etc.} \ethancheck
\end{equation*}
The following lower bound was established in \cite{woerdeman1989minimal,eidelman2014separable}, which we reproduce here for sake of completeness.

\begin{proposition}\label{prop:low_rank_completion_lower_bound}
For every $X\in\mathbb{R}^{m_2\times m_1}$, we have:
$$\rank \left( \twobytwo{A}{B}{X}{C} \right) \ge \rank \begin{bmatrix} A & B \end{bmatrix} + \rank \begin{bmatrix} B \\ C \end{bmatrix} - \rank(B) = \rank(B) + \dim V_{\bar{A}} + \dim W_{\bar{C}}.$$
\end{proposition}
\begin{proof}

Set $s := \rank \left( \onebytwo{A}{B} \right)$ and let $\onebytwo{x_1^T}{y_1^T},\onebytwo{x_2^T}{y_2^T},\ldots,\onebytwo{x_{s}^T}{y_{s}^T}$ be a row basis for $\onebytwo{A}{B}$. 
Now extend the row basis for $\onebytwo{A}{B}$ to a row basis for $\twobytwo{A}{B}{X}{C}$ by adding rows
$$\onebytwo{x_{s+1}^T}{y_{s+1}^T},\ldots,\onebytwo{x_{r}^T}{y_{r}^T},\ethancheck$$
where $r = \rank \left( \twobytwo{A}{B}{X}{C} \right)$. Then $y_1^T,\ldots,y_{r}^T$ must span the row space of $C$. Since $y_{1}^T,\ldots,y_{s}^T$ all lie in the row space of $B$ which shares no vectors in common with $W_{\bar{C}}$, there must be at least $\dim W_{\bar{C}}$ vectors in $y_{s+1}^T,\ldots,y_{r}^T$. Thus,
\begin{equation*}
    \rank \left(\twobytwo{A}{B}{X}{C}\right) - \rank \left(\onebytwo{A}{B}\right) = r - s  \ge \dim W_{\bar{C}} = \rank \left(\twobyone{B}{C}\right) - \rank(B). \ethancheck
\end{equation*}
\end{proof}
In fact, the bound in Proposition \ref{prop:low_rank_completion_lower_bound} can be obtained by judicious choice of $X$. \ethancheck
\begin{alg}[Construction of $2\times 2$ Low-Rank Completion Problem Solution]\label{alg:2x2lrcp} Consider the $2\times 2$ low-rank completion problem  \eqref{eq:lrcp}. 

\begin{enumerate}
    \item \label{step:1alg2}Let $\PAbar$, $\PAB$, and $\PBbar $ be bases for $V_{\bar{A}}$, $V_{AB}$, and $V_{\bar{B}}$ respectively, and similarly let $\QBbar$, $\QBC$,  and $\QCbar$ be bases for $W_{\bar{B}},W_{BC},W_{\bar{C}}$. (Since $V_{\bar{A}}$, $V_{\bar{B}}$, $W_{\bar{B}}$, and $W_{\bar{C}}$ are non-unique, this requires also choosing such complementary subspaces satisfying \eqref{eq:complementary_subspaces}.)
    \item \label{step:2alg2}Conclude $\onebytwo{\PAbar}{\PAB}$ is a column basis for $A$, so there exists a matrix $Q^T_A$ such that
    \begin{equation*}
        A = \onebytwo{\PAbar}{\PAB} Q^T_A = \onebytwo{\PAbar}{\PAB} \twobyone{Q_{\bar{A},A}^T}{Q_{AB,A}^T}. \ethancheck 
    \end{equation*}
    Likewise, we may factor $B$ and $C$ as
    \begin{equation}\label{eq:factor_B_and_C}
        B = \onebytwo{\PAB}{\PBbar } \twobyone{Q_{AB,B}^T}{Q_{\bar{B},B}^T} = \onebytwo{\PBBbar}{\PBBC}\twobyone{Q_{\bar{B}}^T}{Q_{BC}^T}, \quad C = \onebytwo{\PCBC}{\PCCbar }\twobyone{Q_{BC}^T}{Q_{\bar{C}}^T}.\ethancheck
    \end{equation}
    \item\label{step:3alg2} Note we have now have two rank factorizations for $B$, which much necessarily be related by a nonsingular matrix $R\in \real^{\rank(B)\times \rank(B)}$
    \begin{gather}
        B = \onebytwo{\PAB}{\PBbar } R^{-1} \twobyone{Q_{\bar{B}}^T}{Q_{BC}^T} =   \onebytwo{\PBBbar}{\PBBC}R \twobyone{Q_{AB,B}^T}{Q_{\bar{B},B}^T}, \label{eq:B_rank_factorization}
    \end{gather}
    Conformally partition $R$ as
    \begin{gather}
           R = \twobytwo{\RBbarAB}{\RBbarBbar }{\RBCAB }{\RBCBbar }. \nonumber \ethancheck
    \end{gather}
    \item Define the solution set to be
    \begin{subequations}
    \begin{gather}
    \mathscr{S}_{A,B,C} = \left\{ X^{A,B,C}_{\FAbar,\FCbar} : \FAbar, \FCbar \mbox{ free} \right\} \\  X^{A,B,C}_{\FAbar ,\FCbar } = \FAbar \QAbarA + \PCBC \RBCAB \QABA  + \PCCbar \FCbar.
    \end{gather} \label{eq:2x2solutionset}
    \end{subequations}
\end{enumerate}

\end{alg}
\begin{lemma}\label{thm:2x2_low_rank_completion} 
The affine set \eqref{eq:2x2solutionset} computed in Algorithm~\ref{alg:2x2lrcp} describes the full solution set $\mathscr{S}_{A,B,C}$, as defined in \eqref{eq:lrcp_solution_set}, for the low-rank completion problem \eqref{eq:lrcp}.
\end{lemma}

\begin{proof}
Let $\FAbar \in \real^{m_2\times \dim V_{\bar{A}}}$ and $\FCbar  \in \real^{\dim W_{\bar{C}}\times n_1}$ be arbitrary matrices and consider the matrix $M_{\FAbar ,\FCbar }$ defined by
\begin{equation*}
    M_{\FAbar ,\FCbar } = \begin{bmatrix}
    \PAbar & \PAB & \PBbar  & 0 \\
    \FAbar  & \PCBC \RBCAB  & \PCBC \RBCBbar  & \PCCbar 
    \end{bmatrix}\begin{bmatrix}
    \QAbarA & 0 \\
    \QABA  & \QABB  \\
    0 & \QBbarB \\
    \FCbar  & \QCbar 
    \end{bmatrix}.
\end{equation*}
Then we observe that $M_{\FAbar ,\FCbar }$ has been written as the product of two full rank matrices and consequently $\rank (M_{\FAbar ,\FCbar }) = \rank \onebytwo{A}{B} + \rank \twobyone{B}{C} - \rank B$. Moreover, carrying out the matrix multiplication, we see that
\begin{equation}\label{eq:lowrankfact}
    M_{\FAbar ,\FCbar } = \twobytwo{A}{B}{X_{\FAbar ,\FCbar }^{A,B,C}}{C}, \quad X^{A,B,C}_{\FAbar ,\FCbar } = \FAbar \QAbarA + \PCBC \RBCAB \QABA  + \PCCbar \FCbar .\ethancheck 
\end{equation}
Thus $\mathscr{S}_{A,B,C} \supseteq \{ X^{A,B,C}_{\FAbar,\FCbar} : \FAbar, \FCbar \mbox{ free} \}$. Conversely, suppose $M = M(X)$ solves the low rank completion problem. Then $M$ has a rank factorization
\begin{gather*}
    M = \twobytwo{A}{B}{X}{C} = \twobyone{P_1}{P_2}\onebytwo{Q_1^T}{Q^T_2},\\
    \twobyone{P_1}{P_2} \in \real^{m_1+m_2 \times \ropt}, \: \onebytwo{Q_1^T}{Q_2^T} \in \real^{\ropt \times n_1+n_2}.
\end{gather*}
Our goal is to re-write the above into \eqref{eq:lowrankfact} through a sequence of invertible transformations. The columns of $P_1$ span $R(\onebytwo{A}{B})$, so there exists a non-singular matrix $S$ such that 
\begin{equation*}
    P_1S = \begin{bmatrix} \PAbar & \PAB & \PBbar  & P' \end{bmatrix},
\end{equation*}
where the columns in $P'$ are zero or do not lie in $R(A) + R(B)$.\ethancheck Partition $S^{-1}\onebytwo{Q_1^T}{Q_2^T}$ as
\begin{equation*}
    S^{-1}\onebytwo{Q_1^T}{Q_2^T} = \begin{bmatrix} G_{11} & G_{12} \\
    G_{21} & G_{22} \\ G_{31} & G_{32} \\ G_{41} & G_{42} \end{bmatrix}.\ethancheck 
\end{equation*}
Then we have
\begin{align*}
    A &= \PAbar G_{11} + \PAB G_{21} + \PBbar  G_{31} + P'G_{41}, \\
    B &= \PAbar G_{12} + \PAB G_{22} + \PBbar G_{32} + P'G_{42}.\ethancheck 
\end{align*}
Then $A - \PAbar G_{11} - \PAB G_{21} = \PBbar G_{31} + P'G_{41}$. Since the columns of $A - \PAbar G_{11} - \PAB G_{21}$ lie in $R(A)$ and the columns of $\PBbar G_{31} + P'G_{41}$ lie in a complement of $R(A)$, we must have $\PBbar G_{31} + P'G_{41} = 0$. Writing $\PBbar G_{31} = -P'G_{41}$ and using the same technique, we deduce $\PBbar G_{31} = P'G_{41} = 0$. Similarly, we may see that $\PAbar G_{12} = P'G_{42} = 0$. Since $\PAbar$ and $\PBbar$ have full column rank, $G_{12} = 0$ and $G_{31} = 0$. Thus, we have
\begin{align*}
    A &= \onebytwo{P_{\bar A}}{\PAB}\twobyone{G_{11}}{G_{21}} \implies \twobyone{G_{11}}{G_{21}} = \twobyone{\QAbarA}{\QABA },\\
    B &= \onebytwo{\PAB}{P_{\bar B}}\twobyone{G_{22}}{G_{32}} \implies \twobyone{G_{22}}{G_{32}} = \twobyone{\QABB  }{\QBbarB}.\ethancheck 
\end{align*}
In summary, we have shown
\begin{equation*}
    S^{-1}\onebytwo{Q_1^T}{Q_2^T} = \begin{bmatrix} \QAbarA & 0 \\
    \QABA  & \QABB   \\ 0 & \QBbarB \\ G_{41} & G_{42} \end{bmatrix}.\ethancheck 
\end{equation*}
We have that $P_2Q_2^T = (P_2S)(S^{-1}Q_2^T) = C$. Consequently, we see that $\begin{bmatrix} Q_{AB,B} & Q_{\bar{B},B} & G_{42}^T \end{bmatrix}^T$ must be a row basis for $\twobyone{B}{C}$. Thus there exists a nonsingular transformation $T$ such that
\begin{equation*}
    T = \begin{bmatrix}
    I & 0 & 0 & 0 \\
    0 & I & 0 & 0 \\
    0 & 0 & I & 0 \\
    0 & T_1 & T_2 & T_3
    \end{bmatrix}, \quad T^{-1}S^{-1}\onebytwo{Q_1^T}{Q_2^T} = \begin{bmatrix} \QAbarA & 0 \\
    \QABA  & \QABB   \\ 0 & \QBbarB \\ \FCbar  & \QCbar  \end{bmatrix}.\ethancheck 
\end{equation*}
Then we have
\begin{equation*}
    \twobyone{P_1}{P_2}ST = \begin{bmatrix} \PAbar & \PAB + P'T_1 & \PBbar + P'T_2  & P'T_3 \\
    \FAbar  & H_2 & H_3 & H_4\end{bmatrix}.\ethancheck 
\end{equation*}
Then seeing that 
\begin{align*}
    C &= (P_2ST)(T^{-1}S^{-1}Q_2^T) = \begin{bmatrix} H_2 & H_3 & H_4 \end{bmatrix} \begin{bmatrix} \QABB   \\ \QBbarB \\ \QCbar  \end{bmatrix} \\
    &= \begin{bmatrix} \PCBC \RBCAB  & \PCBC \RBCBbar  & \PCCbar  \end{bmatrix}\begin{bmatrix} \QABB   \\ \QBbarB \\ \QCbar  \end{bmatrix},
\end{align*}
we conclude that $H_2 = \PCBC \RBCAB $, $H_3 = \PCBC \RBCBbar $, and $H_4=\PCCbar $ because $\begin{bmatrix} \QABB   \\ \QBbarB \\ \QCbar  \end{bmatrix}$ has full row rank. We also have $B = \PAB\QABB   + \PBbar \QBbarB + P'(T_1\QABB + T_2\QBbarB + T_3\QCbar) = B + P'(T_1\QABB + T_2\QBbarB + T_3\QCbar)$. Thus
\begin{equation*}
    P' \begin{bmatrix} T_1 & T_2 & T_3 \end{bmatrix} \begin{bmatrix} \QABB \\ \QBbarB \\ \QCbar \end{bmatrix} = 0 \implies P' \begin{bmatrix} T_1 & T_2 & T_3 \end{bmatrix} = 0.\ethancheck 
\end{equation*}
Thus $P'T_3 = 0$ so $P' = 0$ since $T_3$ is invertible. We have thus shown that
\begin{align*}
    M &= \twobyone{P_1}{P_2}\onebytwo{Q_1^T}{Q_2^T} = \twobyone{P_1}{P_2}ST(T^{-1}S^{-1})\onebytwo{Q_1^T}{Q_2^T} = \\
    &= \begin{bmatrix}
    \PAbar & \PAB & \PBbar  & 0 \\
    \FAbar  & \PCBC \RBCAB  & \PCBC \RBCBbar  & \PCCbar 
    \end{bmatrix}\begin{bmatrix}
    \QAbarA & 0 \\
    \QABA  & \QABB  \\
    0 & \QBbarB \\
    \FCbar  & \QCbar 
    \end{bmatrix} = M_{\FAbar ,\FCbar }.\ethancheck 
\end{align*}
Thus, the class of matrices $M_{\FAbar ,\FCbar }$ completely characterize the solution set of the low rank completion problem. \qedhere 
%
%

%
%
\end{proof}
\notrunc{
\begin{example}
Consider the low rank completion problem
$$ \begin{bmatrix}
        A & B \\
        ? & C
        \end{bmatrix} =  \begin{bmatrix}
        \begin{bmatrix}
        1 & 0 \\
        0 & 1
        \end{bmatrix}  & \begin{bmatrix}
        0 & 0 \\
        1 & 1
        \end{bmatrix}  \\
        \begin{bmatrix}
        ? & ? \\
        ? & ?
        \end{bmatrix}  & \begin{bmatrix}
        1 & 1 \\
        1 & 0
        \end{bmatrix} 
        \end{bmatrix}. $$
The solution to this particular completion problem is
$$  \begin{bmatrix}
        ? & ? \\
        ? & ?
        \end{bmatrix} = \begin{bmatrix}
        \alpha & 1 \\
        \beta & \gamma
        \end{bmatrix}, \quad \alpha,\beta,\gamma \in \mathbb{R}.  $$
Indeed, \eqref{eq:lowrankfact} generates this solution. We compute
\begin{align*}
A &= \onebytwo{\PAbar}{\PAB} \twobyone{Q_{\bar{A},A}^T}{Q_{AB,A}^T} = \onebytwo{ \twobyone{1}{0}      }{  \twobyone{0}{1}  }  \twobyone{ \onebytwo{1}{0}  }{ \onebytwo{0}{1} }, \\
C &= \onebytwo{\PCBC}{\PCCbar }\twobyone{Q_{BC}^T}{Q_{\bar{C}}^T} = \onebytwo{ \twobyone{1}{0}      }{  \twobyone{0}{1}  }  \twobyone{ \onebytwo{1}{1}  }{ \onebytwo{1}{0} }, \\
B & = \onebytwo{\PAB}{\PBbar } \twobyone{Q_{AB,B}^T}{Q_{\bar{B},B}^T} = \onebytwo{ \twobyone{0}{1}      }{  \twobyone{}{}  }  \twobyone{ \onebytwo{1}{1}  }{ \onebytwo{}{} } \\
& =  \onebytwo{\PBBbar}{\PBBC}\twobyone{Q_{\bar{B}}^T}{Q_{BC}^T} =  \onebytwo{  \twobyone{}{}    }{   \twobyone{0}{1}    }  \twobyone{ \onebytwo{}{} }{   \onebytwo{1}{1}  }, \\
R &= \twobytwo{\RBbarAB}{\RBbarBbar }{\RBCAB }{\RBCBbar } = \twobytwo{()}{()}{(1)}{()}.
  \end{align*}
This yields
\begin{eqnarray*} 
 X_{\FAbar ,\FCbar } & = &  \FAbar  \QAbarA + \PCBC\RBCAB \QABA  + \PCCbar \FCbar  \\
             & = & \twobyone{f_{11}}{ f_{12} } \onebytwo{1}{0}   +  \twobyone{1}{0} (1) \onebytwo{1}{1}   + \twobyone{0}{1} \onebytwo{f_{21}}{ f_{22} }  \\
             & = & \twobytwo{f_{11}+1}{1}{f_{12}+f_{21}}{f_{22}}.
 \end{eqnarray*}
Setting $f_{11} := \alpha - 1$, $f_{12} := \beta$, $f_{21} := 0$, and $f_{22} := \gamma$, we recover all the solutions. We observe that a unique choice of the free variables does not generate a unique solution, as any pair of assignments of the free variables $f_{12}$ and $f_{21}$ respecting the constraint $f_{12} + f_{21} = \beta$ will generate the same solution. \ethancheck 
\end{example}
}
\subsubsection{The solution set subjected to addition of columns and removal of rows} \label{sec:addcolremrows}

Consider the lower triangular portion of the $5\times 5$ CSS construction. We seek to pick $X$ minimizing the rank of the Hankel blocks of 
\begin{equation*}
\begin{bmatrix}
A_{21} &&& \\
A_{31} & A_{32} && \\
A_{41} & A_{42} & A_{43} &\\
X & A_{52} & A_{53} & A_{54}
\end{bmatrix}. \ethancheck
\end{equation*}
%
%
%
%
%
The ``skinniest'' and ``fattest'' Hankel blocks shall pose no major problem as minimizing their ranks will just require that $X$'s rows and columns lie in certain spaces, $R(X) \subseteq R (\begin{bmatrix}
A_{52} & A_{53} & A_{54}
\end{bmatrix})$ and $R(X^T) \subseteq R (\begin{bmatrix}
A_{21} & A_{31} & A_{41}
\end{bmatrix}^T)$. The difficulty will be in minimizing the ranks of the intermediate Hankel blocks
\begin{equation*}
\begin{bmatrix}
A_{31} & A_{32} \\
A_{41} & A_{42} \\
X & A_{52}
\end{bmatrix} \textrm{ and } \begin{bmatrix}
A_{41} & A_{42} & A_{43}\\
X & A_{52} & A_{53}
\end{bmatrix}.\ethancheck
\end{equation*}
Algorithm~\ref{alg:2x2lrcp} already enables us to find solutions for  the low-rank completion problems \eqref{eq:lrcp} for the matrices
\begin{equation*}
\begin{bmatrix}
\twobyone{A_{31}}{A_{41}} & \twobyone{A_{32}}{A_{42}} \\
X & A_{52}
\end{bmatrix}
, \quad  \begin{bmatrix}
A_{41} & A_{42} \\
X & A_{52}
\end{bmatrix}, \mbox{ and } \begin{bmatrix}
A_{41} & \onebytwo{A_{42}}{A_{43}} \\
X & \onebytwo{A_{52}}{A_{53}}
\end{bmatrix}. \ethancheck
\end{equation*}
The question is to what extent we can find \textit{common solutions} of these low rank completion problems. To address this question, we need to tackle the following problems first.
%
\begin{problem}\label{prob:adding_cols}
	How does the solution set of the $2\times 2$ low-rank completion problem \eqref{eq:lrcp} change after an addition of columns to \eqref{eq:lrcpcols}?
\end{problem}
\begin{problem}\label{prob:removing_rows}
	How does the solution set of the $2\times 2$ low-rank completion problem \eqref{eq:lrcprows} change after a removal of rows to \eqref{eq:lrcp}?
\end{problem}

\begin{example}
Counterintuitively, adding columns can both contract the solution set:
\begin{equation*}
M = \twobytwo{1}{0}{x}{0} \textrm{ is solved by all } x\in \real, \: M' = \begin{bmatrix}
1 & 0 & 1 \\
x & 0 & 1
\end{bmatrix} \textrm{ solved only for } x = 1,
\end{equation*}
or expand the solution set:
\begin{equation*}
M = \twobytwo{0}{0}{x}{0} \textrm{ is solved only for } x=0, \: M' = \begin{bmatrix}
0 & 0 & 0\\
x & 0 & 1
\end{bmatrix} \textrm{ solved for all } x \in \real.\ethancheck
\end{equation*}
Similar statement can be made also for the removal of rows.
\end{example}
%
We note that while the above example shows that the solution set of a low rank completion can grow or shrink with the addition of columns, in both examples there remain \textit{common solutions}. We shall show that this is always the case. Our strategy will be as follows. We shall begin with the six bases $\PAbar$, $\PAB $, $\PBbar $, $\QBbar$, $\QBC $, and $\QCbar$ which define the solution set of \eqref{eq:lrcp} by \eqref{eq:2x2solutionset}. We shall then perform an algebraic calculation to compute six bases $\PAbarbullet $, $\PAG $, $\PGbar $, $\QGbar $, $\QGH$, and $\QHbar $ associated with the new problem \eqref{eq:lrcpcols}. From here we shall compare the solution sets of the two problems.  \ethancheck

\begin{remark}
When there is a naming conflict such as needing to have a separate ``$\PAbar$'' for both the original low rank completion problem \eqref{eq:lrcp} and the new problem \eqref{eq:lrcpcols}, we shall use a $\bullet$ to denote the named object corresponding to the new problem \eqref{eq:lrcpcols}. For instance, we will have $R(A) = R(\PAbar) \oplus R(\PAB) = R(\PAbarbullet) \oplus R(\PAG)$, where $\PAG$ is a basis for $R(A) \cap R\left(\twobyone{B}{G}\right)$ and $\PAB$ is a basis for $R(A) \cap R(B)$.
\end{remark}
    
\begin{lemma}[Addition of Columns]\label{lem:addition_of_cols}
    Consider a subset 
    \begin{equation*}
        \mathscr{S}_{A,B,C}^Y := \left\{ \FAbar\QAbarA + \PCBC \RBCAB \QABA + \PCCbar Y : \FAbar \textnormal{ free} \right\},
    \end{equation*}
    of the solution set \eqref{eq:2x2solutionset} of the low rank completion problem \eqref{eq:lrcp} computed in Algorithm~\ref{alg:2x2lrcp} where the second free variable ``$\FCbar$'' is set to an arbitrary fixed value $Y$.
    Suppose $R(B)\cap R(G)= \{0\}$. Then there exists a special choice of the row and column bases $\PAbarbullet$, $\PAG$, $\PGbar$, $\QGbar$, $\QGH$, and $\QHbar$ in Step~\ref{step:1alg2} of Algorithm~\ref{alg:2x2lrcp} applied to \eqref{eq:lrcpcols}
    and a matrix $Z$ depending on $Y$ such that
    \begin{equation*}
        \mathscr{S}_{A,\BG,\CH}^Z := \left\{ \FAbarbullet\QAbarbulletA + \PHGH \RGHAG \QAGA + \PHHbar Z : \FAbarbullet \textnormal{ free} \right\} \subseteq \mathscr{S}_{A,B,C}^Y.
    \end{equation*}
\end{lemma}
Lemma~\ref{lem:addition_of_cols} shows that there exists a nonempty set $\mathscr{S}_{A,\BG,\CH}^Z$ which simultaneously minimizes the rank of \eqref{eq:lrcp} and \eqref{eq:lrcpcols} with the ``$\FCbar$'' free variable in the solution set of \eqref{eq:lrcp} is set to the particular value of $Y$. In preparation of proving Lemma~\ref{lem:addition_of_cols}, we need the following technical result.
\begin{proposition}\label{prop:extend_to_row_basis}
Let $\twobytwo{\Qs{11}}{\Qs{12}}{0}{\Qs{22}}$ be a linearly independent collection of rows in the row space of a matrix $\onebytwo{B}{G}$ and let $\Qs{31}$ be a linearly independent collection of rows such that $\twobyone{\Qs{11}}{\Qs{31}}$ forms a row basis for the matrix $B$. Then there exists a row basis of the form
\begin{equation*}
    \begin{bmatrix}
    \Qs{11} & \Qs{12} \\
    0 & \Qs{22} \\
    \Qs{31} & \Qs{32} \\
    0 & \Qs{42}
    \end{bmatrix}
\end{equation*}
for $\onebytwo{B}{G}$.
\end{proposition}

\begin{proof}
Let $V^T = \onebytwo{V_1^T}{V_2^T}$ be a row basis for $\onebytwo{B}{G}$. Then since $\twobytwo{\Qs{11}}{\Qs{12}}{0}{\Qs{22}}$ is a linearly independent collection of rows spanned by $V^T$, there exists a nonsingular matrix $S$ such that
\begin{equation*}
    SV^T = \begin{bmatrix}
    \Qs{11} & \Qs{12} \\
    0 & \Qs{22} \\
    V_{31}^T & V_{32}^T
    \end{bmatrix}.\ethancheck
\end{equation*}
Then, since $\twobyone{\Qs{11}}{\Qs{13}}$ is a basis for $B$, there exists a nonsingular block triangular matrix $T$ such that
\begin{equation*}
    T = \begin{bmatrix}
    I & 0 & 0 \\
    0 & I & 0 \\
    T_{31} & T_{32} & T_{33} \\
    T_{41} & T_{42} & T_{43}
    \end{bmatrix}, \quad W^T := TSV^T = \begin{bmatrix}
    \Qs{11} & \Qs{12} \\
    0 & \Qs{22} \\
    \Qs{31} & \Qs{32} \\
    \Qs{\ast} & \Qs{42}
    \end{bmatrix},
\end{equation*}
where $\Qs{\ast}$ has trivial row intersection with $\twobyone{\Qs{11}}{\Qs{31}}$ (and thus $B$).\ethancheck Since $TS$ is nonsingular, $W^T$ is a row basis for $\onebytwo{B}{G}$ so there exists a column basis $U = \begin{bmatrix} U_1 & U_2 & U_3 & U_4 \end{bmatrix}$ such that
\begin{equation*}
    \onebytwo{B}{G} = UW^T = \begin{bmatrix} U_1 & U_2 & U_3 & U_4 \end{bmatrix} \begin{bmatrix}
    \Qs{11} & \Qs{12} \\
    0 & \Qs{22} \\
    \Qs{31} & \Qs{32} \\
    \Qs{\ast} & \Qs{42}
    \end{bmatrix}.\ethancheck
\end{equation*}
Thus $B = U_1\Qs{11} + U_3 \Qs{31} + U_4 \Qs{\ast}$ so $B - U_1\Qs{11} - U_3 \Qs{21} = U_4\Qs{\ast}$. Since the left- and right-hand sides of this equation lie row spaces with trivial intersection, we conclude $B - U_1\Qs{11} - U_3 \Qs{21} = U_4\Qs{\ast} = 0$. Since $U_3$ has full column rank, $\Qs{\ast} = 0$ and $W^T$ has the desired structure.\ethancheck
\end{proof}

\begin{proof}[Proof of Lemma~\ref{lem:addition_of_cols}]

Suppose we have computed row and column bases $\PAbar$, $\QBC$, etc. given by Algorithm \ref{alg:2x2lrcp} applied to \eqref{eq:lrcp}. We shall now construct the six bases $\QGH$ \eqref{eq:QGH}, $\QGbar$ \eqref{eq:QGbar}, $\QHbar$ \eqref{eq:QHbar}, $\PAbarbullet$ \eqref{eq:PAbarbullet}, $\PAG$ \eqref{eq:PAG}, and $\PGbar$ \eqref{eq:PAG} and the corresponding complementary bases $\PGGbar$ \eqref{eq:Bcolumnbases}, $\PGGH$ \eqref{eq:Bcolumnbases}, $\PHGH$ \eqref{eq:Ccolumbasis}, $\PHHbar$ \eqref{eq:Ccolumbasis}, $\QAbarbulletA$ \eqref{eq:Arowbasis}, $\QAGA$ \eqref{eq:Arowbasis}, $\QAGG$ \eqref{eq:BGrowbases}, and $\QGbarG$ \eqref{eq:BGrowbases} in Step~\ref{step:1alg2} of Algorithm~\ref{alg:2x2lrcp}. For convenience, we list an equation reference where each of these bases is defined and demarcate different stages of the construction by bolded subheadings. 

\proofmarker{Construction of $\QGH$.} We choose a row basis
\begin{equation}\label{eq:QGH}
\QGH = \begin{bmatrix}\QBCone  & \QGHone  \\
0 & \QGHtwo  \end{bmatrix} \textrm{ for } R\left(\onebytwo{B}{G}^T\right) \cap R\left(\onebytwo{C}{H}^T\right)
\end{equation}
such that $\QBCone $ has full row rank.\ethancheck (To get such a basis, simply choose any basis for this space and apply $QR$.)

\proofmarker{Construction of $\QGbar$.} Extend $\QBCone$ to a basis $\twobyone{\QBCone }{\QBCtwo  }$ for the intersection of the row spaces of $B$ and $C$. Then, since any two row bases are related by left multiplication by a nonsingular matrix, there exists a nonsingular matrix $K$ such that
\begin{equation}\label{eq:splitting_QTBC}
K\QBC  = \twobyone{K_1}{K_2}\QBC  = \twobyone{\QBCone }{\QBCtwo  }.\ethancheck
\end{equation}
Then, we extend $\begin{bmatrix}\QBCone  & \QGHone  \\
0 & \QGHtwo  \end{bmatrix}$ to the following basis of $R\onebytwo{B}{G}^T$ using Proposition \ref{prop:extend_to_row_basis}:
\begin{equation}\label{eq:QGbar}
\left[\arraycolsep=1.4pt\def\arraystretch{1.4}\begin{array}{cc}
\QBCtwo   & \QGbarone \\
\QBbar & \QGbartwo  \\
0 & \QGbarthree  \\\hline
\QBCone  & \QGHone  \\
0 & \QGHtwo 
\end{array} \right] = \arraycolsep=1.4pt\def\arraystretch{1.4}\twobyone{\QGbar }{\QGH}. \ethancheck
\end{equation}
%
\proofmarker{Construction of $\PGGH$ and $\PGGbar$.} Let $\onebytwo{\PGGbar }{\PGGH } = \onebytwo{B}{G}\left( \twobyone{\QGbar }{\QGH} \right)^\dagger$ be the complementary column basis of $\onebytwo{B}{G}$ satisfying
\begin{equation}\label{eq:Bcolumnbases}
\begin{split}
    \onebytwo{B}{G} &= \onebytwo{\PGGbar }{\PGGH }\twobyone{\QGbar }{\QGH} \\
    &= \left[\begin{array}{ccc|cc}
    \PBBCtwo  & \PGGbartwo  & \PGGbarthree   & \PBBCone  & \PGGHtwo  
    \end{array}\right] \left[\arraycolsep=1.4pt\def\arraystretch{1.4}\begin{array}{cc}
        \QBCtwo   & \QGbarone \\
        \QBbar & \QGbartwo  \\
        0 & \QGbarthree  \\\hline
        \QBCone  & \QGHone  \\
        0 & \QGHtwo 
        \end{array}\right].\ethancheck
\end{split}
\end{equation}
(Here, $(\cdot)^\dagger$ denotes the matrix pseudoinverse.) Then 
\begin{equation*}
    B = \begin{bmatrix}
    \PBBCone  & \PBBCtwo  & \PGGbartwo   
    \end{bmatrix}\begin{bmatrix}
    \QBCone  \\ \QBCtwo   \\ \QBbar
    \end{bmatrix} = \begin{bmatrix}
    \begin{bmatrix} \PBBCone  & \PBBCtwo \end{bmatrix}K & \PGGbartwo   
    \end{bmatrix}\begin{bmatrix}
    \QBC  \\ \QBbar
    \end{bmatrix},
\end{equation*}
so $\PGGbartwo  =\PBBbar $ and $\begin{bmatrix} \PBBCone  & \PBBCtwo \end{bmatrix} = \PBBC K^{-1}$ by comparison with \eqref{eq:factor_B_and_C} and the full row rank of $\begin{bmatrix} \QBC \\ \QBbar \end{bmatrix}$. \notrunc{(See \eqref{eq:splitting_QTBC} for the definition of $K$.)}\ethancheck 

\proofmarker{Construction of $\QHbar$.} Similar to the construction of $\QGbar$, we use Proposition \ref{prop:extend_to_row_basis} to construct a basis of $R(\onebytwo{C}{H}^T)$:
\begin{equation}\label{eq:QHbar}
\left[\arraycolsep=1.4pt\def\arraystretch{1.4}\begin{array}{cc}
\QBCtwo   & \QHbarone  \\
\QCbar & \QHbartwo  \\
0 & \QHbarthree \\ \hline
\QBCone  & \QGHone  \\
0 & \QGHtwo 
\end{array}\right] = \twobyone{\QHbar }{\QGH}.\ethancheck
\end{equation}
\proofmarker{Construction of $\PHGH$ and $\PHHbar$.}In the same way as $\onebytwo{B}{G}$, we deduce that $\onebytwo{C}{H}$ has a column basis of the form 
\begin{equation}\label{eq:Ccolumbasis}
    \onebytwo{\PHHbar }{\PHGH }=\left[\begin{array}{ccc|cc} \PCBCtwo  & \PCCbar  & \PHHbarthree & \PCBCone  & \PHGHtwo \end{array}\right]
\end{equation}
such that 
\begin{equation*}
    \onebytwo{C}{H} = \onebytwo{\PHHbar }{\PHGH }\twobyone{\QHbar }{\QGH}. \ethancheck
\end{equation*}
%
\proofmarker{Construction of $\PAbarbullet$.} Using our assumption that $R(B)\cap R(G)=\{0\}$, $R(G)$ can be divided into two parts $R(G) = R(\PGcapA ) \oplus R(\PGcapAbar )$, where $R(\PGcapA ) \subseteq R(A)$ and $R(\PGcapAbar ) \cap (R(A) + R(B)) = \{0\}$. Since $\onebytwo{\PGcapA}{\PAB}$ is a collection of linearly independent columns of $A$, they may be extended to a basis
\begin{equation}\label{eq:PAbarbullet}
\begin{bmatrix}
\PAbarbullet & \PGcapA & \PAB
\end{bmatrix}
\end{equation}
for $R(A)$.

\proofmarker{Construction of $\PAG$ and $\PGbar$.} Since any two column bases can be obtained by left multiplication by a nonsingular matrix, there exists nonsingular transformation $L$ such that
\begin{equation*}
    \onebytwo{\PAbar}{\PAB }\underbrace{\twobytwo{L_{11}}{L_{21}}{L_{21}}{L_{22}}}_{=L} = \onebytwo{\onebytwo{\PAbarbullet }{\PGcapA }}{\PAB }, \quad L = \twobytwo{L_{11}}{L_{21}}{L_{21}}{L_{22}} = \twobytwo{L_{11}}{0}{L_{21}}{I}.\ethancheck
\end{equation*}
$L$ must have the stated block lower triangular structure since
\begin{equation*}
    \onebytwo{\PAbar}{\PAB} \twobyone{L_{21}}{L_{22}} = \onebytwo{\PAbar}{\PAB} \twobyone{0}{I} = \PAB \mbox{ which implies } \twobyone{L_{21}}{L_{22}} = \twobyone{0}{I} 
\end{equation*}
since $\onebytwo{\PAbar}{\PAB}$ has full column rank.
Thus, defining
\begin{equation}\label{eq:PAG}
    \PGbar  := \onebytwo{\PBbar }{\PGcapAbar },\quad \PAG  := \onebytwo{\PAB }{\PGcapA},
\end{equation}
we have that $\PAG $ is a basis for $R(A)\cap R(\onebytwo{B}{G})$ and $R(A) = R(\PAG ) \oplus R(\PAbarbullet )$ and $R(\onebytwo{B}{G}) = R(\PAG ) \oplus R(\PGbar )$.

\proofmarker{Construction of $\QAbarbulletA$ and $\QAGA$.} We have
\begin{equation}\label{eq:Arowbasis}
\begin{split}
    A &= \onebytwo{\PAbar}{\PAB }\twobyone{\QAbarA }{\QABA } = \onebytwo{\PAbar}{\PAB }LL^{-1}\twobyone{\QAbarA }{\QABA } \\
    &= \begin{bmatrix}
    \begin{bmatrix} \PAbarbullet  & \PGcapA \end{bmatrix}  & \PAB 
    \end{bmatrix}\twobytwo{L_{11}^{-1}}{0}{-L_{21}L_{11}^{-1}}{I}\begin{bmatrix}
    \QAbarA  \\ \QABA 
    \end{bmatrix} \\
    &= \begin{bmatrix}
    \begin{bmatrix} \PAbarbullet  & \PGcapA \end{bmatrix}  & \PAB 
    \end{bmatrix}\begin{bmatrix}
    L_{11}^{-1}\QAbarA  \\ \QABA -L_{21}L_{11}^{-1}\QAbarA
    \end{bmatrix} \\
    &= \begin{bmatrix} \PAbarbullet  & \PGcapA  & \PAB 
    \end{bmatrix}\begin{bmatrix}
    \QAbarbulletA  \\ \QGcapAA  \\ \QABbulletA
    \end{bmatrix} \\
    &= \left[\begin{array}{c|cc}
    \PAbarbullet  & \PAB & \PGcapA  
    \end{array}\right]\left[\arraycolsep=1.4pt\def\arraystretch{1.4}\begin{array}{c}
    \QAbarbulletA  \\\hline \QABbulletA \\ \QGcapAA  
    \end{array}\right] \\
    &= \begin{bmatrix}
    \PAbarbullet  & \PAG 
    \end{bmatrix}\begin{bmatrix} \QAbarbulletA  \\ \QAGA  \end{bmatrix},
\end{split}
\end{equation}
where
\begin{equation}\label{eq:QABbulletA}
\QABbulletA := \QABA - L_{21}L_{11}^{-1}\QAbarA = \QABA + J\QAbarA., \quad J = -L_{21}L_{11}^{-1}.\ethancheck   
\end{equation}
\proofmarker{Construction of $\QAGG$ and $\QGbarG$.} Then since $\onebytwo{\PGcapA }{\PGcapAbar }$ is a column basis for $G$, there exists full rank matrices $\QGcapAG $ and $\QGcapAbarG $ such that
\begin{equation*}
    G = \onebytwo{\PGcapA }{\PGcapAbar }\twobyone{\QGcapAG }{\QGcapAbarG }.\ethancheck
\end{equation*}
Then
\begin{equation}\label{eq:BGrowbases}
    \onebytwo{B}{G} = \left[\begin{array}{cc|cc}
    \PAB  & \PGcapA  & \PBbar  & \PGcapAbar 
    \end{array}\right] \left[\arraycolsep=1.4pt\def\arraystretch{1.4}\begin{array}{cc}
    \QABB  & 0 \\
    0 & \QGcapAG  \\ \hline
    \QBbarB  & 0 \\
    0 & \QGcapAbarG 
    \end{array}\right] = \onebytwo{\PAG }{\PGbar }\twobyone{\QAGG }{\QGbarG}.\ethancheck
\end{equation}
\proofmarker{Expression of the solution set of \eqref{eq:lrcpcols}.} We have now computed the six bases $\PAbarbullet $, $\PAG $, $\PGbar $, $\QGbar $, $\QGH$, and $\QHbar $ for the problem \eqref{eq:lrcpcols}, so following the calculations of Algorithm~\ref{alg:2x2lrcp}, the solution set of \eqref{eq:lrcpcols} is $\mathscr{S}_{A,\BG,\CH} := \{ X^{A,\BG,\CH}_{\FAbarbullet,\FHbar} : \FAbarbullet, \FHbar \mbox{ free} \}$, where
\begin{equation}\label{eq:additional_column_solutions}
    X^{A,\BG,\CH}_{\FAbarbullet,\FHbar} = \FAbarbullet\QAbarbulletA  + \PHGH \RGHAG \QAGA  + \PHHbar \FHbar,
\end{equation}
and $R'=\twobytwo{\RGbarAG}{\RGbarGbar}{\RGHAG}{\RGHGbar}$ is the unique nonsingular matrix satisfying 
\begin{align*}
    \onebytwo{B}{G} &= \onebytwo{\PGGbar }{\PGGH }\twobytwo{\RGbarAG}{\RGbarGbar}{\RGHAG}{\RGHGbar}\twobyone{\QAGG }{\QGbarG} \\
    &= \onebytwo{\begin{bmatrix} \PBBCtwo  & \PBBbar  & \PGGbarthree   \end{bmatrix}}{\begin{bmatrix} \PBBCone  & \PGGHtwo   \end{bmatrix}} \times \\
    &\quad \twobytwo{\begin{bmatrix}
    \RGbarAG^{11} & \starorval{\RGbarAG^{12}} \\
    \RGbarAG^{21} & \starorval{\RGbarAG^{22}} \\
    \RGbarAG^{31} & \starorval{\RGbarAG^{32}}
    \end{bmatrix}}{\begin{bmatrix}
    \RGbarGbar^{11} & \starorval{\RGbarGbar^{12}} \\
    \RGbarGbar^{21} & \starorval{\RGbarGbar^{22}} \\
    \RGbarGbar^{31} & \starorval{\RGbarGbar^{32}}
    \end{bmatrix}}{\begin{bmatrix}
    \RGHAG^{11} & \RGHAG^{12} \\
    \RGHAG^{21} & \RGHAG^{22}
    \end{bmatrix}}{\begin{bmatrix}
    \RGHGbar^{11} & \starorval{\RGHGbar^{12}} \\
    \RGHGbar^{21} & \starorval{\RGHGbar^{22}}
    \end{bmatrix}} \times\twobyone{\twobytwo{\QABB }{0}{0}{\QGcapAG }}{\twobytwo{\QBbarB }{0}{0}{\QGcapAbarG }}. \ethancheck
\end{align*}
(See \eqref{eq:Bcolumnbases} and \eqref{eq:BGrowbases} for expressions regarding the block partitioning of $\PGGbar$, $\PGGH$, $\QAGG$, and $\QGbarG$.) We use asterisks to denote matrices whose value shall be immaterial to the ensuing calculations.

\proofmarker{Derivation of a relation between $\RGHAG$ and $\RBCAB$.} We can write $B$ as
\begin{equation}\label{eq:first_B}
\begin{split}
    B &= \begin{bmatrix} \PBBCtwo  & \PBBbar  & \PGGbarthree   & \PBBCone  & \PGGHtwo   \end{bmatrix}\begin{bmatrix}
    \RGbarAG^{11} & \RGbarGbar^{11}  \\
    \RGbarAG^{21} & \RGbarGbar^{21}  \\
    \RGbarAG^{31} & \RGbarGbar^{31}  \\
    \RGHAG^{11} & \RGHGbar^{11}   \\
    \RGHAG^{21} & \RGHGbar^{21}  
    \end{bmatrix}\begin{bmatrix}
    \QABB  \\ \QBbarB 
    \end{bmatrix}\\
    &= \begin{bmatrix}
    \PBBbar & \PBBCone & \PBBCtwo 
    \end{bmatrix}\begin{bmatrix}
    \RGbarAG^{21} & \RGbarGbar^{21} \\
    \RGHAG^{11} & \RGHGbar^{11} \\
    \RGbarAG^{11} & \RGbarGbar^{11}
    \end{bmatrix}\twobyone{\QABB}{\QBbarB} \\
    &\qquad + \onebytwo{\PGGbarthree}{\PGGHtwo} \twobytwo{\RGbarAG^{31}}{\RGbarGbar^{31}}{\RGHAG^{21}}{\RGHGbar^{21}}
    \twobyone{\QABB}{\QBbarB}.\ethancheck
\end{split}
\end{equation}
Recalling $\PBBC = \onebytwo{\PBBCone}{\PBBCtwo}\twobyone{K_1}{K_2}$, we can also write $B$ as
\begin{equation}\label{eq:second_B}
\begin{split}
B &= \begin{bmatrix} \PBBbar  & \PBBC  \end{bmatrix}\twobytwo{\RBbarAB }{\RBbarBbar }{\RBCAB }{\RBCBbar }\begin{bmatrix}
    \QABB  \\ \QBbarB 
    \end{bmatrix} \\
    &= \begin{bmatrix}
    \PBBbar & \PBBCone & \PBBCtwo
    \end{bmatrix}\begin{bmatrix} \RBbarAB & \RBbarBbar \\ K_1 \RBCAB & K_1 \RBCBbar \\ K_2 \RBCAB & K_2 \RBCBbar \end{bmatrix} \twobyone{\QABB}{\QBbarB}.  \ethancheck
\end{split}
\end{equation}
Thus, combining \eqref{eq:first_B} and \eqref{eq:second_B}, we have
\begin{gather*}
    \begin{bmatrix}
    \PBBbar & \PBBCone & \PBBCtwo
    \end{bmatrix}\left(\begin{bmatrix} \RBbarAB & \RBbarBbar \\ K_1 \RBCAB & K_1 \RBCBbar \\ K_2 \RBCAB & K_2 \RBCBbar \end{bmatrix} - \begin{bmatrix}
    \RGbarAG^{21} & \RGbarGbar^{21} \\
    \RGHAG^{11} & \RGHGbar^{11} \\
    \RGbarAG^{11} & \RGbarGbar^{11}
    \end{bmatrix}\right)\twobyone{\QABB}{\QBbarB} \\
    = \onebytwo{\PGGbarthree}{\PGGHtwo} \twobytwo{\RGbarAG^{31}}{\RGbarGbar^{31}}{\RGHAG^{21}}{\RGHGbar^{21}}
    \twobyone{\QABB}{\QBbarB}.\ethancheck
\end{gather*}
Since the column spaces of the left- and right-hand sides have trivial intersection, both the left- and right-hand side must be zero. Since $\begin{bmatrix}
    \PBBbar & \PBBCone & \PBBCtwo
    \end{bmatrix}$ and $\onebytwo{\PGGbarthree}{\PGGHtwo}$ have full column rank and $\twobyone{\QABB}{\QBbarB}$ has full row rank, we conclude  
\begin{align*}
    \twobytwo{\RGbarAG^{31}}{\RGbarGbar^{31}}{\RGHAG^{21}}{\RGHGbar^{21}} &= \twobytwo{0}{0}{0}{0}, \\\begin{bmatrix}
    \RGbarAG^{21} & \RGbarGbar^{21} \\
    \RGHAG^{11} & \RGHGbar^{11} \\
    \RGbarAG^{11} & \RGbarGbar^{11} 
    \end{bmatrix} &= 
    \begin{bmatrix} \RBbarAB & \RBbarBbar \\ K_1 \RBCAB & K_1 \RBCBbar \\ K_2 \RBCAB & K_2 \RBCBbar \end{bmatrix}.\ethancheck
\end{align*}
It follows that 
\begin{equation}\label{eq:RGHAG}
    \RGHAG = \twobytwo{\RGHAG^{11}}{\RGHAG^{12}}{\RGHAG^{21}}{\RGHAG^{22}} = \twobytwo{K_1\RBCAB}{\RGHAG^{12}}{0}{\RGHAG^{22}}.
\end{equation}
\proofmarker{Comparison of solutions to \eqref{eq:lrcp} and \eqref{eq:lrcpcols}.} Substituting expressions for $\PHHbar$ and $\PHGH$ (cf.~\eqref{eq:Ccolumbasis}), $\QAGA$ (cf.~\eqref{eq:Arowbasis}), and $\RGHAG$ (cf.~\eqref{eq:RGHAG}) in \eqref{eq:additional_column_solutions}, we have
\begin{align*}
    X^{A,\BG,\CH}_{\FAbarbullet,\FHbar} &= \FAbarbullet\QAbarbulletA  + \onebytwo{\PCBCone }{\PHGHtwo  }\twobytwo{K_1\RBCAB }{\RGHAG^{12}}{0}{\RGHAG^{22}}\twobyone{\QABbulletA}{\QGcapAA}   \\
    &\qquad+ \begin{bmatrix} \PCBCtwo  & \PCCbar  & \PHHbarthree   \end{bmatrix}\FHbar  \\
    &= \FAbarbullet\QAbarbulletA  + \PCBCone K_1\RBCAB \QABbulletA  \\
    &\qquad+ (\PCBCone \RGHAG^{12} + \PHGHtwo  \RGHAG^{22}  )\QGcapAA + \begin{bmatrix} \PCBCtwo  & \PCCbar  & \PHHbarthree   \end{bmatrix}\FHbar \\
    & = \begin{bmatrix}\FAbarbullet & \PCBCone \RGHAG^{12} + \PHGHtwo  \RGHAG^{22} \end{bmatrix} \begin{bmatrix} \QAbarbulletA \\ \QGcapAA \end{bmatrix}   + \PCBCone K_1\RBCAB \QABA  \\
    &\qquad + \PCBCone K_1\RBCAB J\QAbarA  + \begin{bmatrix} \PCBCtwo  & \PCCbar  & \PHHbarthree   \end{bmatrix}  \FHbar,
\end{align*}
where in the last equality we invoke \eqref{eq:QABbulletA}. Since the variable $\FHbar$ is free, we may reparametrize to a new free variable $\FHbarnum{\bullet}$ related to $\FHbar$ by $\FHbar = \begin{bmatrix} K_2\RBCAB Q^T_{AB,A} \\ 0 \\ 0 \end{bmatrix}+\FHbarnum{\bullet}$ , we have solutions of the form
\begin{align*}
    X^{A,\BG,\CH}_{\FAbarbullet,\FHbar} &=  \begin{bmatrix}\FAbarbullet & \PCBCone \RGHAG^{12} + \PHGHtwo  \RGHAG^{22} \end{bmatrix} \begin{bmatrix} \QAbarbulletA \\ \QGcapAA \end{bmatrix}  +\PCBCone K_1\RBCAB  J\QAbarA  \\
    &\qquad+ \begin{bmatrix} \PCBCone &  \PCBCtwo \end{bmatrix} \begin{bmatrix}
     K_1 \RBCAB \\   K_2 \RBCAB
    \end{bmatrix} \QABA + \begin{bmatrix} \PCBCtwo  & \PCCbar  & \PHHbarthree   \end{bmatrix}  \FHbarnum{\bullet} \\
    &= \begin{bmatrix}\FAbarbullet & \PCBCone \RGHAG^{12} + \PHGHtwo  \RGHAG^{22} \end{bmatrix} \begin{bmatrix} \QAbarbulletA \\ \QGcapAA \end{bmatrix} +\PCBCone K_1\RBCAB  J\QAbarA  \\
    &\qquad+ \PCBC \RBCAB  \QABA + \begin{bmatrix} \PCBCtwo  & \PCCbar  & \PHHbarthree   \end{bmatrix}  \FHbarnum{\bullet}.
\end{align*}
We note that since $\twobyone{\QAbarbulletA}{\QGcapAA} = L_{11}^{-1}\QAbarA$ (cf.~\eqref{eq:Arowbasis}),
\begin{align*}
    X^{A,\BG,\CH}_{\FAbarbullet,\FHbar} &=  \begin{bmatrix}\FAbarbullet & \PCBCone \RGHAG^{12} + \PHGHtwo  \RGHAG^{22} \end{bmatrix} L_{11}^{-1} \QAbarA +\PCBCone K_1\RBCAB  J\QAbarA  \\
    &\qquad+ \PCBC \RBCAB  \QABA + \begin{bmatrix} \PCBCtwo  & \PCCbar  & \PHHbarthree   \end{bmatrix}  \FHbarnum{\bullet} \\
    &=\left( \begin{bmatrix}\FAbarbullet & \PCBCone \RGHAG^{12} + \PHGHtwo  \RGHAG^{22} \end{bmatrix} L_{11}^{-1}   + \PCBCone K_1\RBCAB  J \right) \QAbarA  \\
    &\qquad+ \PCBC \RBCAB  \QABA + \begin{bmatrix} \PCBCtwo  & \PCCbar  & \PHHbarthree   \end{bmatrix}  \FHbarnum{\bullet}.
\end{align*}
The solutions to \eqref{eq:lrcp} are of the form
\begin{align*}
    X^{A,B,C}_{\FAbar,\FCbar} &= \FAbar\QAbarA + \PCBC \RBCAB \QABA + \PCCbar \FCbar.\ethancheck
\end{align*}
%
%

\proofmarker{Determination of $Z$ and conclusion.} 
The restricted class of solutions of interest is $\mathscr{S}_{A,B,C}^Y := \{ X^{A,B,C}_{\FAbar,Y} : \FAbar \mbox{ free}\}$. We set 
\begin{equation*}
    \FHbarnum{\bullet} = \begin{bmatrix}
    0 \\ Y \\ 0
    \end{bmatrix} \mbox{ or equivalently } Z := \FHbar = \begin{bmatrix}
    -K_2\RBCAB\QABA  \\ Y \\ 0
    \end{bmatrix}.
\end{equation*}
Then,
\begin{equation*}
    \mathscr{S}_{A,\BG,\CH}^Z := \{ X^{A,\BG,\CH}_{F,Z} : F \mbox{ free} \} \subseteq \mathscr{S}_{A,B,C}^Y,\mathscr{S}_{A,\BG,\CH}, \ethancheck
\end{equation*}
which was as to be shown.

\end{proof}

We now consider a removal of rows from the top of the matrix. Define $\EA := \twobyone{E}{A}$ and $\FB := \twobyone{F}{B}$.\ethancheck

\begin{lemma}[Removal of Rows]\label{lem:removal_of_rows}
    Consider a subset 
    \begin{equation*}
        \mathscr{S}_{\EA,\FB,C}^Y := \{ \FEbar \QEbarE + \PCFC \RFCEF \QEFE + \PCCbarbullet Y : \FEbar \textnormal{ free} \}
    \end{equation*}
    of the solution set of the low-rank completion problem \eqref{eq:lrcprows} as computed by Algorithm 2 where the free variable ``$\FCbarbullet$'' is set to a fixed value $Y$. Suppose $R(F^T) \cap R(B^T) = \{0\}$. Then there exist choices for the bases $\PAbar$, $\PAB$, $\PBbar$, $\QBbar$, $\QBC$, and $\QCbar$ in Algorithm~\ref{alg:2x2lrcp} and a matrix $Z$ depending on $Y$ such that
    \begin{equation*}
        \mathscr{S}_{A,B,C}^Z := \{ \QAbarA \FAbar + \PCBC \RBCAB  \QABA + \PCCbar Z : \FAbar \textnormal{ free} \} \subseteq \mathscr{S}_{\EA,\FB,C}^Y.
    \end{equation*}
	
\end{lemma}
\begin{proof}
The proof of this lemma follows a similar procedure to Lemma \ref{lem:addition_of_cols} and is therefore omitted.
\end{proof}

\subsubsection{Modification of the Hankel blocks}\label{sec:hankel_block_modification}

Lemma~\ref{lem:addition_of_cols} and \ref{lem:removal_of_rows} require the assumptions that $R(B)\cap R(G)= \{0\}$ and $R(F^T) \cap R(B^T) = \{0\}$. These assumptions may sound restrictive,  but it turns out that we can rewrite our original overlapping Hankel form into an equivalent problem where these properties hold for every pair of overlapping Hankel blocks. This is the purpose of Lemma~\ref{lem:equivhankelproblem}. Before we get that result, consider the following proposition.

\begin{proposition}\label{prop:reparametrization_of_cols}
If $R\left(\twobyone{B}{C}\right) = R\left(\twobyone{B'}{C'}\right)$, then the solution sets \eqref{eq:lrcp_solution_set} for the low rank completion problems \eqref{eq:lrcp} for the matrices
\begin{displaymath}
\begin{bmatrix} A & B \\ X & C \end{bmatrix}, \quad  \begin{bmatrix} A & B' \\ X & C' \end{bmatrix}
\end{displaymath}
are identical. A similar statement can be formulated also whenever $R\left(\onebytwo{A}{B}^T\right) = R\left(\onebytwo{A'}{B'}^T\right)$.
\end{proposition}

\begin{proof}
    For any choice of $X\in\mathbb{R}^{m_2\times m_1}$, $M(A,B,C;X)$ and $M(A,B',C';X)$, namely 
    \begin{equation*}
        \dim \left[ R\left(\twobyone{A}{X}\right) + R\left(\twobyone{B}{C}\right) \right] = \dim \left[ R\left(\twobyone{A}{X}\right) + R\left(\twobyone{B'}{C'}\right) \right].\ethancheck \nithincheck
    \end{equation*}
\end{proof}
\begin{proposition}\label{prop:feasibility_of_reparametrization}
There exists a  $G$ and $H$ such that 
\begin{equation*}
    R\left(\twobytwo{B}{G_{\rm old}}{C}{H_{\rm old}}\right) = R\left(\twobytwo{B}{G}{C}{H}\right), \quad R(B) \cap R(G) = \{0\},
\end{equation*}
Likewise, there exists a $E$ and $F$ such that 
\begin{equation*}
    R\left(\begin{bmatrix}
	E_{\rm old} & F_{\rm old} \\
	A & B 
	\end{bmatrix}^T \right)  =  R\left(\begin{bmatrix}
	E & F \\
	A & B 
	\end{bmatrix}^T \right),\quad  R(B^T) \cap R(F^T) = \{0\}.
\end{equation*}
\end{proposition}
\begin{proof}
Choose a basis $P_1$ for $R(B) \cap R(G_{\rm old})$ and extend to a basis $\onebytwo{P_1}{P_2}$ for $G_{\rm old}$ and $\onebytwo{P_1}{P_3}$ for $B$. Then there exists a nonsingular matrices $S$ and $T$ such that $G_{\rm old}S = \onebytwo{P_1}{P_2}$ and $BT = \onebytwo{P_1}{P_3}$. Then, 
\begin{equation*}
    B\underbrace{T \twobytwo{I}{0}{0}{0}}_{:= T'} \onebytwo{P_1}{0},
\end{equation*}
where the sizes of the zero matrices is chosen so that $BT'$ and $G_{\rm old}S$ have the same size. Then
\begin{equation*}
    \twobytwo{B}{G_{\rm old}}{C}{H_{\rm old}}\underbrace{\twobytwo{I}{-T'}{0}{S}}_{=\Phi} = \twobytwo{B}{\onebytwo{0}{P_2}}{C}{H_{\rm old} S - CT'}.
\end{equation*}
Set $G := \onebytwo{0}{P_2}$ and $H := H_{\rm old} S - CT'$. Since $\Phi$ is nonsingular,
\begin{equation*}
    R\left(\twobytwo{B}{G_{\rm old}}{C}{H_{\rm old}}\right) = R\left(\twobytwo{B}{G}{C}{H}\right).
\end{equation*}
It is clear that $R(B) \cap R(G) = R(B) \cap R(P_2) = \{0\}$ by construction of $P_2$.
\end{proof}

Recall the $k$th Hankel block $\mathcal{H}_k(X)$ as defined by \eqref{eq:hankelblock}. For $k\le n-2$, we can partition the $k$th and $(k+1)$st Hankel blocks by
    \begin{equation}
    \begin{split}
    \mathcal{H}_{k}(X) &= \left[\begin{array}{c|ccc}
    A_{(k+1)1} & A_{(k+1) 2} & \cdots & A_{(k+1)k} \\ \hline
    A_{(k+2)1} & A_{(k+2) 2} & \cdots & A_{(k+2)k} \\
    \vdots & \vdots & \ddots    & \vdots \\
    A_{(n-1)1} & A_{(n-1)2} & \cdots & A_{(n-1)k} \\ \hline
    X & A_{n2} & \cdots & A_{nk} 
    \end{array}\right] = \begin{bmatrix}
    \hat{E}^k & \hat{F}^k \\
    \hat{A}^k & \hat{B}^k \\
    X & \hat{C}^k
    \end{bmatrix}, \\
    \mathcal{H}_{k+1}(X) &= \left[\begin{array}{c|ccc|c}
    A_{(k+2)1} & A_{(k+2) 2} & \cdots & A_{(k+2)k} & A_{(k+2)(k+1)} \\
    \vdots & \vdots & \ddots    & \vdots & \vdots \\
    A_{(n-1)1} & A_{(n-1)2} & \cdots & A_{(n-1)k} & A_{(n-1)(k+1)}\\ \hline
    X & A_{n2} & \cdots & A_{nk}  & A_{n(k+1)}
    \end{array}\right] = \begin{bmatrix}
    \hat{A}^k & \hat{B}^k & \hat{G}^k \\
    X & \hat{C}^k & \hat{H}^k
    \end{bmatrix}.\ethancheck
    \end{split} \label{eq:hankelpartitioning}
    \end{equation}
Similar to Section \ref{sec:addcolremrows}, we denote
    \begin{equation*}
    \hat{\EA}^k = \twobyone{\hat{E}^k}{\hat{A}^k}, \hat{\FB}^k = \twobyone{\hat{F}^k}{\hat{B}^k}, \hat{\BG}^k = \onebytwo{\hat{B}^k}{\hat{G}^k}, \mbox{ and } \hat{\CH}^k = \onebytwo{\hat{C}^k}{\hat{H}^k}.
    \end{equation*}
Having said that, we can prove the following result.

\begin{lemma} \label{lem:equivhankelproblem}
Consider the Hankel blocks $\mathcal{H}_k(X)$ partitioned as in \eqref{eq:hankelpartitioning}. There exist modified Hankel blocks $\mathcal{K}_k(x)$, $1\le k \le n-1$, which can be partitioned as
\begin{equation}\label{eq:modified_hankel_defs}
\mathcal{K}_k(X) = \begin{bmatrix}
        E^k & F^k \\
        A^k & B^k \\
        X & C^k
        \end{bmatrix}, \quad \mathcal{K}_{k+1}(X) = \begin{bmatrix} A^k & B^k & G^k \\
        X & C^k & H^k\end{bmatrix}
\end{equation}
such that
\begin{enumerate}[label={(\roman*)}]
        \item For every choice of $X$ and $1\le k \le n-1$, $\rank \mathcal{H}_k(X) = \rank \mathcal{K}_k(X)$. In particular, $\hat{X}$ minimizes the rank of all Hankel blocks $\mathcal{H}_k(X)$ if, and only if, it minimizes the rank of all Hankel blocks $\mathcal{K}_k(X)$.\label{item:hankel_prop_1}
        \item For every $1\le k \le n-2$, \label{item:hankel_prop_2}
        \begin{enumerate}[label={(\alph*)}]
            \item $R((B^k)^T) \cap R((F^k)^T) = \{0\}$, and \label{item:hankel_prop_2a}
            \item $R(B^k) \cap R(G^k) = \{0\}$. \label{item:hankel_prop_2b}
        \end{enumerate}
    \end{enumerate}
\end{lemma}

\begin{proof} 
    By ``sweeping'' over the rows, we shall first construct intermediate Hankel blocks
    \begin{equation}
    \mathcal{H}_{k}'(X) = \begin{bmatrix}
    \tilde{E}^k & \tilde{F}^k \\
    \tilde{A}^k & \tilde{B}^k \\
    X & \tilde{C}^k
    \end{bmatrix}, \quad \mathcal{H}_{k+1}'(X) = \begin{bmatrix}
    \tilde{A}^k & \tilde{B}^k & \tilde{G}^k \\
    X & \tilde{C}^k & \tilde{H}^k
    \end{bmatrix} \label{eq:hankelprime}
    \end{equation}
    which satisfy \ref{item:hankel_prop_1} and \ref{item:hankel_prop_2}-\ref{item:hankel_prop_2a}, but not \ref{item:hankel_prop_2}-\ref{item:hankel_prop_2b}. The construction procedure is as follows:
    \begin{enumerate}
        \item For each $1\le k \le n-2$,  define 
        \begin{displaymath}
        \onebytwo{\tilde{E}^k}{\tilde{F}^k} := \onebytwo{\hat{E}^k}{\hat{F}^k} + S_1^k\onebytwo{\hat{A}^k}{\hat{B}^k}, 
        \end{displaymath}
        where $S_1^k$ is chosen to ensure $R((\tilde{F}^k)^T) \cap R((\tilde{B}^k)^T) = \{0\}$ (see Proposition~\ref{prop:feasibility_of_reparametrization}).
        \item Set $\tilde{A}^{n-1} := \hat{A}^{n-1}$, $\tilde{B}^{n-1} := \hat{B}^{n-1}$ and $\tilde{A}^k := \twobyone{\tilde{E}^{k+1}}{\tilde{A}^{k+1}}$, $\tilde{B}^k := \twobyone{\tilde{F}^{k+1}}{\tilde{B}^{k+1}}$, $\tilde{C}^k := \hat{C}^k$.
    \end{enumerate}
By construction, \eqref{eq:hankelprime} satisfies \ref{item:hankel_prop_2}-\ref{item:hankel_prop_2a}, and since  
 \begin{equation*}
     R\left(\begin{bmatrix}
    \tilde{E}^k & \tilde{F}^k \\
    \tilde{A}^k & \tilde{B}^k \\
    \end{bmatrix}^T\right) = R\left(\begin{bmatrix}
    \hat{E}^k & \hat{F}^k \\
    \hat{A}^k & \hat{B}^k \\
    \end{bmatrix}^T\right),
 \end{equation*}
    by Proposition \ref{prop:reparametrization_of_cols} it also satisfies \ref{item:hankel_prop_1}. To obtain the desired result,   in much the same way, we now perform sweeps on the columns of $\mathcal{H}_{k}'(X)$ to obtain $\mathcal{K}_k(X)$, i.e. 
    \begin{enumerate}
        \item For each  $1\le k \le n-2$, we define
        \begin{displaymath}\twobyone{G^k}{H^k} := \twobyone{\tilde{G}^k}{\tilde{H}^k} + \twobyone{\tilde{B}^k}{\tilde{C}^k}(S_2^k)^T
        \end{displaymath}
        such that $R(B^k) \cap R(G^k) = \{0\}$ (again using Proposition~\ref{prop:feasibility_of_reparametrization}).
        \item Set ${B}^{1} := \tilde{B}^{1}$, ${C}^{1} := \tilde{C}^{1}$ and $B^{k} := \onebytwo{B^{k-1}}{G^{k-1}}$, $C^k := \onebytwo{C^{k-1}}{H^{k-1}}$, ${A}^k := \tilde{A}^k$.
    \end{enumerate}
    The resulting Hankel blocks $\mathcal{K}_k(X)$ satisfy \ref{item:hankel_prop_1} and \ref{item:hankel_prop_2}, which was to be constructed.
\end{proof}

\subsection{Proof of the main result}\label{sec:proof_of_main_result}
We are now ready to prove the main result.

\begin{proof}[Proof of Theorem~\ref{thm:rankcompletion}]

   We shall inductively show that the sets $\mathcal{S}_k$, as defined in \eqref{eq:solutionsets}, i.e.
   \begin{displaymath}
    \mathcal{S}_k := \left\{ \hat{X}\in\mathbb{R}^{m_2\times m_1}:\quad  \rank  \mathcal{H}_k(\hat{X}) = \min_{ X \in \mathbb{R}^{N_n \times N_1} } \rank \mathcal{H}_\ell(X),\qquad \ell=1,2\ldots,k  \right\} 
   \end{displaymath}
   remain nonempty for $k=1,2,\ldots,n-1$.

   First of all, we observe that by Lemma~\ref{lem:equivhankelproblem}, the rank minimization of $\mathcal{H}_k(X)$ can be replaced with that of $\mathcal{K}_k(X)$, without affecting the final solution set. Hence, 
   \begin{displaymath}
    \mathcal{S}_k := \left\{ \hat{X}\in\mathbb{R}^{m_2\times m_1}:\quad  \rank  \mathcal{K}_k(\hat{X}) = \min_{ X \in \mathbb{R}^{N_n \times N_1} } \rank \mathcal{K}_\ell(X),\qquad \ell=1,2\ldots,k  \right\}. 
   \end{displaymath}
   We now proceed by induction. As a base case, we observe that
    \begin{equation*}
         \mathscr{S}_{\EA^1,\FB^1,C^1}^{Y_1} = \left\{ F \begin{bmatrix} A_{21} \\ A_{31} \\ \vdots \\ A_{(n-1)1} \end{bmatrix} : F \mbox{ free} \right\}  = \mathcal{S}_1 
    \end{equation*}
    where $Y_1$ is an empty matrix of the appropriate size.
    
    Inductively, suppose that every element of the nonempty set $ \mathscr{S}_{\EA^k,\FB^k,C^k}^{Y_k} \subseteq \mathcal{S}_k$ minimizes the first $k$ Hankel blocks for some matrix $Y_k$. Then by Lemma \ref{lem:removal_of_rows}, there exists $W$ such that the nonempty set $\mathscr{S}_{A^k,B^k,C^k}^W$ satisfies $\mathscr{S}_{A^k,B^k,C^k}^W \subseteq \mathscr{S}_{\EA^k,\FB^k,C^k}^{Y_k} \subseteq \mathcal{S}_k$. Then, by Lemma \ref{lem:addition_of_cols}, there exists $Y_{k+1}$ such that the nonempty set $\mathscr{S}_{A^k,\BG^k,\CH^k}^{Y_{k+1}}$ satisfies $\mathscr{S}_{A^k,\BG^k,\CH^k}^{Y_{k+1}} \subseteq \mathscr{S}_{A^k,B^k,C^k}^W \subseteq \mathcal{S}_{k}$. But $\mathscr{S}_{A^k,\BG^k,\CH^k}^{Y_{k+1}} \subseteq\mathscr{S}_{A^k,\BG^k,\CH^k}$ and $\EA^{k+1} = A^k$, $\FB^{k+1} = \BG^{k}$, and $C^{k+1} = \CH^k$ (cf. \eqref{eq:modified_hankel_defs}) as well so 
    \begin{displaymath}\mathscr{S}_{\EA^{k+1}, \FB^{k+1}, C^{k+1}}^{Y_{k+1}} = \mathscr{S}_{A^k,\BG^k,\CH^k}^{Y_{k+1}} \subseteq \mathscr{S}_{\EA^{k+1}, \FB^{k+1}, C^{k+1}} \cap \mathcal{S}_{k} = \mathcal{S}_{k+1}.
    \end{displaymath}
    Therefore, we conclude there exists an element $\hat{X} \in \mathscr{S}_{\EA^{n-1},\FB^{n-1},C^{n-1}}^{Y_{n-1}} \subseteq \mathcal{S}_{n-1}$ such that
    \begin{equation*}
        \rank \mathcal{H}_\ell(\hat{X}) = \rank \mathcal{K}_\ell(\hat{X}) \le \rank \mathcal{K}_{\ell}(X) = \rank \mathcal{H}_\ell(X)
    \end{equation*}
    for every matrix $X$ and $1\le \ell \le n-1$. \nithincheck
    
    
    
    
    
    
    
\end{proof}




\begin{remark}\label{rmk:complexity}
As presented above, the solution of the overlapping Hankel block low-rank completion problem has time complexity $\mathcal{O}(N^4)$, as it requires performing dense linear algebraic calculations $\mathcal{O}(N^3)$ on $\mathcal{O}(N)$ different Hankel blocks. However, there is reason to be optimistic. If at any point during the computation, we find that there is a \textit{unique} solution minimizing the first $k$ Hankel blocks, then Theorem \ref{thm:rankcompletion} ensures that this solution in fact minimizes all the Hankel blocks and we can terminate our calculation early.
If a unique solution can be found by only considering $\mathcal{O}(N)$ Hankel blocks, then the complexity of solving the low rank completion problem reduces to $\mathcal{O}(N^2)$, as the first Hankel blocks are skinny with sizes $\mathcal{O}(N) \times \mathcal{O}(1)$. However, to achieve this improved time complexity, the two stages of the proof above (precomputation and addition of columns/removal of rows) need to be intermingled, which we believe will be possible. We conjecture that with a more careful algorithm, the complexity of solving this problem can be made to be $\mathcal{O}(N^2)$ in the worst case.
\end{remark}

\begin{conjecture}
There exists an algorithm which computes a candidate solution $\hat{X}$ to the overlapping Hankel block low rank completion problem in $\mathcal{O}(N^2)$ time.
\end{conjecture}

\begin{example}
Consider the CSS construction for the matrix $A = 100 I + xy^T + e_5e_2^T + e_8e_5^T$ for 
\begin{align*}
    x &= \begin{bmatrix} 2 & 1 & 2 & 1 & 1 & 2 & 1 & 1 & 1 & 2 \end{bmatrix}^T,\\
    y &= \begin{bmatrix} 1 & 3 & 1 & 2 & 1 & 3 & 2 & 1 & 1 & 1\end{bmatrix}^T,
\end{align*}
and $(e_i)_j = \delta_{ij}$. Choosing $N_1 = \cdots = N_5 = 2$, \eqref{eq:ohblrcp} reduces to minimizing the Hankel block ranks in
\begin{equation*}
    \left[\begin{array}{cc|cc|cc|cc}
    2 & 6 &&&&&& \\
    1 & 3 &&&&&& \\\hline 
    1 & 4 & 1 & 2 &&&&\\
    2 & 6 & 2 & 4 &&&&\\\hline
    1 & 3 & 1 & 2 & 1 & 3 &&\\
    1 & 3 & 1 & 2 & 2 & 3 &&\\\hline
    ? & ? & 1 & 2 & 1 & 3 & 2 & 1\\
    ? & ? & 2 & 4 & 2 & 6 & 4 & 2
    \end{array}\right]
\end{equation*}
After applying Lemma~\ref{lem:equivhankelproblem}, we reduce to an equivalent problem
\begin{equation*}
    \left[\begin{array}{cc|cc|cc|cc}
    0 & 0 &&&&&& \\
    0 & 0 &&&&&& \\\hline 
    0 & 1 & 0 & 0 &&&&\\
    0 & 0 & 0 & 0 &&&&\\\hline
    1 & 3 & 1 & 2 & 0 & 0 &&\\
    1 & 3 & 1 & 2 & 1 & 0 &&\\\hline
    ? & ? & 1 & 2 & 0 & 0 & 0 & 0\\
    ? & ? & 2 & 4 & 0 & 0 & 0 & 0
    \end{array}\right]
\end{equation*}
It is clear that it is sufficient to consider the row rank completion problems for the second and third Hankel blocks. For the second Hankel block, we have
\begin{equation*}
    \PAbar = \begin{bmatrix} 1 \\ 0 \\ 3 \\ 3 \end{bmatrix}, \PAB = \begin{bmatrix} 0 \\ 0 \\ 1 \\ 1 \end{bmatrix}, \QBC = \begin{bmatrix} 1 & 2\end{bmatrix},
\end{equation*}
with $\PBbar$, $\QBbar$, and $\QCbar$ being empty matrices of the appropriate sizes. Then
\begin{equation*}
    A = \onebytwo{\PAbar}{\PAB}\twobyone{\QAbarA}{\QABA} = \left[\begin{array}{c|c} 1 & 0 \\
    0 & 0 \\
    3 & 1 \\
    3 & 1\end{array}\right]\left[\begin{array}{cccc} 0 & 1 \\\hline 1 & 0 \end{array} \right]
\end{equation*}
and
\begin{equation*}
    B = \begin{bmatrix} \PBBbar & \PBBC \end{bmatrix} \begin{bmatrix} \RBbarAB & \RBbarBbar \\ \RBCAB & \RBCBbar \end{bmatrix} \begin{bmatrix} \QABB \\ \QBbarB \end{bmatrix} = \left[ \begin{array}{c|c} \: & 0 \\ \: & 0 \\ \: & 1 \\ \: & 1\end{array} \right] \left[\begin{array}{c|c} \: & \: \\\hline
    1 \: \end{array} \right] \left[\begin{array}{cc} 1 & 2 \\\hline \: & \: \end{array} \right].
\end{equation*}
We easily conclude $\PCBC = \onebytwo{1}{2}^T$. Thus,
\begin{align*}
    \mathcal{S}_2 &= \left\{ \FAbar \QAbarA + \PCBC \RBCAB \QABA + \PCCbar \FCbar : \FAbar, \FCbar \mbox{ free} \right\} \\
                  &= \left\{ \twobyone{f_1}{f_2} \onebytwo{0}{1} + \twobytwo{1}{0}{2}{0} : f_1,f_2 \in \real \right\}.
\end{align*}
We then perform the removal of rows calculation to
\begin{equation*}
    \left[\begin{array}{cc|cc}
    1 & 3 & 1 & 2 \\
    1 & 3 & 1 & 2 \\\hline
    ? & ? & 1 & 2 \\
    ? & ? & 2 & 4 
    \end{array}\right],
\end{equation*}
where we see that 
\begin{equation*}
    \hat{X} = \twobytwo{1}{3}{2}{6}
\end{equation*}
is the unique solution. We conclude that the optimal Hankel block ranks are $1$, $2$, $2$, and $1$.
\end{example}

\begin{remark}\label{rmk:oneHankel}
Following Lemma~\ref{thm:2x2_low_rank_completion}, the condition for uniqueness for the low-rank completion problem \eqref{eq:lrcp} is $R(A) \subseteq R(B)$ and $R(C^T) \subseteq R(B^T)$. By Theorem~\ref{thm:rankcompletion}, it is guaranteed that if there is a \textit{unique} rank minimizer $\hat{X}$ for \textit{any} Hankel block $\mathcal{H}_k(X)$, then $\hat{X}$ minimizes the rank of \textit{all} the Hankel blocks. Even if the solution set to the low-rank completion problem \eqref{eq:lrcp} is of the form $\mathscr{S}_{A,B,C} = \{ \FAbar \QAbarA + \PCBC \RBCAB \QABA + \PCCbar \FCbar : \FAbar, \FCbar \mbox{ free} \}$, then $X - \hat{X} \in \{ \FAbar \QAbarA + \PCCbar \FCbar : \FAbar, \FCbar \mbox{ free} \}$ for any global rank minimizer $\hat{X}$. Thus, if $\QAbarA$ has $r$ rows and $\PCCbar$ has $s$ columns, then any $X\in \mathscr{S}_{A,B,C}$ is an approximate solution to \eqref{eq:ohblrcp} in the sense that

\begin{equation*}
    \mathcal{H}_k(\hat{X}) \le \mathcal{H}_k(X) \le \mathcal{H}_k(\hat{X}) + r + s, \quad k=1,2,\ldots,n-1.
\end{equation*}

Thus, for many problems, an exact or good approximate solution to \eqref{eq:ohblrcp} can be obtained by considering just one Hankel block, without the heavy lifting (and significant computational cost) of Lemmas~\ref{lem:addition_of_cols}, \ref{lem:removal_of_rows}, and \ref{lem:equivhankelproblem}. \nithin{Given a typical problem of our interest, what is the probability of having some Hankel block with a unique completion? Does it not indirectly tell you that if this occurs, the matrix does not allow for a lot compression in the CSS representation with respect to SSS?}
\end{remark}

 \notrunc{More detailed proofs of Lemmas \ref{lem:addition_of_cols} and \ref{lem:removal_of_rows} conducive to direct implementation in a computer programming language can be found in Appendix \ref{app:numerical_details}.\ethancheck}


\section{Numerical examples}\label{sec:numerical_results}

Our goal in investigating CSS representations is not to propose them as a replacement for SSS representations, but rather to use them as a stepping stone to studying $\graph$-SS representations for more complicated graphs. However, if computing the CSS representation is infeasible, then we have no hope for more complicated graph structures. In this section, we demonstrate that computing the CSS numerically \textit{is} tractable, and achieves the linear solve complexity even when the low-rank completion problem is solved inexactly (following Remark~\ref{rmk:oneHankel}).

The construction of a CSS representation comes at the computational cost of solving the overlapping Hankel block low rank completion problem ($\mathcal{O}(N^4)$ worst case, see Remark \ref{rmk:complexity}) on top of the standard procedure for finding an SSS representation ($\mathcal{O}(N^2)$ assuming $r_i^g,r_i^h = \mathcal{O}(1)$).  Despite the requirement of additional terms $U_0$ and $Q_0$, the remaining terms of a CSS representation could potentially be factors smaller than those of the SSS representation. However, as discussed in Remark~\ref{rmk:CSS_vs_SSS}, a matrix which has a significantly smaller CSS representation compared to an SSS representation can be permuted and repartitioned so that the permuted SSS representation has comparable size to the CSS representation in the natural ordering. Given that our ultimate goal is to consider more complicated graphs for which a complexity-reducing permutation is hard to find or does not exist, we shall only consider matrices in their given orderings for the duration of this section.

A direct execution of the techniques used to prove Theorem~\ref{thm:rankcompletion} to solve the problem \eqref{eq:ohblrcp} is prone to severe numerical issues, due to the inherent numerical instability of computing subspace intersections. Indeed, we observe such significant numerical issues in our own implementation. Since our method for solving \eqref{eq:ohblrcp} revolves so heavily around subspace computations, we believe the development of an efficient and stable algorithm to exactly solve \eqref{eq:ohblrcp} in all cases will be a challenging task, which we defer to future study. For these reasons, we shall only consider the $k$th Hankel block for solving \eqref{eq:ohblrcp} for $k = \lceil n/2 \rceil$, which provides an exact or good approximant of the solution in light of Remark~\ref{rmk:oneHankel}.

Our numerical tests are done in Julia, with times being reported as the average over five runs. The SSS construction algorithm (which is used as a subroutine by the CSS construction) is taken from \cite{chandrasekaran2002fast} and uses SVD's to compute low-rank approximations with an absolute singular value threshold of $10^{-8}$. The CSS construction was by Algorithm~\ref{alg:CSSconstruct} with Algorithm~\ref{alg:2x2lrcp} used on the $\lceil n/2\rceil$-th Hankel block. Once the representation has been constructed, $Ax = b$ is solved by using Julia's sparse direct solvers on the lifted sparse system \eqref{eq:EDV_sparse}.

As a favorable example, we shall consider matrices $A_{r,b}$ which are sums of a tridiagonal matrix $T$, a semi-separable matrix of rank $r$, and a full-rank corner block perturbation of size $b\times b$. The entries of the banded matrix are taken to be uniform $[0,1]$-valued random numbers, the corner block perturbations are taken to have random $[0,1]$-valued entries, and the semi-separable matrix is constructed by summing the lower triangular portion of a rank-$r$ matrix $R_1$ and the upper triangular portion of a rank-$r$ matrix $R_2$, where $R_1$ and $R_2$ are defined to be the product of a $N\times r$ and an $r\times N$ random $[0,1]$-valued matrix.
\begin{equation*}
    A_{k,r,b} = T + \mbox{lower}(R_1) + \mbox{upper}(R_2) + \CB(C_{b\times b}, D_{b\times b}).
\end{equation*}
For this study, we take $r = 10$ and $b = N^{1/2}$. For our partitions, we use $N_1 = \cdots = N_{N/4} = 4$ for SSS and $N_1 = N_k = \sqrt{N}$ and $N_2 = \cdots = N_{n-1}$, $k = (N - 2\sqrt{N})/4 + 2$, for CSS. This example is particularly favorable for CSS because it is of the form SSS-10 plus a corner block perturbation. Thus, if the low-rank completion problem is solved optimally, we will have $r_{n-1}^g ,r_1^h = \mathcal{O}(N^{1/2})$ but $r_i^g = r_i^h = r = 10$ for $1 < i < n-1$. The SSS representation, on the other hand, has $r_i^g , r_i^h= \mathcal{O}(N^{1/2})$ for all $i$, leading to a significant loss in performance.

\begin{figure}
    \centering
    \begin{subfigure}[b]{0.85\textwidth}
        \includegraphics[width =\textwidth]{\MyPath/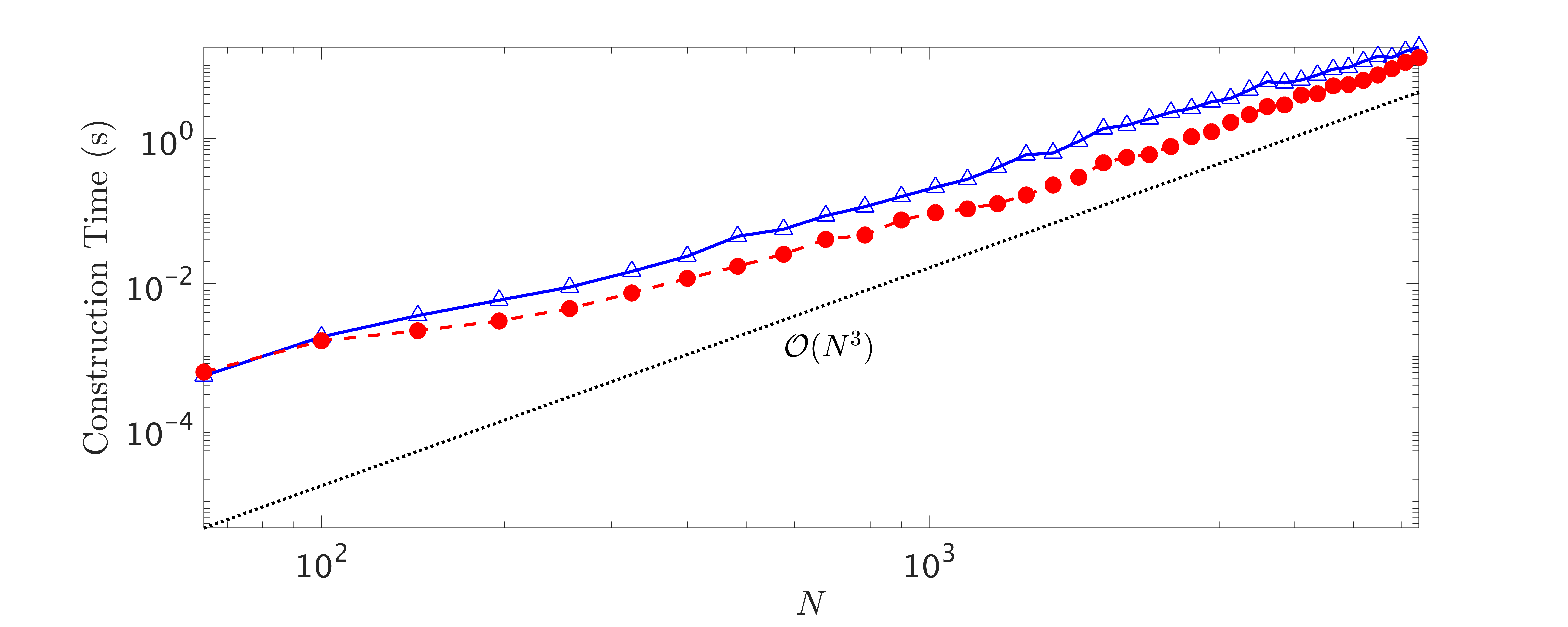} 
        \caption{} \label{subfig:pss_construction}
    \end{subfigure}
    
    \begin{subfigure}[b]{0.85\textwidth}
        \includegraphics[width = \textwidth]{\MyPath/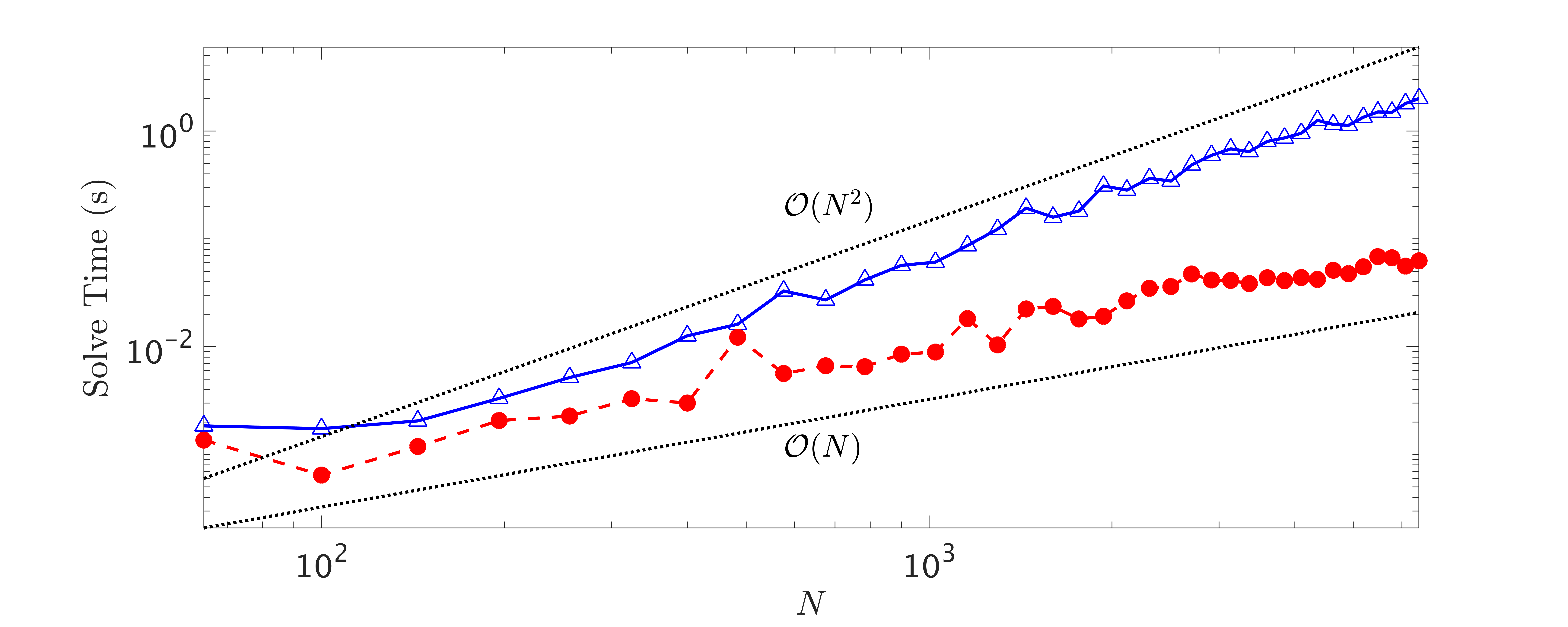}
        \caption{} \label{subfig:pss_solve}
    \end{subfigure}
    
    \begin{subfigure}[b]{0.85\textwidth}
        \includegraphics[width = \textwidth]{\MyPath/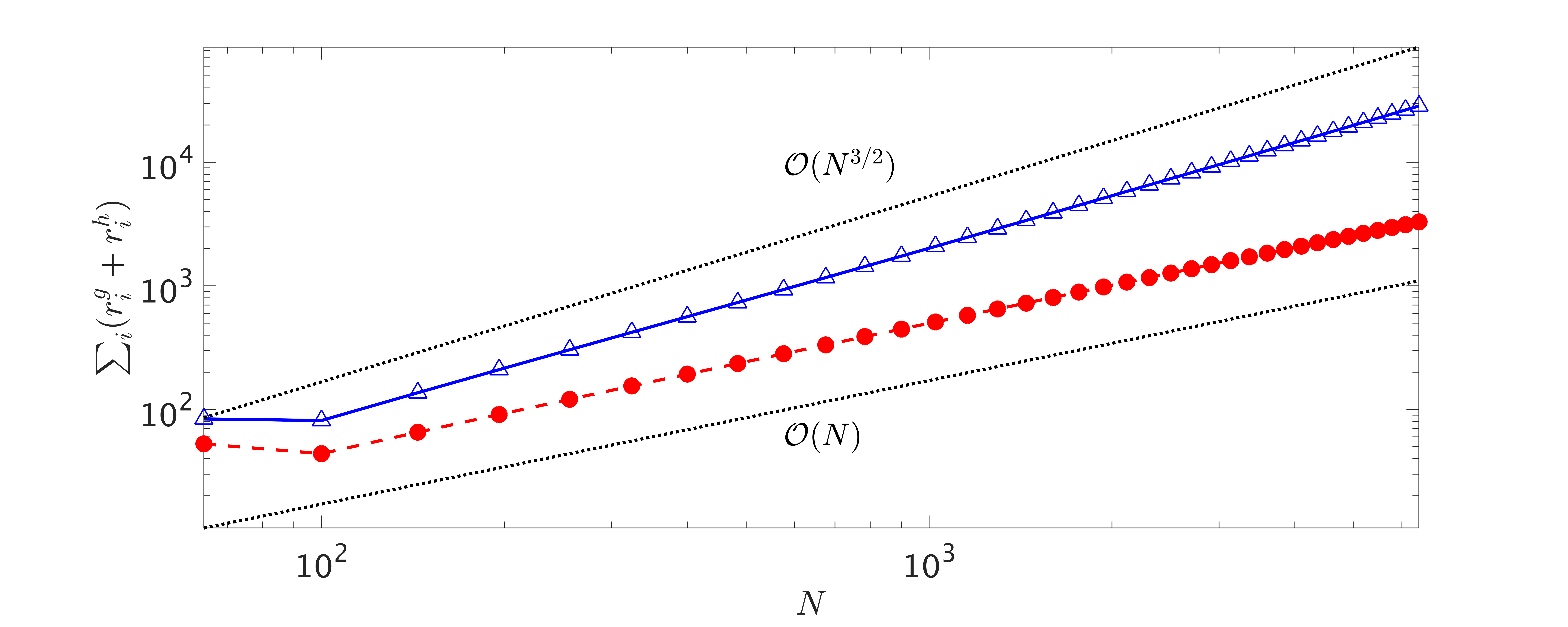}
        \caption{}
        \label{subfig:pss_size}
    \end{subfigure}
    \caption{Construction time, solve time, and representation size $\sum_i (r_i^g + r_i^h)$ for SSS (blue triangles) and CSS (red circles) representations for the class of matrices $A_{10,\sqrt{N}}$.}\label{fig:perturbed_semi_separable}
\end{figure}

The results are shown in Figure~\ref{fig:perturbed_semi_separable}. As expected, we observe a $\mathcal{O}(N^3)$ construction time for the CSS representation, as Algorithm~\ref{alg:2x2lrcp} involves dense matrix calculations on the $k$th Hankel block, which is of size $\mathcal{O}(N)\times \mathcal{O}(N)$. Despite the added overhead of solving the low-rank completion problem, the CSS construction is faster than the SSS construction since the ranks $r_i^g$ and $r_i^h$ are so much smaller than in the SSS case ($\mathcal{O}(1)$ vs $\mathcal{O}(N^{1/2})$ for $i > 1$). Possible further speedups for the CSS construction could be realized by using Algorithm~\ref{alg:2x2lrcp} on one of the first ``skinny'' $\mathcal{O}(N)\times \mathcal{O}(1)$ Hankel blocks, rather then the median $k = \lceil n/2 \rceil$-th Hankel block as we do here. (This benefit may be offset by the ranks $r_i^g$ and $r_i^h$ being higher than if the $\lceil n/2 \rceil$-th Hankel block is used.)

The solve time and size of the representation are also no surprise. For SSS, the size of the representation is $\sum_i (r_i^g + r_i^h) = \sum_i \mathcal{O}(N^{1/2}) = \mathcal{O}(N^{3/2})$ and the solve time is $\sum_i [(r_i^g)^2 + (r_i^h)^2] = \sum_i \mathcal{O}(N) = \mathcal{O}(N^2)$. For CSS, we have $r_{n-1}^g,r_1^h = \mathcal{O}(N^{1/2})$, but $r_i^g = r_i^h = r$ for $1 < i < n-1$ assuming the optimal solution to \eqref{eq:ohblrcp} is found. Thus, if the approximate solution to \eqref{eq:ohblrcp} is ``good enough'' (in the sense that $r_i^g, r_i^h = \mathcal{O}(1)$ for $i > 1$), then we would expect the size of the representation and the solve time to be linear time $\mathcal{O}(N)$. Indeed, as shown in Figures \ref{subfig:pss_solve} and \ref{subfig:pss_size}, this is exactly what we observe.

Clearly, the perturbed semi-separable matrices $A_{r,b}$ are specifically engineered to highlight the CSS representation's strengths over the SSS representations. We have also evaluated the performance of the CSS and SSS representations on a Cauchy kernel on the circle
\begin{equation*}
    B_{jk} = \begin{cases} \left|\exp\left(\frac{2\pi \sqrt{-1}j}{N}\right) - \exp\left(\frac{2\pi \sqrt{-1}k}{N}\right)\right|^{-1}, & j \ne k, \\
    0, & j = k.
    \end{cases}
\end{equation*}
We used the same block partitioning as with the previous example. For this example, the SSS and CSS representations had comparable solve times, with the CSS solve time consistently 10\%-20\% higher than the SSS representation---the smaller size of the CSS representation was offset by the additional fill-in of performing Gaussian elimination on a cycle graph rather than a line graph. The construction time for the CSS representation grew much faster than the construction time for SSS representation as $N$ was increased. We conclude that computing a near-optimal CSS representation is numerically tractable and has comparable solve time to the SSS representation (and much better solve time in specially constructed examples).

\section{Concluding remarks} \label{sec:overview}

We have proposed two representations, DV and $\graph$-SS representations, which both have fast solvers in time complexity given by sparse Gaussian elimination on the graph $\graph$. The $\graph$-SS representation has a linear time multiplication algorithm, so if a compact $\graph$-SS representation of $A^{-1}$ can be produced, $A^{-1}b$ could be computed in linear time.

We studied $\graph$-SS representation in the special case in which $\graph$ is a cycle graph, the so-called CSS representation. When we introduced $\graph$-SS representations, we ended our discussion with three lingering questions: \ref{item:GSS1} the existence and compactness of $\graph$-SS representations in terms of the GIRS property, \ref{item:GSS2} the tractability of construction, and \ref{item:GSS3} the algebraic properties of these representations. We now have the results to answer these questions for the CSS representation.

The existence of a ``best'' CSS representation was shown in Theorem~\ref{thm:CSS_construction}. It was shown in Proposition~\ref{prop:GIRSimpliesCSS} that the size of the entries in this representation for a GIRS-$c$ matrix are bounded by $2c$, and the factor of $2$ was shown to be tight in Example~\ref{ex:css_inversion}. These results give a complete answer to question \ref{item:GSS1}. In Sections~\ref{sec:css} and \ref{sec:hankelminimization}, a constructive proof was provided to compute this representation. While the resulting algorithm to compute the truly optimal representation had a $\mathcal{O}(N^4)$ worst case runtime, it was shown empirically in Section \ref{sec:numerical_results} that a relaxation (Remark~\ref{rmk:oneHankel}) of the problem could achieve a nearly optimal solution in $\mathcal{O}(N^3)$ time. We believe that faster exact or approximate algorithms for the CSS construction algorithm likely exist, which could be the subject of future investigation. We thus answer \ref{item:GSS2} in the affirmative, though faster algorithms may exist.

Question \ref{item:GSS3} is particularly important. Our eventual goal is to construct a $\graph$-SS representation of $A^{-1}$, so that the fast $\graph$-SS multiplication algorithm would give a linear-time solver for repeated right-hand side problems. One possible strategy to do this would be to first construct a $\graph$-SS representation for $A$ and them compute from this a $\graph$-SS representation of $A^{-1}$ without ever forming $A^{-1}$ explicitly. (This is particularly attractive for $A$ sparse.) To realize this strategy an important ingredient would be the size of the representation to not dramatically increase under inversion. In Section~\ref{sec:edv}, we make were unable to make statement of the size of the representation under inversion. Fortunately, in Theorem~\ref{thm:CSSInversion} we show that, under a reasonable set of assumptions, the size of the CSS representation can increase by a factor of no more than $6$ under inversion, a constant we believe is unlikely to be tight. (Though a factor two increase in the size of the minimal CSS representation is possible by Example~\ref{ex:css_inversion}.) This leaves open the feasibility of constructing a compact $\graph$-SS representation of $A^{-1}$ via a compact $\graph$-SS representation of $A$ without ever forming $A^{-1}$ explicitly.

Constructing DV and $\graph$-SS representations will be hard in general. While, DV and $\graph$-SS are guaranteed to exist for a large class of GIRS matrices (those possessing Hamiltonian paths), the existing construction provides no guarantee that these representations will be compact, with the sizes of the matrices in the representations bounded by the GIRS constant $c$. As an added complication, the efficiency of the $\graph$-SS representation may depend on the choice of Hamiltonian path. For the 2D mesh graph, there are many candidates worth trying, but for arbitrary graphs, finding even one Hamiltonian path is NP hard. Despite these challenges, we conjecture that constructing compact DV and $\graph$-SS is possible, encouraged by the tractability of SSS construction and CSS construction, and construction of DV representations for sparse matrices and their inverses. \ethancheck

\begin{conjecture}\label{conj:idv_compact}
If $(A,\graph)$ is GIRS-$c$, then $(A,\graph)$ possesses a DV representation with $r_i \le 2c$ for every $i \in \nodeset_\graph$.\ethancheck
\end{conjecture}

\begin{conjecture}\label{conj:edv_compact}
There exists a constant $1\le \gamma \ll N$ dependeing on $\graph$ such that if $(A,\graph)$ is GIRS-$c$, then $(A,\graph)$ possesses a $\graph$-SS representation with $r_i^g,r_i^h \le c\gamma$ for every $i \in \nodeset_\graph$. 
\end{conjecture}

Conjecture \ref{conj:edv_compact} is admittedly imprecise, but is purposely vague to handle the possibility that there may not exist a universal constant $\gamma$ for all graph structures. We hope that for a collection of graphs belonging to a ``family'', such a uniform constant $\gamma$ may be possible, which we conjecture to be the case for mesh graphs.

\begin{conjecture}
There exists a constant $\gamma \ge 1$ such that for all $M\times M$ 2D mesh graphs $\graph$ there exists a Hamiltonian path $\dirpath$ such that if $(A,\graph)$ is GIRS-$c$, then $(A,\graph)$ possesses a $\graph$-SS representation with $r_i^g,r_i^h \le c\gamma$ for every $i \in \nodeset_\graph$.
\end{conjecture}

The construction of the minimal DV and $\graph$-SS representations for a general graph $\graph$ is a hard problem, which we demonstrated by example by constructing the minimal $\graph$-SS representation for the GIRS matrices on the cycle graph. The techniques used could potentially be used to prove or disprove Conjecture \ref{conj:edv_compact}, or at least provide heuristics for constructing $\graph$-SS representations for general graphs. As a next step, we plan to study $\graph$-SS representations for graphs $\graph$ consisting of a cycle with the addition of one or more edges, followed by the 2D and 3D mesh graphs.

\begin{figure}
\begin{center}
\begin{tikzpicture}[>=triangle 60]
  \matrix[matrix of math nodes, column sep={140pt,between origins},row
    sep={60pt,between origins},nodes={asymmetrical rectangle}] (s)
  {
     |[name=CEDVI]| \stackrel{\mbox{Compact $\graph$-SS}}{\mbox{for inverse}} &|[name=CDVI]| \stackrel{\mbox{Compact DV}}{\mbox{for inverse}} & |[name=IGIRS]| \stackrel{\mbox{Inverse is}}{\mbox{GIRS}} \\
    |[name=CEDV]| \mbox{Compact $\graph$-SS} &|[name=CDV]| \mbox{Compact DV} & |[name=GIRS]| \mbox{GIRS}\\
    |[name=GIRSLC]| \mbox{GIRS (line or cycle)} &|[name=SPARSE]| \mbox{Sparse Matrix} & \\
    |[name=GEN]| \stackrel{\mbox{General Matrix}}{\scriptsize\graph \mbox{ vertex-disjoint path connected}}\\
    |[name=NEDV]| \mbox{Noncompact $\graph$-SS} &|[name=NDV]| \mbox{Noncompact DV} \\
  };
  \draw[->] (GIRSLC) edge node[auto] {$\stackrel{\mbox{Theorem \ref{thm:SSStoGIRS}}}{\mbox{Theorem \ref{thm:CSS_construction}}}$} (CEDV) 
            (CDV) edge node[auto] {Proposition \ref{prop:idv_algebra}} (CDVI) 
            (GIRS) edge node[auto] {Proposition \ref{prop:girs_algebra}}  (IGIRS) 
            (CEDV) edge node[auto] {Eq.\eqref{eq:GSSISDV}} (CDV) 
            (CDV) edge node[auto] {Proposition \ref{prop:idv_girs}} (GIRS)
            (SPARSE) edge node[auto] {Example \ref{prop:sparse_idv}} (CDV)
            (GEN) edge node[auto] {Proposition \ref{prop:edv_universal}} (NEDV)
            (NEDV) edge node[auto] {Eq.\eqref{eq:GSSISDV} } (NDV)
  ;
  \draw[->,gray,dotted] (GIRS) edge[out=270,in=270] node[auto] {Conjecture \ref{conj:edv_compact}} (CEDV)
                        (GIRS) edge[out=270,in=270] node[auto] {Conjecture \ref{conj:idv_compact}} (CDV)
                        (IGIRS) edge[out=150,in=30] node[auto] {Conjecture \ref{conj:edv_compact}} (CEDVI) 
  ;
\end{tikzpicture}
\end{center}
\caption{Summary of relations between GIRS matrices, $\graph$-SS, and DV representations. In the above, a ``compact'' representation means we have some non-trivial \textit{a priori} bound on the size of the representation and saying a matrix is GIRS implies we have a bound on the constant $c$.}\label{fig:girs_summary}
\end{figure}
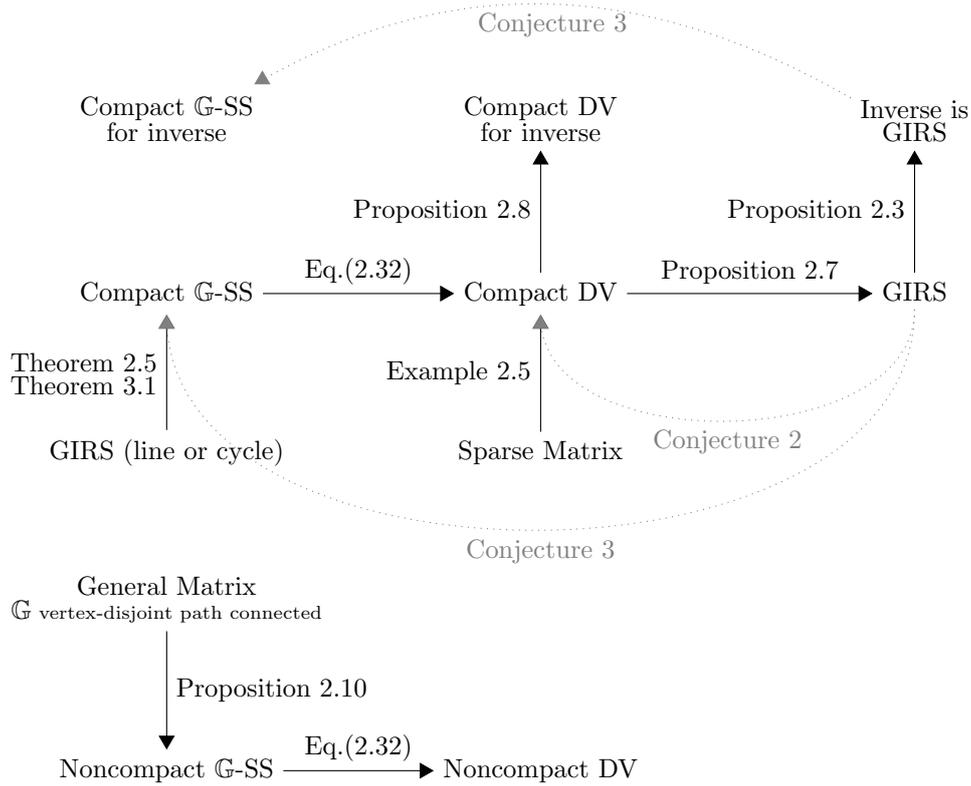

Still, with many important questions remaining on the table, we have established many important connections between GIRS matrices and $\graph$-SS and DV representations, which we summarize in Figure \ref{fig:girs_summary}.

\bibliographystyle{plain}
\bibliography{references}

\appendix

\notrunc{
\section{Details for Construction}\label{app:numerical_details}

While the proof of Theorem \ref{thm:rankcompletion} is constructive, many details necessary to implement this construction in computer code (such as the entire proof of Lemma \ref{lem:removal_of_rows}, do to it’s similarity to Lemma \ref{lem:addition_of_cols}) are missing in our exposition. We spell these details out here.

There are three major calculations to elaborate on. First is the precomputation step in which we perform linear combinations of the rows and columns of the triangular array to ensure the conditions $R(B) \cap R(G) = \{0\}$ and $R(B^T) \cap R(F^T) = \{0\}$.
The next steps are the removal of rows (Lemma \ref{lem:removal_of_rows}) and addition of columns (Lemma \ref{lem:addition_of_cols}).

\subsection{Intersection of subspaces}\label{sec:intersection_of_subspaces}

All of these calculations (as well as the solution of the $2\times 2$ low rank completion problem in Section \ref{sec:2x2_low_rank_completion}) rely very heavily on computing the intersection of subspaces. Subspace computations are inherently unstable as, for any linearly independent vectors $x$ and $y$, $R(x) \cap R(x+\epsilon y) = \{0\}$ if, and only if, $\epsilon = 0$—the intersection of subspaces can be changed by an arbitrarily small perturbation. In view of this, the numerically sensible thing to do is to compute the so-called principal angles $\theta_1 \le \theta_2 \le \cdots \le \theta_k$ between the subspaces, and consider the intersection of the subspaces to be spanned by the first $\ell$ principal vectors such that $\theta_\ell \le \epsilon$ for some tolerance $\epsilon>0$. The principal angles $\theta_1,\ldots,\theta_k$ and principal vectors $u_1,\ldots,u_k$ and $v_1,\ldots,v_k$ between subspaces $U$ and $V$ are defined recursively by
\begin{equation*}
    u_\ell,v_\ell = \argmax_{\stackrel{u\in \mathcal{B}(U_\ell^\perp)}{v\in \mathcal{B}(V_\ell^\perp)}} u^Tv,
\end{equation*}
where $U_\ell = R \left(\begin{bmatrix} u_1 & \cdots & u_\ell \end{bmatrix}\right)$, $U_\ell^\perp$ is the orthogonal complement of $U_\ell$ in $U$, and $\mathcal{B}(U_\ell^\perp)$ denotes the unit ball in $U_\ell^\perp$.

We compute the principal angles and directions using the classical algorithm of Bj{\"o}rk and Golub \cite{bjorck1973numerical}.
In particular, we use SVDs for all rank factorizations to handle any potential rank deficiencies as detailed in \cite[Sec.~5]{bjorck1973numerical} and use the method described in \cite[Sec.~2]{bjorck1973numerical} to compute the \textit{sines} of the principal angles rather then the cosines. (This is better numerically since, for small $\theta$, $\cos(\theta) = 1 - \theta^2/2$, which is the addition of a very small $O(\theta^2)$ perturbation to a large constant, and is thus prone to severe cancellation issues.)
Using these techniques, we are able to compute for given matrices $A$ and $B$, nonsingular matrices $S_A$ and $S_B$ and full column rank matrices $\PAbar$, $\PAB$, and $\PBbar$ such that:
\begin{enumerate}
    \item $AS_A = \begin{bmatrix} \PAB & \PAbar & 0 \end{bmatrix}$ and $BS_B = \begin{bmatrix}
    \PAB & \PBbar & 0
    \end{bmatrix}$,
    \item $\PAB$ is a basis for $R(A) \cap R(B)$, and
    \item $R(A) = R(\PAB) \oplus R(\PAbar)$ and $R(B) = R(\PAB) \oplus R(\PBbar)$.
\end{enumerate}
To be explicit, we describe the details of the algorithm completely here below. 

\begin{alg}[Construction of $\PAbar$, $\PAB$, and $\PBbar$] Given two matrices $A\in \mathbb{R}^{m\times n_A}$ and $ B\in \mathbb{R}^{m\times n_B} $, execute the following steps:
\begin{enumerate}
    \item  Compute SVDs $A = P_A D_A Q_A^T$ and $B = P_B D_B Q^T_B$. Then $P_AP_A^T$ is an orthogonal projector onto $R(A)$, so $(I-P_AP_A^T)$ is a orthogonal projector onto $(R(A))^\perp$.
    \item Define $N := (I-P_AP_A^T)P_B$, and compute an SVD $N = P_N S Q_N^T$. The diagonal entries of $S$ are the sines of the principal angles between $R(A)$ and $R(B)$, $S = \diag(\sin \theta_1,\ldots,\sin \theta_k)$. 
    \item  If we take a nonconventional ordering and assume that the singular values in $S$ are ordered least to greatest, we then can find $\ell$ such that $0\le \theta_1 \le \cdots \le \theta_\ell \le \epsilon < \theta_{\ell+1} \le \cdots \le \theta_k$.
    \item Define $M := P_A^TP_B$ and compute an SVD $M = P_M C Q_M^T$. Here, $C$ contains the cosines of the principal angles, which we take to have the usual ordering $\cos \theta_1 \ge \cdots \ge \cos\theta_k$.
    \item Determine the numerical ranks $r_A$ and $r_B$ of $A$ and $B$, and set
    \begin{equation*}
        D_A^{\mbox{inv}} = \diag(\sigma_1^{-1},\ldots,\sigma_{r_A}^{-1},1,\ldots,1)
    \end{equation*}
    and likewise for $B$.
    \item Define $S_A := Q_A D_A^{\mbox{inv}} \twobytwo{P_M}{0}{0}{I}$ and $S_B = Q_B D_B^{\mbox{inv}} \twobytwo{Q_M}{0}{0}{I}$, where the identity matrices are of the appropriate sizes.
\end{enumerate}
The first $k-\ell$ columns of $AS_A$ and $BS_B$ is defined to be $\PAB$. The next $r_A + \ell - k$ and $r_B + \ell - k$ columns of $AS_A$ and $BS_B$ are $\PAbar$ and $\PBbar$ respectively, where $r_A$ and $r_B$ are the numerical ranks of $A$ and $B$.
\end{alg}
\begin{remark}
Often, we shall often need also to determine matrices $\QABA$, $\QAbarA$, $\QABB$, and $\QBbarB$ such that
\begin{equation*}
    A = \onebytwo{\PAB}{\PAbar} \twobyone{\QABA}{\QAbarA}, \quad B = \onebytwo{\PAB}{\PBbar}\twobyone{\QABB}{\QBbarB}.
\end{equation*}
These can be found be simply looking at the appropriate rows of $S_A^{-1}$ and $S_B^{-1}$. By applying this algorithm to $A^T$ and $B^T$ this approach can be used to compute the interesection of the row spaces instead, 
\end{remark}

\newcommand{\QAB}{\Qs{\intersect{A}{B}}}
\newcommand{\QAbar}{\Qs{\comp{A}}}

\begin{equation*}
    S_A A = \begin{bmatrix}
    \QAB \\ \QAbar \\ 0
    \end{bmatrix}, \quad S_B B = \begin{bmatrix}
    \QAB \\ \QBbar \\ 0
    \end{bmatrix}
\end{equation*}

\subsection{Precomputation}

Our goal is elaborate on the ``precomputation'' step in the proof of Theorem~\ref{thm:rankcompletion}, in which the original Hankel blocks $\mathcal{H}_k(X)$ are modified to new Hankel blocks $\mathcal{K}_k(X)$ such that 

\begin{enumerate}
    \item For every choice of $X$ and $1\le k \le n-1$, $\rank \mathcal{H}_k(X) = \rank \mathcal{K}_k(X)$. In particular, $\hat{X}$ minimizes the rank of all Hankel blocks $\mathcal{H}_k(X)$ if, and only if, it minimizes the rank of all Hankel blocks $\mathcal{K}_k(X)$.
    \item For every $1\le k \le n-2$, $R((B^k)^T) \cap R((F^k)^T) = \{0\}$ and $R(B^k) \cap R(G^k) = \{0\}$,
\end{enumerate}
where
\begin{equation*}
    \mathcal{K}_k(X) = \begin{bmatrix}
        E^k & F^k \\
        A^k & B^k \\
        X & C^k
        \end{bmatrix}, \quad \mathcal{K}_{k+1}(X) = \begin{bmatrix}
        A^k & B^k & G^k \\
        X & C^k & H^k
        \end{bmatrix}.
\end{equation*}
We shall use the same notation as the proof of Theorem~\ref{thm:rankcompletion}, including the partitioning laid out in \eqref{eq:hankel_partitioning}. The one detail which deserves elaboration is the choice of the matrices $S_1$ and $S_2$. The matrix $S_1$ is to be chosen such that for $\onebytwo{\tilde{E}^k}{\tilde{F}^k} := \onebytwo{\hat{E}^k}{\hat{F}^k} + S_1\onebytwo{\hat{A}^k}{\hat{B}^k}$, $R((\tilde{F}^k)^T) \cap R((\hat{B}^k)^T) = \{0\}$. The construction of $S_2$ is identical if we consider the tranpose of the procedure to construct $S_1$, so we shall focus on $S_1$.

To simplify our notation, we suppress the dependence on $k$ and the use of ornamentations, and write $E' := \tilde{E}^k$, $F' := \tilde{F}^k$, $E := \hat{E}^k$, $F := \hat{F}^k$, $A := \hat{A}^k$, and $B := \hat{B}^k$. First, use the intersection of row spaces algorithm to find matrices $S_B$ and $S_F$ such that
\newcommand{\QFB}{\Qs{\intersect{F}{B}}}
\newcommand{\QjustFbar}{\Qs{\comp{F}}}
\begin{equation*}
    S_BB = \begin{bmatrix} \QFB \\ \QBbar \\ 0 \end{bmatrix}, \quad S_FF = \begin{bmatrix} \QFB \\ \QjustFbar \\ 0 \end{bmatrix}
\end{equation*}
Thus 
\begin{equation*}
    F +\underbrace{\left[-S_F^{-1}\begin{bmatrix} I&& \\ &0& \\ &&0\end{bmatrix}S_B\right]}_{:=S_1}B = \underbrace{S_F^{-1}\begin{bmatrix} 0 \\ \QjustFbar \\ 0 \end{bmatrix}}_{:=F'}.
\end{equation*}
By construction, it is easy to see that $R(F'^T) \cap R(B^T) = 0$. 

\subsection{Removal of rows}\label{sec:removal_of_rows}

In this section, we provide a detailed constructive proof of Lemma \ref{lem:removal_of_rows}, amenable to direct implementation in a programming language. We make no claims that this is an efficient way of carrying out this calculation, but hopefully provide enough details to produce working code. For expositional clarity, we distinguish different stages of the construction by numbered steps. Throughout, we shall denote by $\ast$ a matrix of the appropriate size whose precise value will be immaterial in the ensuing calculation. We denote by $S^R$ and $S^L$ an arbitrary right and left inverse of a matrix $S$ possessing full row or column rank, $SS^R = I$ and $S^LS = I$. 

Suppose we have the complete solution set to a low rank completion problem
\begin{equation*}
    \begin{bmatrix}
    E & F \\
    A & B \\
    ? & C
    \end{bmatrix}
\end{equation*}
in terms of the six complementary subspaces $\PEbar$, $\PEF$, $\PFbar$, $\QFbar$, $\QFC$, and $\QCbarbullet$. We are interested in a subclass of solutions
\begin{equation*}
    \mathscr{S}_{\EA,\FB,C}^Y = \{ G\QEbarE + \PCFC \RFCEF \QEFE + \PCCbarbullet Y : G \textnormal{ free} \}. 
\end{equation*}
We seek six bases $\PAbar$, $\PAB$, $\PBbar$, $\QBbar$, $\QBC$, and $\QCbar$ characterizing the low rank completion problem solutions of $\twobytwo{A}{B}{?}{C}$ and a matrix $Z$ for which $\mathscr{S}_{A,B,C}^Z \subseteq \mathscr{S}_{\EA,\FB,C}^Y$.

We now begin the construction in earnest.

\paragraph{1. Determination of $\PAbar$ and $\PAB$ from $\PEbar$ and $\PEF$.} Partition $\PEF$, $\PEbar$, and $\PFbar$ as 
\begin{equation*}
    \PEF  = \twobyone{\PEF '}{\PAB '}, \quad \PEbar  = \twobyone{\PEbar '}{\PAbar '}, \quad \PFbar  = \twobyone{\PFbar '}{\PBbar '}.
\end{equation*}
First, find a nonsingular matrix $S_0$ such that
\begin{equation*}
    \PEF S_0 = \begin{bmatrix}
    \PEF^1 & \PEF^2 \\
    \PAB^1 & 0
    \end{bmatrix}.
\end{equation*}
(Use an SVD to write $\PAB' = \onebytwo{\PAB^1}{\ast} \twobytwo{\Sigma}{0}{0}{0} Q^T$ and define $S_0 := Q\twobytwo{\Sigma^{-1}}{0}{0}{I}$.) Use the intersection of column spaces algorithm on $\PAbar'$ and $\PAB'$ to find matrices $S_1$ and $S_2$ such that
\begin{equation*}
    \PAbar'S_1 = \begin{bmatrix} \Ps{\ast} & \PAbar^\bullet & 0 \end{bmatrix}, \quad \PAB' S_2 = \begin{bmatrix} \Ps{\ast} & \ast & 0 \end{bmatrix}.
\end{equation*}
Then, choosing a matrix $S_3$ such that $\PAB^1 S_2S_3 = \onebytwo{\Ps{\ast}}{\ast}S_3 = \begin{bmatrix} \Ps{\ast} & 0 & 0 \end{bmatrix}$, then there exists a permutation $\pi_1$ such that
\begin{equation*}
    \onebytwo{\PAbar'}{\PAB'}\underbrace{\twobytwo{S_1}{}{-S_2S_3}{S_0}\pi_1}_{:=S_4} = \begin{bmatrix} 0_1 & \PAbar^\bullet & 0_2 & \PAB^1 & 0_3 \end{bmatrix}\pi_1 = \begin{bmatrix} \PAbar^\bullet & 0_1 & 0_2 & \PAB^1 & 0_3 \end{bmatrix},
\end{equation*}
where we notate $0_1$, $0_2$, and $0_3$ as zero matrices of the appropriate sizes. Note that, even with the multiplication of the permutation $\pi_1$, $S_4$ remains block lower triangular. Now, use the intersection of column space algorithm with $\onebytwo{\PAbar^\bullet}{\PAB^1}$ and $B$ to find a nonsingular matrix $S_2$ such that
\begin{equation*}
    \onebytwo{\PAbar^\bullet}{\PAB^1}S_5 = \begin{bmatrix}
    \PAB' & \PAbar
    \end{bmatrix}.
\end{equation*}
Then use the intersection of column space algorithm on $\PAB'$ and $\PAB^1 $ to find a nonsingular matrix $S_3$ such that
\begin{equation*}
    \PAB'S_6 = \begin{bmatrix}
    \PAB^{1\prime }& \PAB^2 \end{bmatrix}.
\end{equation*}
Then
\begin{equation*}
    \PAB'\underbrace{S_6\twobytwo{(\PAB^{1\prime })^L\PAB^1}{}{}{I}}_{:=S_7} = \onebytwo{\PAB^1}{\PAB^2}.
\end{equation*}
Then we have
\begin{equation*}
    \onebytwo{\PAbar^\bullet}{\PAB^1}\underbrace{S_5\twobytwo{S_7}{}{}{I}}_{:= M} = \begin{bmatrix} \PAB^1 & \PAB^2 & \PAbar \end{bmatrix}.
\end{equation*}
Denoting 
\begin{equation*}
    M = \begin{bmatrix} M_{11} & M_{12} & M_{13} \\ M_{21} & M_{22} & M_{23} \end{bmatrix},
\end{equation*}
we see right away that $M_{11} = 0$ and $M_{21} = I$. Then12
\begin{align*}
    \onebytwo{\PAbar'}{\PAB'}\underbrace{S_4 \begin{bmatrix} M_{12} & 0 & M_{13} & 0 & 0 \\
    0 & 0 & 0 & 0 & I \\
    M_{22} & I & M_{23} & 0 & 0 \\
    0 & 0 & 0 & I & 0
    \end{bmatrix}}_{:=\Phi} &= \begin{bmatrix} \PAbar^\bullet & 0_1 & \PAB^1 & 0_2 \end{bmatrix}\begin{bmatrix} M_{12} & 0 & M_{13} & 0 & 0 \\
    0 & 0 & 0 & 0 & I \\
    M_{22} & I & M_{23} & 0 & 0 \\
    0 & 0 & 0 & I & 0
    \end{bmatrix} \\
    &= \begin{bmatrix} \PAB^2 & \PAB^1 & \PAbar & 0_2 & 0_1 \end{bmatrix},
\end{align*}
where we again notate $0_1$ and $0_2$ for two different zero matrices of potentially different sizes. Define $\PAB := \onebytwo{\PAB^2}{\PAB^1}$. Note that, by construction, $\Phi$ is nonsingular and has the following block nonzero structure:
\begin{equation*}
    \onebytwo{\PAbar '}{\PAB '}\underbrace{\begin{bmatrix}
    \Phi_{11} & 0 & \Phi_{31} & 0 & \Phi_{15} \\
    \Phi_{21} & \Phi_{22} & \Phi_{23} & \Phi_{24} & \Phi_{25}
    \end{bmatrix}}_{\Phi} = \begin{bmatrix}
    \PAB^{2} & \PAB^{1} & \PAbar & 0_2 & 0_1
    \end{bmatrix}.
\end{equation*}
\paragraph{2.  Determination of $\QAbarA$ and $\QABA$.} We have
\begin{align*}
    \twobyone{E}{A} &= \onebytwo{\PEbar}{\PEF}\twobyone{\QEbarE}{\QEFE} = \twobytwo{\PEbar '}{\PEF '}{\PAbar '}{\PAB '}\twobyone{\QEbarE}{\QEFE} = \twobytwo{\PEbar '}{\PEF '}{\PAbar '}{\PAB '} \Phi\Phi^{-1} \twobyone{\QEbarE}{\QEFE} \\
    &= \begin{bmatrix}
    \ast & \PEF^1 & \ast & \PEF^2 & \ast \\
    \PAB^{2} & \PAB^{1}  & \PAbar  & 0 & 0
    \end{bmatrix} \begin{bmatrix}
    \QABAnum{2} \\ \QABAnum{1} \\ \QAbarA \\ \QEFEnum{2} \\ \QstarA
    \end{bmatrix}.
\end{align*}
Define $\QABA := \twobyone{\QABAnum{2}}{\QABAnum{1}}$.

\paragraph{3.  Determination of $\PBbar$, $\QABB$, and $\QBbar$.} Next, use the intersection of column spaces on $\PBbar'$ and $\PAB'$ to find a matrices $S_8$ and $S_9$ such that
\begin{equation*}
    \PBbar'S_8 = \begin{bmatrix} \Ps{\ast\ast} & \PBbar^\bullet & 0 \end{bmatrix}, \quad \PAB'S_9 = \begin{bmatrix} \Ps{\ast\ast} & \ast & 0 \end{bmatrix}.
\end{equation*}
Then, similar to $S_3$ and $\pi_1$ above, there exists a matrix $S_{10}$ such that $\PAB^1S_9S_{10} = \begin{bmatrix} \Ps{\ast\ast} & 0 & 0 \end{bmatrix}$ and a permutation $\pi_2$ such that
\begin{equation*}
    \onebytwo{\PBbar'}{\PAB'}\underbrace{\twobytwo{S_{8}}{}{-S_9S_{10}}{S_0}\pi_2}_{:=S_{11}} = \begin{bmatrix} 0_1 & \PBbar^\bullet & 0_2 & \PAB^1 & 0_3 \end{bmatrix}\pi_2 = \begin{bmatrix} \PBbar^\bullet & 0_1 & 0_2 & \PAB^1 & 0_3 \end{bmatrix}
\end{equation*}
Then use the intersection of column space algorithm with $\onebytwo{\PBbar^\bullet}{\PAB^1}$ and $A$ to find a nonsingular matrix $S_{12}$ such that
\begin{equation*}
    \onebytwo{\PBbar^\bullet}{\PAB^1}S_{12} = \onebytwo{\PAB''}{\PBbar}.
\end{equation*}
Consequently, noting that $\PAB$, $\PAB''$, and $\onebytwo{\PAB^1}{\PAB^2}$ are all column bases for the same space, namely $R(A) \cap R(B)$, we conclude
\begin{equation*}
    \onebytwo{\PBbar^\bullet}{\PAB^1}\underbrace{S_{12}\twobytwo{(\PAB'')^L\onebytwo{\PAB^1}{\PAB^2}}{}{}{I}}_{:=N} = \begin{bmatrix} \PAB^1 & \PAB^2 & \PBbar \end{bmatrix}, \quad N = \begin{bmatrix} 0 & N_{12} & N_{13} \\ 
    I & N_{22} & N_{23}\end{bmatrix}.
\end{equation*}
Thus
\begin{align*}
    \onebytwo{\PBbar'}{\PAB'}\underbrace{S_{11}\begin{bmatrix} N_{12} & 0 & N_{13} & 0 & 0 \\
    0 & 0 & 0 & 0 & I \\
    N_{22} & I & N_{23} & 0 & 0 \\
    0 & 0 & 0 & I & 0
    \end{bmatrix}}_{:=\Psi} &= \begin{bmatrix} \PBbar^\bullet & 0_1 & \PAB^1 & 0_2 \end{bmatrix}\begin{bmatrix} N_{12} & 0 & N_{13} & 0 & 0 \\
    0 & 0 & 0 & 0 & I \\
    N_{22} & I & N_{23} & 0 & 0 \\
    0 & 0 & 0 & I & 0
    \end{bmatrix} \\
    &= \begin{bmatrix} \PAB^2 & \PAB^1 & \PAbar & 0_2 & 0_1 \end{bmatrix}.
\end{align*}
We note that by construction, $\Psi$ has the block nonzero structure
\begin{equation*}
    \onebytwo{\PBbar '}{\PAB '}\underbrace{\begin{bmatrix}
    \Psi_{11} & 0 & \Psi_{13} & 0 & \Psi_{15} \\
    \Psi_{21} & \Psi_{22} & \Psi_{23} & \Psi_{24} & \Psi_{25}
    \end{bmatrix}}_{=\Psi} = \begin{bmatrix}
    \PAB^{2} & \PAB^{1} & \PBbar  & 0_2 & 0_1
    \end{bmatrix}.
\end{equation*}
%
%
Then
\begin{align*}
    \twobyone{F}{B} &= \onebytwo{\PFbar }{\PEF }\twobyone{\QFbarF}{\QEFF} = \onebytwo{\PFbar }{\PEF }\Psi\Psi^{-1}\twobyone{\QFbarF}{\QEFF} \\
    &= \begin{bmatrix}
    \ast & \PEF^1 & \ast & \PEF^2 & \ast \\
    \PAB^{2} & \PAB^{1}   & \PBbar  & 0 & 0
    \end{bmatrix}\begin{bmatrix}
    \QABBnum{2} \\ \QABBnum{1} \\ \QBbarB \\ \QEFFnum{2} \\ \QstarB
    \end{bmatrix}.
\end{align*}
Define $\QABB = \twobyone{\QABBnum{2}}{\QABBnum{1}}$.
\ethan{\textbf{Given}: $B, F, C$ such that $R(B^T) \cap R(F^T) = \{0\}$ and row bases $\QFC$ for $R(\twobyone{F}{B}^T) \cap R(C^T)$ and $\twobyone{\QFC}{\QFbar}$ for $\twobyone{F}{B}$ \textbf{find} row bases $\QBC$ for $R(B^T)\cap R(C^T)$, $\twobyone{\QBC}{\QBbar}$ for $R(B^T)$, $\twobyone{\QBC}{\QFCnum{2}}$ for $R(\twobyone{F}{B}^T) \cap R(C^T)$, and $\begin{bmatrix} \QBC \\ \QFCnum{2} \\ \QBbar \\ \QFbarnum{2}\end{bmatrix}$ for $R(\twobyone{F}{B}^T)$.}
\paragraph{4.  Determination of $\QBbar$, $\QBC$, $\QCbar$, $\PBBbar$, $\PBBC$, $\PCBC$, $\PCCbar$.} We now consider $\begin{bmatrix}
F \\ B \\ C
\end{bmatrix}$. First, use the intersection of row spaces algorithm with $B$ and $C$ to find a nonsingular matrix $T_1$ such that
\begin{equation*}
    T_1B = \begin{bmatrix} \QBC \\ \QBbar  \\ 0 \end{bmatrix}.
\end{equation*}
Next, use the row intersection algorithm with $\QFC $ and $\QBC$ to find a matrix $T_2$ such that
\begin{equation*}
    T_2\QFC  = \twobyone{\Qs{\ast} }{\QFCnum{2}}.
\end{equation*}
Then, since $R(\QFC) \cap R(\QBC) = R(\QBC)$, we have
\begin{equation*}
    T_3 = \twobyone{\QBC}{\QFCnum{2}}\left(\twobyone{\Qs{\ast}}{\QFCnum{2}}\right)^RT_2, \mbox{ we have }T_3\QFC = \twobyone{\QBC}{\QFCnum{2}}.
\end{equation*}
We then use the row intersection algorithm again with $\twobyone{F}{B}$ and $\begin{bmatrix} \QFC \\ \QBbar \end{bmatrix}$ to find a matrix $T_4$ such that
\begin{equation*}
    T_4\twobyone{F}{B}  = \begin{bmatrix} \ast \\ \QFbarnum{2} \\ 0 \end{bmatrix}.
\end{equation*}
Since both $\begin{bmatrix}
    \QBbar  \\ \QFbarnum{2}  \\ \QBC  \\ \QFCnum{2}
    \end{bmatrix}$ and
$\twobyone{\QFbar }{\QFC }$ are row bases for $\begin{bmatrix} F \\ B \\ C \end{bmatrix}$, they are related by a nonsingular matrix $\Omega$:
\begin{equation*}
    \underbrace{\begin{bmatrix}
    \Omega_{11} & \Omega_{12} \\
    0 & \Omega_{22}
    \end{bmatrix}}_{=\Omega}\twobyone{\QFbar }{\QFC } =
    \begin{bmatrix} \begin{bmatrix} \QBbar  \\ \QFbarnum{2}\end{bmatrix}  \\ \vspace{-2mm} \\ \begin{bmatrix} \QBC  \\ \QFCnum{2} \end{bmatrix} \end{bmatrix}, \mbox{ where } \Omega_{22} = T_3.
\end{equation*}
Notice that we have
\begin{align*}
    C &= \onebytwo{\PCFC }{\PCCbarbullet }\twobyone{\QFC }{\QCbarbullet} = \onebytwo{\PCFC T_3^{-1}}{\PCCbarbullet }\twobyone{\twobyone{\QBC }{\QFCnum{2} }}{\QCbarbullet } \\
    &= \onebytwo{\PCBC }{\underbrace{\onebytwo{\PCCbarnum{1}}{\PCCbarbullet }}_{:=\PCCbar}}\twobyone{\QBC }{\QCbar },
\end{align*}
where $\QCbar := \twobyone{\QFCnum{2}}{\QCbarbullet}$.
%


%
Thus
\begin{align*}
    \twobyone{F}{B} &= \onebytwo{\PFFbar }{\PFFC }\twobyone{\QFbar }{\QFC } = \onebytwo{\PFFbar }{\PFFC }\Omega^{-1}\Omega \twobyone{\QFbar }{\QFC } \\
    &= \begin{bmatrix}
    \ast & \ast & \ast & \ast \\
    \PBBbar  & 0 & \PBBC & 0
    \end{bmatrix}\begin{bmatrix}
    \QBbar  \\ \QFbarnum{2}  \\ \QBC  \\ \QFCnum{2}
    \end{bmatrix}
\end{align*}
\paragraph{5.  Studying the relations between the ``$R$'' matrices for the new and old problem.} We have
\newcommand{\permutedblockRremovalofrows}{\twobytwo{\RFbarFbar}{\RFbarEF}{\RFCFbar}{\RFCEF}}
\begin{align*}
    \twobyone{F}{B} &= \onebytwo{\PFFbar }{\PFFC } \twobytwo{\RFbarEF}{\RFbarFbar}{\RFCEF}{\RFCFbar}\twobyone{\QEFF }{\QFbarF} \\
    &= \onebytwo{\PFFbar }{\PFFC } \permutedblockRremovalofrows \twobyone{\QFbarF}{\QEFF}. \\
\end{align*}
Inserting nonsingular transformations, we may write
\begin{align*}
    \twobyone{F}{B} &= \onebytwo{\PFFbar }{\PFFC }\Omega^{-1}\Omega \permutedblockRremovalofrows\Psi\Psi^{-1}\twobyone{\QFbarF}{\QEFF} \\
    &= \begin{bmatrix}
    \begin{bmatrix}
    \ast & \ast \\
    \PBBbar  & 0
    \end{bmatrix} & \begin{bmatrix}
    \ast & \ast \\
    \PBBC & 0
    \end{bmatrix}
    \end{bmatrix}\twobytwo{\Omega_{11}}{\Omega_{12}}{}{T_3}\permutedblockRremovalofrows\Psi\begin{bmatrix}
    \QABBnum{2} \\ \QABBnum{1} \\ \QBbarB \\ \QEFFnum{2}  \\ \QstarB
    \end{bmatrix}. \\
\end{align*}    
We now use the the known structure of $\Psi$:
\begin{align*}
    \twobyone{F}{B} &= \begin{bmatrix}
    \begin{bmatrix}
    \ast & \ast \\
    \PBBbar  & 0
    \end{bmatrix} & \begin{bmatrix}
    \ast & \ast \\
    \PBBC & 0
    \end{bmatrix}
    \end{bmatrix}\twobytwo{\Omega_{11}}{\Omega_{12}}{}{T_3}\permutedblockRremovalofrows\\
    &\qquad \times \begin{bmatrix}
    \Psi_{11} & 0 & \Psi_{13} & 0 & \Psi_{15} \\
    \Psi_{21} & \Psi_{22} & \Psi_{23} & \Psi_{24} & \Psi_{25}
    \end{bmatrix}\begin{bmatrix}
    \QABBnum{2} \\ \QABBnum{1} \\ \QBbarB \\ \QEFFnum{2} \\ \QstarB
    \end{bmatrix} \\
\end{align*}
This yields 
\begin{equation*}
    \twobyone{F}{B} = \begin{bmatrix}
    \begin{bmatrix}
    \ast & \ast \\
    \PBBbar  & 0
    \end{bmatrix} & \begin{bmatrix}
    \ast & \ast \\
    \PBBC & 0
    \end{bmatrix}
    \end{bmatrix}\begin{bmatrix}
    \Gamma_{11} & \Gamma_{12} & \Gamma_{13} & \Gamma_{14} & \Gamma_{15} \\
    \Gamma_{21} & \Gamma_{22} & \Gamma_{23} & \Gamma_{24} & \Gamma_{25} 
    \end{bmatrix} \begin{bmatrix}
    \QABBnum{2} \\ \QABBnum{1} \\ \QBbarB \\ \QEFFnum{2}  \\ \QstarB
    \end{bmatrix},
\end{equation*}
where
\begin{align*}
    \Gamma_{11} &= \Omega_{11}(\RFbarFbar \Psi_{11} + \RFbarEF \Psi_{21}) + \Omega_{12}(\RFCFbar \Psi_{11} + \RFCEF \Psi_{21}) \\
    \Gamma_{12} &= \Omega_{11}\RFbarEF \Psi_{22} + \Omega_{12}\RFCEF \Psi_{22} \\
    \Gamma_{13} &= \Omega_{11}(\RFbarFbar \Psi_{13} + \RFbarEF \Psi_{23}) + \Omega_{12}(\RFCFbar \Psi_{13} + \RFCEF \Psi_{23}) \\
    \Gamma_{14} &= \Omega_{11}\RFbarEF \Psi_{24} + \Omega_{12}\RFCEF \Psi_{24} \\
    \Gamma_{15} &= \Omega_{11}(\RFbarFbar \Psi_{15} + \RFbarEF \Psi_{25}) + \Omega_{12}(\RFCFbar \Psi_{15} + \RFCEF \Psi_{25}) \\
    \Gamma_{21} &= T_3(\RFCFbar \Psi_{11} +\RFCEF \Psi_{21}) \\
    \Gamma_{22} &= T_3\RFCEF \Psi_{22} \\
    \Gamma_{23} &= T_3(\RFCFbar \Psi_{13} +\RFCEF \Psi_{23}) \\
    \Gamma_{24} &= T_3\RFCEF \Psi_{24} \\
    \Gamma_{25} &= T_3(\RFCFbar \Psi_{15} + \RFCEF \Psi_{25}). \\
\end{align*}
\newcommand{\tGamma}{\tilde{\Gamma}}
Denote by $\tilde{X}$ the matrix $X$ multiplied by $\onebytwo{I}{0}$. For instance,  $\tilde{\Omega}_{11} := \onebytwo{I}{0}\Omega_{11}$, $\tGamma_{11} = \onebytwo{I}{0}\Omega_{11}$, and $\tilde{T}_3 := \onebytwo{I}{0}T_3$, etc. Then
\begin{align*}
    B &= \onebytwo{\PBBbar }{\PBBC}\begin{bmatrix}
    \tGamma_{11} & \tGamma_{12} & \tGamma_{13} & \tGamma_{14} & \tGamma_{15} \\
    \tGamma_{21} & \tGamma_{22} & \tGamma_{23} & \tGamma_{24} & \tGamma_{25} 
    \end{bmatrix}\begin{bmatrix}
    \QABBnum{2} \\ \QABBnum{1} \\ \QBbarB \\ \QEFFnum{2}  \\ \QstarB
    \end{bmatrix} \\
    &= \onebytwo{\PBBbar }{\PBBC}\begin{bmatrix}
    \tGamma_{11} & \tGamma_{12} & \tGamma_{13} \\
    \tGamma_{21} & \tGamma_{22} & \tGamma_{23}
    \end{bmatrix}\begin{bmatrix}
    \QABBnum{2} \\ \QABBnum{1} \\ \QBbarB
    \end{bmatrix},
\end{align*}
and we conclude $\tGamma_{14} = 0$, $\tGamma_{15} = 0$, $\tGamma_{24} = 0$, and  $\tGamma_{25} = 0$. Thus, recalling $\QABB = \twobyone{\QABBnum{2}}{\QABBnum{1}}$, we have
\begin{equation*}
    \twobytwo{\RBbarAB}{\RBbarBbar }{\RBCAB }{\RBCBbar } = \twobytwo{ \onebytwo{ \tGamma_{11} } {\tGamma_{12}}} {\tGamma_{13}}{ \onebytwo{ \tGamma_{21} } {\tGamma_{22}}} {\tGamma_{23}}.
\end{equation*}
\paragraph{6.  Comparison of solution sets.} Thus the complete solution set to the new problem is $\mathscr{S}_{A,B,C} = \{ X^{A,B,C}_{\FAbar,\FCbar}  : \FAbar, \FCbar \mbox{ free}\}$ where
\begin{align*}
    X^{A,B,C}_{\FAbar,\FCbar} &= \FAbar \QAbarA + \PCBC \RBCAB \QABA + \PCCbar \FCbar \\
    &= \FAbar \QAbarA + \PCBC \onebytwo{\tGamma_{21}}{\tGamma_{22}}\twobyone{ \QABAnum{2}}{\QABAnum{1}} + \onebytwo{\PCCbarnum{1}}{\PCCbarbullet } \twobyone{\FCbarnum{1}}{\FCbarnum{2}}.
\end{align*}
And the complete solution set to the original problem is $\mathscr{S}_{\EA,\FB,C} = \{ X^{\EA,\FB,C}_{\FEbar,\FCbarbullet} : \FEbar, \FCbarbullet \mbox{ free} \}$ where
\begin{align*}
    X^{\EA,\FB,C}_{\FEbar,\FCbarbullet} &= \FEbar \QEbarE + \PCFC \RFCEF \QEFE  + \PCCbarbullet  \FCbarbullet \\
    &= \FEbar \QEbarE + \onebytwo{\PCBC }{\PCCbarnum{1}}T_3\RFCEF \QEFE + \PCCbarbullet \FCbarbullet \\
    &= \onebytwo{\FEbar}{\onebytwo{\PCBC }{\PCCbarnum{1}}T_3\RFCEF}\twobyone{Q^T_{\bar{E},E}}{\QEFE }+ \PCCbarbullet \FCbarbullet \\
    &= \onebytwo{\FEbar}{\onebytwo{\PCBC }{\PCCbarnum{1}}T_3\RFCEF}\begin{bmatrix}
    \Phi_{11} & 0 & \Phi_{13} & 0 & \Phi_{15} \\
    \Phi_{21} & \Phi_{22} & \Phi_{23} & \Phi_{24} & \Phi_{25}
    \end{bmatrix}\begin{bmatrix}
    \QABAnum{2} \\ \QABAnum{1} \\ \QAbarA \\ \QEFEnum{2} \\ \QstarA
    \end{bmatrix}\\
    &\qquad + \PCCbarbullet \FCbarbullet \\
    &= \FEbar \begin{bmatrix}
    \Phi_{11} & \Phi_{13} & \Phi_{15}
    \end{bmatrix}\begin{bmatrix}
    \QABAnum{2} \\ \QAbarA \\ \QstarA
    \end{bmatrix}+ \PCCbarbullet  \FCbarbullet  \\
    &\qquad+ \onebytwo{\PCBC }{\PCCbarnum{1}}T_3\RFCEF\begin{bmatrix}\Phi_{21} & \Phi_{22} & \Phi_{23} & \Phi_{24} & \Phi_{25} \end{bmatrix}\begin{bmatrix}
    \QABAnum{2} \\ \QABAnum{1} \\ \QAbarA \\ \QEFEnum{2} \\ \QstarA
    \end{bmatrix}.
\end{align*}
We note that by construction $\onebytwo{\Phi_{22}}{\Phi_{24}} = \onebytwo{\Psi_{22}}{\Psi_{24}} = S_0$ and $\begin{bmatrix}
    \Phi_{11} & \Phi_{13} & \Phi_{15}
    \end{bmatrix}$
is a nonsingular matrix. Make the reparametrization
\begin{equation*}
    \FEbar' = \FEbar \begin{bmatrix}
    \Phi_{11} & \Phi_{13} & \Phi_{15}
    \end{bmatrix}  + \onebytwo{\PCBC }{\PCCbarnum{1}}T_3\RFCEF  \begin{bmatrix}
    \Phi_{21} & \Phi_{23} & \Phi_{25}.
    \end{bmatrix}
\end{equation*}
Denoting $\hat{T}_3 := \onebytwo{0}{I}T_3$, we have
\begin{align*}
    X^{\EA,\FB,C}_{\FEbar,\FCbarbullet} &= \FEbar'\begin{bmatrix}
    \QAbarA \\ \QABAnum{2} \\ \QstarA
    \end{bmatrix} + \onebytwo{\PCBC }{\PCCbarnum{1}}T_3\RFCEF\begin{bmatrix}\Psi_{22} & \Psi_{24} \end{bmatrix}\begin{bmatrix} \QABAnum{1} \\ \QEFEnum{2} \end{bmatrix}+ \PCCbarbullet  \FCbarbullet \\
    &= \FEbar'\begin{bmatrix}
    \QAbarA \\ \QABAnum{2} \\ \QstarA
    \end{bmatrix} + \PCBC \tilde{T}_3\RFCEF \begin{bmatrix}\Psi_{22} & \Psi_{24} \end{bmatrix}\twobyone{\QABAnum{1}}{\QEFEnum{2}}  \\
    &\qquad + \PCCbarnum{1} \hat{T}_3\RFCEF S_0 \begin{bmatrix} \QABAnum{1} \\ \QEFEnum{2} \end{bmatrix} + \PCCbarbullet \FCbarbullet.
\end{align*}
Reparametrizing a final time and using the fact that $\tGamma_{24} = \tilde{T}_3 \RFCEF \Psi_{24} = 0$, we get our solution set
\begin{align*}
    X^{\EA,\FB,C}_{\FEbar,\FCbarbullet} &= \FEbar''\begin{bmatrix}
    \QAbarA \\ \QABAnum{2} \\ \QstarA
    \end{bmatrix} + \PCBC \tGamma_{22}\QABAnum{1} + \PCCbarnum{1} \hat{T}_3 \RFCEF S_0\twobyone{\QABAnum{1}}{\QEFEnum{2}} + \PCCbarbullet \FCbarbullet.
\end{align*}
\paragraph{7.  Determination of $Z$.} Suppose that for our old solution we must have $\FCbarbullet = Y$ to satisfy earlier Hankel blocks. Then in our new solution we must set
\begin{equation*}
    Z = \FCbar = \twobyone{\FCbarnum{1}}{\FCbarnum{2}} = \twobyone{Y_2R_3S_0\twobyone{\QABAnum{1}}{\QEFEnum{2}}}{Y}.
\end{equation*}
One easily verifies that under this choice $\mathscr{S}_{A,B,C}^Z \subseteq \mathscr{S}_{\EA,\FB,C}^Y$.

\subsection{Addition of columns}

Here, we present a detailed constructive proof of Lemma \ref{lem:addition_of_cols}, using the same conventions outlined at the beginning of Section \ref{sec:removal_of_rows}. Before we begin the construction in earnest, we present a more explicit proof of Proposition \ref{prop:extend_to_row_basis}.

\begin{proof}[Proof of Proposition \ref{prop:extend_to_row_basis}, in more detail]

Let $V^T = \onebytwo{V^T_1}{V^T_2}$ be an row basis for $\onebytwo{B}{G}$. Perform the intersection of row spaces on $\twobytwo{\Qs{11}}{\Qs{12}}{0}{\Qs{22}}$ and $V^T$ to find a nonsingular matrix $S_1$ such that

\begin{equation*}
    S_1V^T = \begin{bmatrix}
    V_{11}^T & V_{12}^T \\ V_{31}^T & V_{32}^T
    \end{bmatrix}.
\end{equation*}

Since the row space of $\twobytwo{\Qs{11}}{\Qs{12}}{0}{\Qs{22}}$ is contained in the row space of $V^T$, $\onebytwo{V_{11}^T}{V_{12}^T}$, being a basis for the intersection, must be a basis for the row space of $\twobytwo{\Qs{11}}{\Qs{12}}{0}{\Qs{22}}$, which has full row rank. Thus, there exists a nonsingular matrix 

\begin{equation*}
    S_2 = \twobytwo{\Qs{11}}{\Qs{12}}{0}{\Qs{22}}\left(\onebytwo{V_{11}^T}{V_{12}^T}\right)^R \mbox{ such that } S_2\onebytwo{V_{11}^T}{V_{12}}^T = \twobytwo{\Qs{11}}{\Qs{12}}{0}{\Qs{22}}.
\end{equation*}

Then

\begin{equation*}
    \underbrace{\twobytwo{S_2}{}{}{I}S_1}_{:= S_3} V^T = \begin{bmatrix}
    \Qs{11} & \Qs{12} \\
    0 & \Qs{22} \\
    V_{31}^T & V_{32}^T
    \end{bmatrix}.
\end{equation*}

Perform the intersection of row spaces on $\Qs{11}$ and $V_{31}^T$ to obtain nonsingular matrices $S_4$ and $S_5$ such that

\begin{equation*}
    S_4V_{31}^T = \begin{bmatrix} Y_1^T \\ Y_2^T \\ 0 \end{bmatrix}, \quad S_5 \Qs{11} = \begin{bmatrix}
    Y_1^T \\ Y_3^T \\ 0
    \end{bmatrix}.
\end{equation*}

Then, there exists a permutation $\pi$ such that

\begin{equation*}
    \pi\left(S_4 V_{31}^T - \underbrace{\begin{bmatrix} I && \\ &0& \\ &&0 \end{bmatrix}S_5}_{:=S_6}\Qs{11}\right) = \pi \begin{bmatrix}
    0 \\ Y_2^T \\ 0
    \end{bmatrix} = \begin{bmatrix}
    Y_2^T \\ 0
    \end{bmatrix}.
\end{equation*}

Thus

\begin{equation*}
    \underbrace{\begin{bmatrix}
    I&& \\ &I& \\ - \pi S_6 & 0 & \pi S_4
    \end{bmatrix}S_3}_{:=S_7}V^T = \begin{bmatrix}
    \Qs{11} & \Qs{12} \\
    0 & \Qs{22} \\
    Y_2^T & Y_4^T \\
    0 & \Qs{42}
    \end{bmatrix}.
\end{equation*}

We now know that $\twobyone{\Qs{11}}{Y_2^T}$ is a linearly independent collection of rows which span $B$. Thus, $\twobyone{\Qs{11}}{Y_2^T}$ and $\twobyone{\Qs{11}}{\Qs{31}}$ are both bases for the same matrix so there exists a unique nonsingular matrix $T$ such that

\begin{equation*}
    T\twobyone{\Qs{11}}{Y_2^T} = \twobyone{\Qs{11}}{\Qs{31}}, \mbox{ where } T = \twobyone{\Qs{11}}{\Qs{31}}\left(\twobyone{\Qs{11}}{Y_2^T}\right)^R.
\end{equation*}

It is easy to see that $T$ must have the block lower triangular structure

\begin{equation*}
    T = \twobytwo{I}{}{T_{21}}{T_{22}}, \quad T_{22} \mbox{ nonsingular}.
\end{equation*}

Thus

\begin{equation*}
    \underbrace{\begin{bmatrix}
    I &&& \\
    &I&& \\
    T_{21} & 0 & T_{22} & 0 \\
    &&&I
    \end{bmatrix}S_7}_{:= K}V^T = \begin{bmatrix}
    \Qs{11} & \Qs{12} \\
    0 & \Qs{22} \\
    \Qs{31} & \Qs{32} \\
    0 & \Qs{42}
    \end{bmatrix}.
\end{equation*}

\end{proof}

Note that this proof not only guarentees the existence of this basis, but produces a nonsingular matrix $K$ which converts a given basis $V^T$ to a basis of the desired form.

With that complete, we move on to the construction itself. Suppose we have the complete solution set to a low rank completion problem

\begin{equation*}
    \begin{bmatrix}
    A & B \\
    ? & C
    \end{bmatrix}
\end{equation*}

in terms of the six complementary subspace $\PAbar $, $\PAB $, $\PBbar $, $\QBbar $, $\QBC $, and $\QCbar $. We seek a solution class $\mathscr{S}_{A,\BG,\CH}^Z$ to the low rank completion problem

\begin{equation*}
    \begin{bmatrix}
    A & B & G \\
    ? & C & H
    \end{bmatrix}
\end{equation*}

belonging to the class $\mathscr{S}_{A,B,C}^Y = \{ \FAbar \QAbarA + \PCBC \RBCAB \QABA + \PCCbar Y : \FAbar \mbox{ free} \}$.

\paragraph{1.  Determination of $\QGH$.} Begin by using the intersection of row spaces algorithm on $\onebytwo{B}{G}$ and $\onebytwo{C}{H}$ to find nonsingular matrices $S_1$ and $S_2$ such that

\begin{equation*}
    S_1\onebytwo{B}{G} = \begin{bmatrix} \QGHsym{\intersect{B}{C}} & \QGHsym{\intersect{G}{H}} \\ \QGbarsym{B} & \QGbarsym{G} \\  0 & 0 \end{bmatrix}, \quad S_2\onebytwo{C}{H} = \begin{bmatrix}
    \QGHsym{\intersect{B}{C}} & \QGHsym{\intersect{G}{H}} \\ \QHbarsym{C} & \QHbarsym{H} \\ 0 & 0 \end{bmatrix},
\end{equation*}

where $\onebytwo{\QGHsym{\intersect{B}{C}}}{\QGHsym{\intersect{G}{H}}}$ is a basis for the intersection of the row spaces of $\onebytwo{B}{G}$ and $\onebytwo{C}{H}$. Apply a nonsingular transformation $S_3$ so that

\begin{equation*}
    \QGH := S_3\onebytwo{\QGHsym{\intersect{B}{C}}}{\QGHsym{\intersect{G}{H}}} = \twobytwo{\QBCnum{1}}{\QGHnum{1}}{0}{\QGHnum{2}}, \quad \QBCnum{1} \mbox{ full rank}.
\end{equation*}

\paragraph{2.  Determination of $\QGbar$.} Use the intersection of row spaces algorithm on $\QBC $ and $\QBCnum{1}$ to find a nonsingular matrix $S_4$ such that

\begin{equation*}
    S_4\QBC = \twobyone{\Qs{\ast}}{\QBCnum{2}}.
\end{equation*}

where $R(\QBCnum{2}) \cap R(\QBCnum{1}) = \{0\}$. Since $R(\QBCnum{1}) \subseteq R(\QBC)$, $\Qs{\ast}$ and $\QBCnum{1}$ must be related by a nonsingular transformation $\QBCnum{1} = (\QBCnum{1}(\Qs{\ast})^R)\Qs{\ast}$. Then defining the nonsingular matrix $S_5 := \twobytwo{\QBCnum{1}(\Qs{\ast})^R}{}{}{I}S_4$ such that 

\begin{equation*}
    S_5\QBC  = \twobyone{\QBCnum{1}}{\QBCnum{2}}.
\end{equation*}

We now use Proposition \ref{prop:extend_to_row_basis} with the basis $V^T = \begin{bmatrix} \QGHsym{\intersect{B}{C}} & \QGHsym{\intersect{G}{H}} \\ \QGbarsym{B} & \QGbarsym{G}\end{bmatrix}$, $\Qs{11} = \QBCnum{1}$, $\Qs{12} = \QGHnum{1}$, $\Qs{22} = \QGHnum{2}$, $\Qs{31} = \twobyone{\QBCnum{2}}{\QBbar}$ to find a matrix $S_6$ such that

\begin{equation*}
   S_{6}\begin{bmatrix} \QGHsym{\intersect{B}{C}} & \QGHsym{\intersect{G}{H}} \\ \QGbarsym{B} & \QGbarsym{G}\end{bmatrix} = \left[\arraycolsep=1.4pt\def\arraystretch{1.4}\begin{array}{cc}
   \QBCnum{1} & \QGHnum{1}  \\
   0 & \QGHnum{2} \\\hline
   \QBCnum{2} & \QGbarnum{1} \\
   \QBbar  & \QGbarnum{2}  \\ 
   0 & \QGbarnum{3} 
   \end{array}\right]= \arraycolsep=1.4pt\def\arraystretch{1.4}\begin{bmatrix}
   \QGH \\ \QGbar
   \end{bmatrix}.
\end{equation*}

Thus

\begin{equation*}
   \underbrace{\begin{bmatrix} S_6 & 0 \\0 & I\end{bmatrix} S_1}_{:= S_7}\onebytwo{B}{G} = \left[\arraycolsep=1.4pt\def\arraystretch{1.4}\begin{array}{cc}
   \QBCnum{1} & \QGHnum{1}  \\
   0 & \QGHnum{2} \\\hline
   \QBCnum{2} & \QGbarnum{1} \\
   \QBbar  & \QGbarnum{2}  \\ 
   0 & \QGbarnum{3} \\\hline
   0 & 0
   \end{array}\right]= \arraycolsep=1.4pt\def\arraystretch{1.4}\begin{bmatrix}
   \QGH  \\\QGbar  \\  0
   \end{bmatrix}.
\end{equation*}

\paragraph{3.  Determination of $\QHbar$.} Now use Proposition \ref{prop:extend_to_row_basis} again with the basis $V^T = \begin{bmatrix} \QGHsym{\intersect{B}{C}} & \QGHsym{\intersect{G}{H}} \\ \QHbarsym{C} & \QHbarsym{H}\end{bmatrix}$, $\Qs{11} = \QBCnum{1}$, $\Qs{12} = \QGHnum{1}$, $\Qs{22} = \QGHnum{2}$, $\Qs{31} = \twobyone{\QBCnum{2}}{\QCbar}$ to find a matrix $S_8$ such that

\begin{equation*}
   S_{8}\begin{bmatrix} \QGHsym{\intersect{B}{C}} & \QGHsym{\intersect{G}{H}} \\ \QHbarsym{C} & \QHbarsym{H}\end{bmatrix} = \arraycolsep=1.4pt\def\arraystretch{1.4}\begin{bmatrix}
   \QBCnum{1} & \QGHnum{1}  \\
   0 & \QGHnum{2}  \\\hline
   \QBCnum{2} & \QHbarnum{1} \\
   \QCbar  & \QHbarnum{2}  \\ 
   0 & \QHbarnum{3}
   \end{bmatrix} = \begin{bmatrix}
   \QGH \\ \QHbar
   \end{bmatrix}.
\end{equation*}

Thus

\begin{equation*}
  \underbrace{\begin{bmatrix} S_8 & 0 \\0 & I\end{bmatrix} S_2}_{:= S_9}\onebytwo{C}{H} = \arraycolsep=1.4pt\def\arraystretch{1.4}\begin{bmatrix}
   \QBCnum{1} & \QGHnum{1}  \\
   0 & \QGHnum{2}  \\\hline
   \QBCnum{2} & \QHbarnum{1} \\
   \QCbar  & \QHbarnum{2}  \\ 
   0 & \QHbarnum{3}  \\\hline
   0 & 0
   \end{bmatrix} = \begin{bmatrix}
   \QGH \\ \QHbar \\ 0
   \end{bmatrix}.
\end{equation*}

\paragraph{4.  Determination of $\PGGbar$, $\PGGH$, $\PHHbar$, and $\PHGH$.} Since $\begin{bmatrix}
\QGbar  \\ \QGH 
\end{bmatrix}$ is a row basis for $\begin{bmatrix} B & G\end{bmatrix}$ there exists a corresponding column basis 

\begin{equation*}
    \begin{bmatrix}
\PGGbar & \PGGH
\end{bmatrix} = \begin{bmatrix}
\PBBC^2 & \PBBbar' & \PGGbarthree & \PBBC^1 & \PGGHtwo
\end{bmatrix}
\end{equation*}

such that $\begin{bmatrix} B & G\end{bmatrix} = \begin{bmatrix} \PGGbar & \PGGH \end{bmatrix}\begin{bmatrix}
\QGbar  \\ \QGH 
\end{bmatrix}$. ($\onebytwo{\PGGH}{\PGGbar}$ is simply the first columns of $S_{7}^{-1}$.) We necessarily have $\begin{bmatrix}
\PBBC^1 & \PBBC^2
\end{bmatrix} = \PBBC S_5^{-1}$ and $\PBBbar'=\PBBbar$.

Similarly, since $\begin{bmatrix}
\QHbar  \\ \QGH 
\end{bmatrix}$ is a row basis for $\begin{bmatrix} C & H\end{bmatrix}$ there exists a column basis 

\begin{equation*}
\begin{bmatrix}
\PHHbar & \PHGH
\end{bmatrix} = \begin{bmatrix}
\PCBC^2 & \PCCbar & \PHHbarthree & \PCBC^1 & \PHGHtwo
\end{bmatrix}
\end{equation*}

such that $\begin{bmatrix} C & H\end{bmatrix} = \begin{bmatrix} P_{\bar{H},H} & P_{GH,H} \end{bmatrix}\begin{bmatrix}
\QHbar  \\ \QGH 
\end{bmatrix}$, where $\begin{bmatrix} P_{BC1,C} & P_{BC2,C}
\end{bmatrix} = P_{BC,C}S_5^{-1}$.

\paragraph{4.  Determination of $\PAbarbullet$, $\PAG$, and $\PGbar$.} Now we focus our attention on $\begin{bmatrix} A & B & G\end{bmatrix}$. By our construction the range spaces of $G$ and $B$ have trivial intersection. Thus, using the intersection of column spaces algorithm, there exists a nonsingular matrix $T_1$ such that

\begin{equation*}
    GT_1 = \begin{bmatrix}
    \PGcapA  & \PGcapAbar  & 0
    \end{bmatrix}.
\end{equation*}

Next use the intersection of column spaces algorithm with $A$ and $\onebytwo{\PAB }{\PGcapA }$ to obtain a nonsingular matrices $T_2$ such that

\begin{equation*}
    AT_2 = \begin{bmatrix} \ast & \PAbarbullet  & 0 \end{bmatrix}.
\end{equation*}

Then there must exist a nonsingular matrix $T_3 = \twobytwo{L_{11}}{}{L_{21}}{I}$ such that 

\begin{equation*}
    \onebytwo{\PAbar }{\PAB }T_3 = \begin{bmatrix}
    \begin{bmatrix} \PAbarbullet  & \PGcapA  \end{bmatrix} & \PAB 
    \end{bmatrix}.
\end{equation*}

(Choose $\twobyone{L_{11}}{L_{21}} = \onebytwo{\PAbar}{\PAB}^L\onebytwo{\PAbarbullet}{\PGcapA}$.) Define $\PAG  := \begin{bmatrix}
\PAB  & \PGcapA 
\end{bmatrix}$ and $\PGbar  := \begin{bmatrix}
\PBbar  & \PGcapAbar 
\end{bmatrix}$.

\paragraph{5.  Determination of $\QAbarbulletA$, $\QAGA$, $\QGbarG$, $\QAGG$.} Then

\begin{align*}
    \begin{bmatrix}
    A & B & G
    \end{bmatrix} &= \begin{bmatrix}
    \PAbarbullet  & \PAB  & \PGcapA  & \PBbar  & \PGcapAbar 
    \end{bmatrix}\begin{bmatrix}
    \QAbarbulletA & 0 & 0 \\
    \QABbulletA & \QABB & 0 \\
    \QGcapAA & 0 & \QGcapAG \\
    0 & \QBbarB  & 0 \\
    0 & 0 & \QGcapAbarG
    \end{bmatrix} \\
    &= \begin{bmatrix}
    \PAbarbullet  & \PAG  & \PGbar 
    \end{bmatrix}\begin{bmatrix}
    \QAbarbulletA  & 0 \\
    \QAGA  & \QAGG  \\
    0 & \QGbarG  
    \end{bmatrix},
\end{align*}

where 

\begin{equation*}
    T_3^{-1}\twobyone{\QAbarA }{\QABA } = \begin{bmatrix}
    \QAbarbulletA  \\ \QGcapAA  \\ \QABbulletA 
    \end{bmatrix}, \QABbulletA = \QABA + \underbrace{(-L_{21}L_{11}^{-1})}_{:=J}\QAbarA \quad T_1^{-1} = \begin{bmatrix}
    \QGcapAG  \\ \QGcapAbarG  \\ \ast
    \end{bmatrix}.
\end{equation*}

\paragraph{6.  Studying the relation between the ``$R$'' matrices for the new and old problem.}  We have

\begin{align*}
    \begin{bmatrix}
    B & G
    \end{bmatrix} &= \begin{bmatrix}
    \PGGbar & \PGGH
    \end{bmatrix} \twobytwo{\RGbarAG}{\RGbarGbar}{\RGHAG}{\RGHGbar} \begin{bmatrix}
    \QAGG  \\ \QGbarG  
    \end{bmatrix} \\
    &= \resizebox{0.7\textwidth}{!}{$\begin{bmatrix}
    \PBBCtwo & \PBBbar & \PGGbarthree & \PBBCone & \PGGHtwo
    \end{bmatrix}  \twobytwo{\begin{bmatrix}
    \RGbarAG^{11} & \RGbarAG^{12} \\
    \RGbarAG^{21} & \RGbarAG^{22} \\
    \RGbarAG^{31} & \RGbarAG^{32}
    \end{bmatrix}}{\begin{bmatrix}
    \RGbarGbar^{11} & \RGbarGbar^{12} \\
    \RGbarGbar^{21} & \RGbarGbar^{22} \\
    \RGbarGbar^{31} & \RGbarGbar^{32}
    \end{bmatrix}}{\begin{bmatrix}
    \RGHAG^{11} & \RGHAG^{12} \\
    \RGHAG^{21} & \RGHAG^{22}
    \end{bmatrix}}{\begin{bmatrix}
    \RGHGbar^{11} & \RGHGbar^{12} \\
    \RGHGbar^{21} & \RGHGbar^{22}
    \end{bmatrix}} \begin{bmatrix}
    \QABB & 0 \\
    0 & \QGcapAG  \\
    \QBbarB  & 0 \\
    0 & \QGcapAbarG 
    \end{bmatrix}$}
\end{align*}

Then

\begin{equation*}
    B = \begin{bmatrix}
    \PBBCtwo  & \PBBbar  & \PGGbarthree  & \PBBCone  & \PGGHtwo 
    \end{bmatrix}\begin{bmatrix}
    \RGbarAG^{11} & \RGbarGbar^{11}  \\
    \RGbarAG^{21} & \RGbarGbar^{21}  \\
    \RGbarAG^{31} & \RGbarGbar^{31}  \\
    \RGHAG^{11} & \RGHGbar^{11}   \\
    \RGHAG^{21} & \RGHGbar^{21}  
    \end{bmatrix}\begin{bmatrix}
    \QABB\\
    \QBbarB 
    \end{bmatrix}
\end{equation*}

We conclude that $\RGbarAG^{31} = 0$, $\RGbarGbar^{31} = 0$, $\RGHAG^{21}=0$, and $\RGHGbar^{21}=0$ so

\begin{align*}
    B &= \begin{bmatrix}
    \PBBCtwo  & \PBBbar  & \PBBCone 
    \end{bmatrix}\begin{bmatrix}
    \RGbarAG^{11} & \RGbarGbar^{11}  \\
    \RGbarAG^{21} & \RGbarGbar^{21}  \\
    \RGHAG^{11} & \RGHGbar^{11}  
    \end{bmatrix}\begin{bmatrix}
    \QABB  \\ \QBbarB 
    \end{bmatrix} \\
    &= \begin{bmatrix}
    \PBBbar  & \PBBC 
    \end{bmatrix}\begin{bmatrix}
    \RGbarAG^{21} & \RGbarGbar^{21} \\
    S_5^{-1}\begin{bmatrix} \RGHAG^{11} \\ \RGbarAG^{11}\end{bmatrix} & S_5^{-1}\begin{bmatrix} \RGHGbar^{11} \\ \RGbarGbar^{11}\end{bmatrix}
    \end{bmatrix}\begin{bmatrix}
    \QABB  \\ \QBbarB 
    \end{bmatrix} \\
    &= \begin{bmatrix}
    \PBBbar  & \PBBC 
    \end{bmatrix}\begin{bmatrix}
    \RBbarAB & \RBbarBbar \\
    \RBCAB & \RBCBbar
    \end{bmatrix}\begin{bmatrix}
    \QABB  \\ \QBbarB 
    \end{bmatrix},
\end{align*}

so we conclude $\RBbarAB = \RGbarAG^{21}$, $\RBbarBbar = \RGbarGbar^{21}$, $\RBCAB = S_5^{-1}\begin{bmatrix} \RGHAG^{11} \\ \RGbarAG^{11}\end{bmatrix} $, and $\RBCBbar = S_5^{-1}\begin{bmatrix} \RGHGbar^{11} \\ \RGbarGbar^{11}\end{bmatrix}$.

\paragraph{7.  Comparison of solution sets.} The solution to the original problem $\twobytwo{A}{B}{?}{C}$ is $\mathscr{S}_{A,B,C} = \{ X^{A,B,C}_{\FAbar,\FCbar} : \FAbar, \FCbar \mbox{ free}\}$ where

\begin{equation*}
    X^{A,B,C}_{\FAbar,\FCbar} = \FAbar \QAbarA  + \PCBC \RBCAB \QABA  + \PCCbar\FCbar .
\end{equation*}

The solution set the new problem $\begin{bmatrix} A & B & G \\ ? & C & H \end{bmatrix}$ is $\mathscr{S}_{A,\BG,\CH} = \{ X^{A,\BG,\CH}_{\FAbarbullet,\FHbar} : \FAbarbullet, \FHbar \mbox{ free} \}$ to the new problem is 

\begin{align*}
    X^{A,\BG,\CH}_{\FAbarbullet,\FHbar} &= \FAbarbullet \QAbarbulletA  + \PHGH \begin{bmatrix}
    \RGHAG^{11} & \RGHAG^{12} \\
    \RGHAG^{21} & \RGHAG^{22}
    \end{bmatrix}\QAGA  + \PHHbar \FHbar \\
    &= \FAbarbullet\QAbarbulletA  + \onebytwo{\PCBCone}{\PHGHtwo}\begin{bmatrix}
    \RGHAG^{11} & \RGHAG^{12} \\
    0 & \RGHAG^{22}
    \end{bmatrix}\twobyone{\QABbulletA }{\QGcapAA } \\
    &\qquad + \begin{bmatrix}
    \PCBCtwo & \PCCbar & \PHHbarthree
    \end{bmatrix}\begin{bmatrix}
    \FHbarnum{1} \\
    \FHbarnum{2} \\ 
    \FHbarnum{3}
    \end{bmatrix} \\
    &= \FAbar \QAbarbulletA + \PCBCone \RGHAG^{11}(\QABA + J\QAbarA) + (\PCBCone \RGHAG^{12} + \PHGHtwo\RGHAG^{22}) \QGcapAA \\
    &\qquad + \begin{bmatrix}
    \PCBCtwo & \PCCbar & \PHHbarthree
    \end{bmatrix}\begin{bmatrix}
    \FHbarnum{1} \\
    \FHbarnum{2} \\ 
    \FHbarnum{3}
    \end{bmatrix} \\
\end{align*}

Denote $V := \PCBCone\RGHAG^{11}J\QAbarA + (\PCBCone \RGHAG^{12} + \PHGHtwo\RGHAG^{22}) \QGcapAA$ which lies in the row space of $\QAbarA$. Make the change of variables

\begin{equation*}
    \begin{bmatrix}
    \FHbarnum{1} \\
    \FHbarnum{2} \\ 
    \FHbarnum{3}
    \end{bmatrix} = \begin{bmatrix}
    \FHbarnum{1\bullet} \\
    \FHbarnum{2\bullet} \\ 
    \FHbarnum{3\bullet}
    \end{bmatrix} + \begin{bmatrix}
    \RGbarAG^{11}\QABA \\ 0 \\ 0
    \end{bmatrix}.
\end{equation*}

Then noting

\begin{equation*}
    \onebytwo{\PBBCtwo}{\PBBCone}\underbrace{\twobyone{\RGbarAG^{11}}{\RGHAG^{11}}}_{=S_5\RBCAB}\QABA = \onebytwo{\PBBCtwo}{\PBBCone}S_5 \RBCAB \QABA = \PBBC\RBCAB\QABA. 
\end{equation*}

Thus,

\begin{equation*}
    X^{A,\BG,\CH}_{\FAbarbullet,\FHbar} = \FAbarbullet \QAbarbulletA  + \PCBC \RBCAB \QABA + \begin{bmatrix}
    \PCBCtwo & \PCCbar & \PHHbarthree
    \end{bmatrix}\begin{bmatrix}
    \FHbarnum{1\bullet} \\
    \FHbarnum{2\bullet} \\ 
    \FHbarnum{3\bullet}
    \end{bmatrix}
\end{equation*}

\paragraph{8.  Determination of $Z$.} From here, it is easy to see that if we set

\begin{equation*}
    Z := \begin{bmatrix}
    \FHbarnum{1} \\
    \FHbarnum{2} \\ 
    \FHbarnum{3}
    \end{bmatrix} = \begin{bmatrix}
    \RGbarAG^{11}\QABA \\ Y \\ 0
    \end{bmatrix},
\end{equation*}

then indeed we have $\mathscr{S}_{A,\BG,\CH}^Z \subseteq \mathscr{S}_{A,B,C}^Y$.






\section{Learning $\graph$-semi-separable representations}\label{sec:learning_gss}

Determination of minimal $\graph$-SS representations is a challenging problem. Even the simplest nontrivial $\graph$-SS representation, the CSS representation, requires the solution of a difficult low-rank completion problem (Section~\ref{sec:css}). Thus, as an alternate approach, we consider the problem of determining a $\graph$-SS representation by using a black box optimization routine. 

The strategy is thus: make some guesses about the interaction ranks $r_i^{g}$ and $r_i^{h}$, compute some pairs $(x_i,Ax_i)$ and $(x_i,A^Tx_i)$, and then use a black box optimization routine to minimize the squared error by choosing the entries of the generator matrices:
\begin{equation}
    \mathscr{L} = \frac{1}{m} \sum_{i=1}^m \left( \left\| \hat{A}x_i - Ax_i \right\|^2 + \left\| \hat{A}^Tx_i - A^Tx_i \right\|^2 \right),
\end{equation}
where $\hat{A}$ is the putative $\graph$-SS representation of $A$. Of course, $\mathscr{L}$ and $\hat{A}$ are parametric functions of the generators $(\diagmat{U}, \diagmat{W}, \diagmat{V}, \diagmat{D}, \diagmat{P}, \diagmat{R}, \diagmat{Q})$.

In order to train the $\graph$-SS representation, we need a fast multiplication algorithm for $A^T$. We claim that the following recurrences constitute a fast multiplication algorithm:
\begin{subequations}\label{eq:gss_transpose_recurrences}
\begin{align}
    g_j &= \sum_{i \in \dirpath_{\upstream}(j)} U_{i,j}^Tx_i + \sum_{i \in \dirpath_{\upstream}(j)} W_{i,j}^T g_i, \label{eq:transpose_g_form1}\\
    h_j &= \sum_{i \in \dirpath_{\downstream}(j)} P_{i,j}^Tx_i + \sum_{i\in\dirpath_\downstream(j)} R_{i,j}^T h_i, \\
    b_j &= D_j^Tx_j + V_jg_j + Q_jh_j.
\end{align}
\end{subequations}
\begin{proposition}
The recurrences \eqref{eq:gss_transpose_recurrences} compute the product $b = A^Tx$.
\end{proposition}

\begin{proof}

By Equation~\eqref{eq:GSSexplicit}, we have that the block entries of $A^T$ are given by
\begin{equation*}
    [A^T]_{j i} = [A]^T_{ij}
    = \begin{cases}
         \sum_{s \in \dirpath_\upstream(i)} Q_j \Phi_{\upstream}(s,j)^T P_{i,s}^T & i \succ j, \\
         D_i^T & i = j, \\
         \sum_{s \in \dirpath_{\downstream}(i)} V_{j} \Phi_{\downstream}(s,j)^TU_{i,s}^T & i \prec j.
      \end{cases}
\end{equation*}
Then to compute the $j$th block component of $b = A^Tx$, we have
\begin{equation*}
    b_j = \sum_{i=1}^n A_{j i}^Tx_i = Q_j \sum_{i \succ j} \sum_{s \in \dirpath_\upstream(i)}  \Phi_{\upstream}(s,j)^T P_{i,s}^T x_i + D_j^Tx_{j} + V_{j}\sum_{i \prec j} \sum_{s \in \dirpath_{\downstream}(i)}  \Phi_{\downstream}(s,j)^TU_{i,s}^Tx_i.
\end{equation*}
Thus, it is sufficient to show that
\begin{align}
    h_j &= \sum_{i \succ j} \sum_{s \in \dirpath_\upstream(i)}  \Phi_{\upstream}(s,j)^T P_{i,s}^T x_i, \nonumber\\
    g_j &= \sum_{i \prec j} \sum_{s \in \dirpath_{\downstream}(i)}  \Phi_{\downstream}(s,j)^TU_{i,s}^Tx_i. \label{eq:tranpose_g_form2}
\end{align}
We shall prove the latter by induction on $j$. (The proof of the former is similar.) Choose an indexing such that our Hamiltonian path is $1 \to 2 \to \cdots \to n$. As a base case $g_1 = 0$ using both formulas \eqref{eq:transpose_g_form1} and \eqref{eq:tranpose_g_form2}.

Now suppose that the formulas agree for $1,\ldots,j-1$. Then we have
\begin{equation*}
    g_j = \sum_{i \in \dirpath_{\upstream}(j)} \underbrace{\Phi_\downstream(j,j)^T}_{=I}U_{i,j}^Tx_i + \sum_{i \prec j} \sum_{\stackrel{s \in \dirpath_{\downstream}(i)}{s \ne j}} \Phi_\downstream(s,j)^T U_{i,s}^Tx_i.
\end{equation*}
Under the assumption that $s \ne j$, $\Phi_\downstream(s,j)$ is zero unless $s \prec j$. Note that the recurrence \eqref{eq:transition} can be written in alternate form
\begin{equation*}
    \Phi_{\downstream}(s,j) = \sum_{i \in \dirpath_{\upstream}(j)} = \Phi_{\downstream}(s,k)W_{k,j} \mbox{ for } s \prec j.
\end{equation*}
Thus, we may expand using this identity to obtain
\begin{align*}
    g_j &= \sum_{i \in \dirpath_{\upstream}(j)} U_{i,j}^Tx_i + \sum_{i \prec j} \sum_{\stackrel{s \in \dirpath_{\downstream}(i)}{s \prec j}} \Phi_\downstream(s,j)^T U_{i,s}^Tx_i \\
        &= \sum_{i \in \dirpath_{\upstream}(j)} U_{i,j}^Tx_i + \sum_{i \prec j} \sum_{\stackrel{s \in \dirpath_{\downstream}(i)}{s \prec j}} \sum_{k\in \dirpath_\upstream(j)} W_{k,j}^T \Phi_\downstream(s,k)^T U_{i,s}^Tx_i \\
        &= \sum_{i \in \dirpath_{\upstream}(j)} U_{i,j}^Tx_i + \sum_{k\in \dirpath_\upstream(j)} W_{k,j}^T \sum_{i \prec j} \sum_{\stackrel{s \in \dirpath_{\downstream}(i)}{s \prec j}} \Phi_\downstream(s,k)^T U_{i,s}^Tx_i.
\end{align*}
Now suppose that $i,s \prec j$, $s \in \dirpath_{\downstream}(i)$, but $i \not\prec k$. Then $s \not\prec k$ since if $s \prec k$, then $i \prec s \prec j$ a contradiction. Thus, in this case we have $\Phi_{\downstream}(s,k) = 0$. Note that we may drop the condition $s \prec j$ from the summation because $\Phi_{\downstream}(s,k)$ will be zero if $s \not\prec j$. Thus, using the induction hypothesis, we obtain
\begin{equation*}
    g_j = \sum_{i \in \dirpath_{\upstream}(j)} U_{i,j}^Tx_i + \sum_{k\in \dirpath_\upstream(j)} W_{k,j}^T \sum_{i \prec k} \sum_{s\in\dirpath_{\downstream}(i)} \Phi_\downstream(s,k)^T U_{i,s}^Tx_i = \sum_{ \in \dirpath_{\upstream}(j)} U_{i,j}^Tx_i + \sum_{k\in \dirpath_\upstream(j)} W_{k,j}^T g_k,
\end{equation*}
which was as to be shown.


\end{proof}

}
\section{Properties of CSS representations}\label{app:css_properties}

Just as with SSS and the line graph, there is a close connection between CSS representations and GIRS matrices with the cycle graph. \ethancheck  From the general analysis on $\graph$-SS representations in Section~\ref{sec:edv}, we know that a CSS representation with dimensions $r_i^g,r_i^h \le r$ for every $1\le i \le n-1$ is a GIRS-$4r$ matrix. This result can however be sharpened specifically for  CSS matrices.
\begin{proposition} \label{prop:CSSimpliesGIRS}
Let $A=\CSS(\diagmat{U}, \diagmat{W}, \diagmat{V}, \diagmat{D}, \diagmat{P},  \diagmat{R}, \diagmat{Q})$ be a CSS representation. Then $(A,\graph)$ is a GIRS-$r$ matrix with 
\begin{equation} 
r =  \max_{i=2,\ldots,n-2} \{r_i^h + r_1^h,r_i^g + r_{n-1}^g\}. \label{eq:whatever}
\end{equation} 
\end{proposition}
\begin{proof}
The proof of this result follows an analogous path as in the proof of Theorem~\ref{thm:SSStoGIRS}, but $A_{u(\complementgraph, \subgraph), \subgraph}$ and $A_{d(\complementgraph, \subgraph), \subgraph}$ are now factorized as
\begin{displaymath}
A_{u(\complementgraph, \subgraph), \subgraph}  =  \begin{bmatrix}
	P_{t+1} & 0  \\ 
	P_{t+2}R_{t+1} & 0 \\
	P_{t+3} R_{t+2} R_{t+1} & 0 \\
	\vdots & \vdots \\
	P_{n} R_{n-1} R_{n-2} \cdots R_{t+1} & P_0
	\end{bmatrix}  \begin{bmatrix} R_{t} \cdots R_{s+2} R_{s+1} Q_{s}^T  & \cdots &  R_{t} R_{t-1}  Q_{t-2}^T & R_{t} Q_{t-1}^T & Q_t^T \\
	Q^T_1  & \cdots &  0 &  0 & 0
	\end{bmatrix} 
\end{displaymath}
if $s=1$, and
\begin{displaymath}
A_{d(\complementgraph, \subgraph), \subgraph}  =   \begin{bmatrix}
	U_{1} W_{2} W_{3} \cdots W_{s-1} & U_0 \\ 
		\vdots & \vdots \\
	U_{s-3} W_{s-2} W_{s-1} & 0 \\
	U_{s-2} W_{s-1} & 0  \\
	U_{s-1}  & 0
	\end{bmatrix}  \begin{bmatrix} V_s^T & W_{s} V_{s+1}^T    & W_{s} W_{s+1}  V_{s+2}^T  &  \cdots  & W_{s}  W_{s+1} \cdots W_{t-1} V_{t}^T  \\
	0 & 0 & 0 & \cdots & 0 \end{bmatrix}
\end{displaymath}
if $t=n$, respectively.
\end{proof}
At this point, the converse result showing that GIRS implies $\graph$-SS is not known generally, but we do know that for at least CSS representations this proposition still passes through.
\begin{proposition} \label{prop:GIRSimpliesCSS}
If $(A,\graph)$ is a GIRS-$c$ matrix associated with the cycle graph, then $(A,\graph)$ admits a CSS representation $A=\CSS(\diagmat{U}, \diagmat{W}, \diagmat{V}, \diagmat{D}, \diagmat{P},  \diagmat{R}, \diagmat{Q})$ with dimensions $r_i^g,r_i^h \le 2c$ for every $1\le i \le n-1$.\ethancheck
\end{proposition}
\begin{proof}
    The proof of this proposition is through the GIRS properties of SSS matrices (Theorem \ref{thm:SSStoGIRS}) and the reduction of CSS matrices to SSS matrices.
\end{proof}
In the case of SSS, the size of the representation is exactly preserved under taking inverses (see Proposition \ref{prop:SSSalgebra}). As the following example shows, this need not be the case for CSS matrices. \ethancheck 
\begin{example}\label{ex:css_inversion}
Consider the matrix \eqref{eq:circular_tridiagonal}. By solving the low rank completion problem, we are able to get $\rg_{A,1} = \rh_{A,2} = 2 > \rg_{A,2} = \rh_{A,1} = 1$. However, for the inverse, 
\begin{displaymath}
A^{-1} =\resizebox{0.9\textwidth}{!}{$\frac{1}{\det A} \left(
\begin{array}{c|c|c }
\begin{array}{c c} 3 a^4 b-4 a^2 b^3+b^5   &  -a^5-a^4 b+3 a^3 b^2-a b^4 \\
-a^5-a^4 b+3 a^3 b^2-a b^4 & a^4 b-a^2 b^3 \end{array} & \begin{array}{c c} -2 a^4 b+a^3 b^2+a^2 b^3 & a^5+a^4 b-a^3 b^2-a^2 b^3 \\
-a^4 b+a^3 b^2 & -a^4 b+a^3 b^2 \end{array} & \begin{array}{c c} a^4 b-2 a^3 b^2+a b^4 & -a^5-a^4 b+3a^2 b^3-b^5  \\
a^4 b-a^2 b^3 & a^4 b-2 a^3 b^2+a b^4 \end{array}  \\
\hline
\begin{array}{c c} -2 a^4 b+a^3 b^2+a^2 b^3 & -a^4 b+a^3 b^2 \\
 a^5+a^4 b-a^3 b^2-a^2 b^3 & -a^4 b+a^3 b^2 \end{array} & \begin{array}{c c}  -a^4 b+a^3 b^2 & -a^4 b+a^3 b^2 \\
 2 a^4 b-2 a^2 b^3 & -a^5-a^4 b+2 a^3 b^2 \end{array} & \begin{array}{c c} a^4 b-a^2 b^3 & a^4 b-2 a^3 b^2+a b^4  \\
-a^4 b+a^3 b^2 & a^5+a^4 b-a^3 b^2-a^2 b^3  \end{array}  \\
\hline
\begin{array}{c c}  a^4 b-2 a^3 b^2+a b^4 & a^4 b-a^2 b^3  \\
 -a^5-a^4 b+3 a^2 b^3-b^5 & a^4 b-2 a^3 b^2+a b^4 \end{array} & \begin{array}{c c} -a^4 b+a^3 b^2 & -a^4 b+a^3 b^2 \\
a^5+a^4 b-a^3 b^2-a^2 b^3 & -2 a^4 b+a^3 b^2+a^2 b^3 \end{array}  & \begin{array}{c c} a^4 b-a^2 b^3 & -a^5-a^4 b+3 a^3 b^2-a b^4 \\
-a^5-a^4 b+3 a^3 b^2-a b^4 & 3 a^4 b-4
a^2 b^3+b^5
 \end{array} \end{array} \right)$},
\end{displaymath}
we have $\rg_{A^{-1},1} \rg_{A^{-1},2} = \rh_{A^{-1},1} = \rh_{A^{-1},2} = 2$. The CSS representation of the inverse $A^{-1}$ is slightly larger than the CSS representation of the original matrix $A$. By the same token, the CSS representation of $(A^{-1})^{-1}$ is smaller than the CSS representation of $A^{-1}$. Thus, the minimal CSS representation of the inverse of a matrix may be larger or smaller than that of the original matrix. \ethancheck \nithincheck
\end{example}

The previous example should not make us too pessimistic, as we shall soon show that the size of the CSS representation of an inverse cannot differ too much from the original matrix. Indeed, by Proposition~\ref{prop:CSSimpliesGIRS} and Proposition~\ref{prop:GIRSimpliesCSS}, we already know that the dimensions $\tilde{r}^g_i, \tilde{r}^h_i$ of the inverse CSS  representation are bounded by $2r$, where $r$ is defined by\eqref{eq:whatever}. However, this bound can be misleading particularly in situations where $\rh_1>> \rh_i$ for $i=2,\ldots,n-1$ and $\rg_{n-1} >> r^{g}_{i}$ for $i=1,\ldots,n-2$. The GIRS property here  does not fully shed light on why CSS may be a more efficient representation than SSS in certain scenarios. If the off-diagonal corner blocks are of high rank $R$ (such as in integral equations discretized on the circle), then the CSS representation need only suffer this cost twice in $r_{n-1}^g$ and $r_1^h$. However, for SSS, the high rank corner blocks affect all of matrices in the representation. The reason the GIRS property does not explain this difference is that the GIRS property only allows us to bound the \textit{maximum rank} appearing in the representation, whereas the complexity of computing matrix-vector products and solving linear systems depends on the sum of ranks and sum of squares of the ranks, respectively. \ethancheck 

Nonetheless,  a sharper bound on the dimensions of the inverse of CSS representation can be obtained as follows. For that, we particularly consider the specific scenario where $A$ can be written as the sum of SSS matrix with small GIRS constant $r$ (subsequently referred to as an SSS-$r$ matrix) and a corner block perturbation \eqref{eq:splitting}. 

\begin{theorem}\label{thm:CSSInversion}
Let $r_i^g$ and $r_i^h$ represent the sizes occurring in a minimal CSS representation of a matrix $A$ and $\tilde{r}_i^g$ and $\tilde{r}_i^h$ represent the sizes of the minimal CSS representation of $A^{-1}$. Suppose in addition that $A$ can be written as a sum of an SSS-$r$ matrix $B$ and a corner block perturbation, i.e. $A = B + \CB(Y,X)$, where $\rank X, \rank Y \le R$. Then 
\begin{equation*}
    \tilde{r}_1^h,\tilde{r}_{n-1}^g \le 6r + 4R, \quad \tilde{r}_i^h\le 6r \mbox{ for } 2\le i \le n-1, \quad \tilde{r}_i^g \le 6r  \mbox{ for }  1\le i \le n-2.
\end{equation*}
\end{theorem}
To prove this Theorem~\ref{thm:CSSInversion}, we will be needing the following lemma.
\begin{lemma}\label{lem:SSStimesCornerBlock}
Let $A$ be an SSS-$r$ matrix. Then for matrices $E_1,E_2,E_3,E_4$ of the appropriate sizes, there exists an SSS matrix $B$ with GIRS constant $2r$ and matrices $F_1,F_2$ of the appropriate sizes such that
\begin{equation*}
    A\underbrace{\begin{bmatrix}
    E_1 & 0 & \cdots & 0 & E_2 \\
    0 & 0 & \cdots & 0 & 0 \\
    \vdots & \vdots & \cdots & \vdots \\
    0 & 0 & \cdots & 0 & 0 \\
    E_3 & 0 & \cdots & 0 & E_4
    \end{bmatrix}}_{=E} = B + \CB(F_1,F_2),
\end{equation*}
and $\rank F_1 \le \rank E_2 + \rank E_4$ and $\rank F_2 \le \rank E_1 + \rank E_2$ where $\CB$ is defined as in \eqref{eq:corner_block}. \ethancheck %
\end{lemma}

\begin{proof}
    Write $A$ and $E$ as
    \begin{equation*}
        A = \begin{bmatrix} D_1 & A_{12} & A_{13} \\ A_{21} & D_2 & A_{23} \\ A_{31} & A_{32} & D_3 \end{bmatrix}, \quad E = \begin{bmatrix} E_1 & 0 & E_2 \\ 0 & 0 & 0 \\ E_3 & 0 & E_4\end{bmatrix}.
    \end{equation*}
    Then
    \begin{align*}
        AE &= \begin{bmatrix}
        D_1E_1 + A_{13}E_3 & 0 & D_1E_2 + A_{13}E_4 \\
        A_{21}E_1 + A_{23}E_3 & 0 & A_{21}E_2 + A_{23}E_4 \\
        A_{31}E_1 + D_3E_3 & 0 & A_{31}E_2 + D_3E_4
        \end{bmatrix} \\
        &= \underbrace{\begin{bmatrix}
        D_1E_1 + A_{13}E_3 & 0 & 0 \\
        A_{21}E_1 + A_{23}E_3 & 0 & A_{21}E_2 + A_{23}E_4 \\
        0 & 0 & A_{31}E_2 + D_3E_4
        \end{bmatrix}}_{:=B} + \CB(\underbrace{D_1E_2 + A_{13}E_4}_{:=F_1},\underbrace{A_{31}E_1 + D_3E_3}_{:=F_2})
    \end{align*}
    Since $A_{21}$ and $A_{23}$ both have rank $\le r$, $A_{21}E_1 + A_{23}E_3$ and $A_{21}E_2 + A_{23}E_4$ have rank $\le 2r$ so $B$ is an SSS-$2r$ matrix. The rank bounds for $F_1$ and $F_2$ are obvious. \ethancheck  \nithincheck
    
\end{proof}


%

\begin{proof}[Proof of Theorem~\ref{thm:CSSInversion}]
    
    By the Sherman–Morrison–Woodbury formula,
    \begin{equation*}
        A^{-1} = B^{-1} - B^{-1}\begin{bmatrix}
        0 & I \\
        0 & 0 \\
        I & 0
        \end{bmatrix}\underbrace{\left(I + \begin{bmatrix} X & 0 & 0 \\ 0 & 0 & Y\end{bmatrix}B^{-1}\begin{bmatrix} 0 & I \\ 0 & 0 \\ I & 0 \end{bmatrix} \right)^{-1}}_{:=S}\begin{bmatrix} X & 0 & 0 \\ 0 & 0 & Y \end{bmatrix} B^{-1}.
    \end{equation*}
    Writing $S^{-1} = \twobytwo{a}{b}{c}{d}$, we have
    \begin{equation}\label{eq:CSSInverseSchur}
        A^{-1} = B^{-1} - B^{-1}\begin{bmatrix} cX & 0 & dY \\
        0 & 0 & 0 \\
        aX & 0 & bY
        \end{bmatrix}B^{-1}
    \end{equation}
    Since SSS-$r$ matrices are closed under inversion $B^{-1}$ is an SSS-$r$ matrix. By Lemma \ref{lem:SSStimesCornerBlock},
    \begin{equation*}
        B^{-1}\begin{bmatrix} cX & 0 & dY \\
        0 & 0 & 0 \\
        aX & 0 & bY
        \end{bmatrix} = C + \CB(F_1,F_2),
    \end{equation*}
    where $C$ is an SSS-$2r$ matrix and $\rank F_1,\rank F_2 \le 2R$. Then
    \begin{equation*}
        B^{-1}\begin{bmatrix} cX & 0 & dY \\
        0 & 0 & 0 \\
        aX & 0 & bY
        \end{bmatrix}B^{-1} = CB^{-1}+\CB(F_1,F_2)B^{-1}.
    \end{equation*}
    By Proposition \ref{prop:SSSalgebra}, $CB^{-1}$ is an SSS-$3r$ matrix. By Lemma \ref{lem:SSStimesCornerBlock},
    \begin{equation*}
        \CB(F_1,F_2)B^{-1} = D + \CB(G_1,G_2)
    \end{equation*}
    for $D$ an $SSS$-$2r$ matrix and $\rank G_1,\rank G_2 \le 4R$. Thus, using Proposition \ref{prop:SSSalgebra} one final time,
    \begin{equation*}
        A^{-1} = B^{-1} - B^{-1}\begin{bmatrix} cX & 0 & dY \\
        0 & 0 & 0 \\
        aX & 0 & bY
        \end{bmatrix}B^{-1} = \underbrace{B^{-1} + CB^{-1} + D}_{\mbox{SSS-$6r$}} + \CB(G_1,G_2)
    \end{equation*}
    The stated rank bounds follow from the CSS construction algorithm, Theorem \ref{thm:CSS_construction}.
\end{proof}
\begin{remark}
We believe that the constants $6$ and $4$ in the bounds in Theorem \ref{thm:CSSInversion} are unlikely to be tight. Further refinements are likely possible by using the specific structure of the Schur complement expression \eqref{eq:CSSInverseSchur} rather than using the algebraic properties of the SSS representation (Proposition \ref{prop:SSSalgebra}) out of the box.  Also, note that since inversion is an involution, the possibility that the CSS representation of the inverse may be larger than the original matrix necessarily implies that it may be smaller as well!
\end{remark}

\end{document}